\documentclass[notitlepage,10pt]{article}
\usepackage{amssymb,amsbsy,amsmath,mathrsfs,enumitem,ulem,mathabx}
\usepackage{pgfplots}

\catcode`\@=11
\@addtoreset{equation}{section}

\catcode`\@=12

\newtheorem{Theorem}{Theorem}[section]
\newtheorem{Lemma}[Theorem]{Lemma}
\newtheorem{Proposition}[Theorem]{Proposition}
\newtheorem{Corollary}[Theorem]{Corollary}
\newtheorem{Remark}[Theorem]{Remark}

\def\sech{\mathrm{\hspace{1.5pt}sech\hspace{1.5pt}}} 
\def\R{{\mathbb R}}
\def\E{{\mathcal E}}
\def\F{{\mathcal F}}
\def\G{{\mathcal G}}
\def\L{{\mathcal L}}
\def\M{{\mathcal M}}
\def\N{{\mathcal N}}
\def\Xe{X_{\eps,a}}
\def\Xet{(X_{\eps,a})_t}
\def\Xett{(X_{\eps,a})_{tt}}
\def\Xeth{(X_{\eps,a})_\theta}
\def\Xetth{(X_{\eps,a})_{t\theta}}
\def\Xethth{(X_{\eps,a})_{\theta\theta}}
\def\Ne{N_{\eps,a}}
\def\Net{(N_{\eps,a})_t}
\def\Nett{(N_{\eps,a})_{tt}}
\def\Neth{(N_{\eps,a})_\theta}
\def\Netth{(N_{\eps,a})_{t\theta}}
\def\Nethth{(N_{\eps,a})_{\theta\theta}}
\def\eps{\varepsilon}
\def\normal{N}
\def\bta{\beta}
\def\ak{a_{k}}
\def\ek{\varepsilon_{k}}
\def\phik{\varphi_{k}}
\def\nk{n_{k}}
\def\Proof{\noindent\textit{Proof. }}
\def\qed{$~\square$\goodbreak \medskip}

\title{On the non-existence of compact surfaces of
genus one\\ with prescribed, almost constant mean curvature, 
\\ close to the singular limit}
\author{Paolo Caldiroli\footnote{Dipartimento di Matematica ``Giuseppe Peano'', Universit\`a degli Studi di Torino, via Carlo Alberto, 10 -- 10123 Torino, Italy. Email: \tt{paolo.caldiroli@unito.it}}, Alessandro Iacopetti\footnote{Dipartimento di Matematica ``Federigo Enriques'', Universit\`a degli Studi di Milano, via Saldini, 50 -- 20133 Milano, Italy. Email: \tt{alessandro.iacopetti@unimi.it}}, Monica Musso\footnote{Department of Mathematical Sciences, University of Bath, North Rd, Bath, BA27AY, United Kingdom. Email: \tt{m.musso@bath.ac.uk}}}

\date{}

\begin{document}
\maketitle

\begin{abstract}
\noindent
In Euclidean 3-space endowed with a Cartesian reference system we consider a class of surfaces, called Delaunay tori, constructed by bending segments of Delaunay cylinders with neck-size $a$ and $n$ lobes along circumferences centered at the origin. Such surfaces are complete and compact, have genus one and almost constant, say 1, mean curvature, when $n$ is large. Considering a class of mappings $H\colon\mathbb{R}^{3}\to\mathbb{R}$ such that $H(X)\to 1$ as $|X|\to\infty$ with some decay of inverse-power type, we show that for $n$ large and $|a|$ small, in a suitable neighborhood of any Delaunay torus with $n$ lobes and neck-size $a$ there is no parametric surface constructed as normal graph over the Delaunay torus and whose mean curvature equals $H$ at every point. 
\smallskip

\smallskip

\end{abstract}

\section{Introduction and main results}

In this paper we will study the problem of compact parametric surfaces in Euclidean 3-space with genus one and prescribed, almost constant mean curvature. 

As a starting point, we can refer to the deep and elegant result by A.D. Alexandrov \cite{Alexandrov} who in the late 50's proved that the only compact, connected, embedded surfaces with (non zero) constant mean curvature (shortly, CMC) are Euclidean spheres. In particular, there is no compact, embedded CMC surface with genus one. By the way, the embeddedness assumption is essential, since immersed CMC tori exist, as proved by H. Wente in 1984 \cite{Wente}.

The basic issue we are going to address concerns the rigidity of the Alexandrov result when one prescribes the mean curvature function as a perturbation of a constant. 

Besides an intrinsic interest in Geometry (see problem 59 raised by S.T. Yau in \cite{Yau}), considering prescribed, \textit{non constant} mean curvature functions, possibly close to a constant, is meaningful in capillarity theory, in order to take account of external forces, like pressure differences or gravity (see, e.g., \cite{Finn}).

One can easily construct compact, embedded surfaces with genus one, having almost constant mean curvature. For example, this can be accomplished by taking circular collars made by almost tangent spheres with the same radius attached each other by very small necks shaped on suitably dilated catenoids (see \cite{CiraoloVezzoni}). Equivalently, one can take a piece of unduloid with a very small neck-size and bend it along a circumference. This ``multi-bubble'' construction will be performed later on. In fact, in \cite{CiraoloMaggi} G. Ciraolo and F. Maggi proved that a compact, embedded surface with prescribed almost constant mean curvature has to be close to an arrangement of almost tangent Euclidean spheres, and qualitative estimates in terms of the oscillation of the mean curvature are provided. 

Considering these facts, we can ask if, for a class of functions $H\colon\R^{3}\to\R$ close to 1 and satisfying suitable reasonable conditions, there exist parametric surfaces with genus one, mean curvature equal to $H$ at every point, and shaped on collars of almost tangent round spheres, as described above. 

Let us point out that in the setting of CMC surfaces there are many remarkable results about the existence of complete, compact/non-compact, embedded/non-embedded surfaces with prescribed genus and a given number of ends (in the case of unbounded surfaces). Such surfaces are constructed by gluing segments or halves of Delaunay surfaces, according to techniques introduced by Kapouleas (\cite{Kapouleas1990, Kapouleas1991}) and subsequently adjusted or revised by other authors (e.g., \cite{BreinerKapouleas2014, GrosseBrauckmann1993, GrosseKusnerSullivan2007, MazzeoPacard2001, MazzeoPacardPollack2001}). 

In fact, a similar construction technique, with suitable adjustments, works even in other contexts. For example, in the purely PDE setting of the nonlinear stationary Schr\"odinger equation
\begin{equation}
\label{NLS}
-\Delta u+V(|x|)u=u^{p},\quad u>0\,,\quad u\in H^{1}(\R^{N})
\end{equation}
where $1<p<\frac{N+2}{N-2}$ if $N\ge 3$, $1<p<\infty$ if $N=2$, and $V$ is a radial potential such that
\begin{equation}
\label{V-shape}
V(r)=1+\frac{A}{r^{\bta}}+O\left(\frac{1}{r^{\bta+\nu}}\right),\quad\text{as }r\to\infty
\end{equation}
for some $A>0$, $\bta>1$ and $\nu>0$, J. Wei and S. Yan in \cite{WeiYan} proved that \eqref{NLS} admits infinitely many non-radial positive solutions with a large number of bumps near the infinity, centered at the vertices of a regular polygon, being each bump close to the ground state of $-\Delta u+u=u^{p}$. See also \cite{delPinoWeiYao} and \cite{MollePassaseo} for the case of non radially symmetric potentials. 

As a further example, let us discuss another problem, concerning planar closed curves with prescribed curvature, given by a radial mapping $K\colon\R^{2}\to\R$ such that
\begin{equation}
\label{K-shape}
K(x)=1+\frac{A}{|x|^{\bta}}+O\left(\frac{1}{|x|^{\bta+1}}\right),\quad\text{as }|x|\to\infty\,,\quad x\in\R^{2}
\end{equation}
for some $A\ne 0$ and $\bta>1$. In \cite{CaldiroliCora} the authors constructed infinitely many closed curves with curvature $K$, obtained by a very large number of approximate unit circumferences which curl in such a way that their center moves along a circle of very large radius. 

In this work we show that, quite surprisingly and somehow contrary to what published in \cite{CaldiroliMusso}, for a class of mean curvature functions of the same form as \eqref{V-shape} and \eqref{K-shape}, and in fact quite more general, a similar multi-bubble construction cannot be accomplished, at least, to some extent. 

Let us set out the details of our result. We will focus on a class of prescribed mean curvatures $H$ of the form
\begin{equation}
\label{eq:H1}
H(X)=1+\frac{A(X/|X|)}{|X|^{\bta}}+O\left(\frac{1}{|X|^{\bta+\nu}}\right),\quad\text{as }|X|\to\infty\,,\quad X\in\R^{3}
\end{equation}
with $A\colon\mathbb{S}^{2}\to\mathbb{R}$ regular enough, $\bta>0$ and $\nu>0$. The choice of this class of functions is motivated by previous works on the $H$-bubble problem (see, e.g., \cite{CaldiroliMusina2002}). Also the decay of inverse-power type seems reasonable for a wide class of geometric problems (see, e.g., \cite{BonheureIacopetti, TreibergsWei}). 

As reference surfaces we consider Delaunay tori, which are nothing but long segments of Delaunay cylinders bended along a circle, in order to give a compact surface with genus one. In terms of the parameterization, we can describe this class of surfaces as follows. We firstly introduce the Delaunay cylinders, in isothermal coordinates, given by 
\begin{equation}
\label{eq:Xa-vector}
X_{a}(t,\theta)=\left[\begin{array}{c}
x_{a}(t)\cos\theta\\
x_{a}(t)\sin\theta\\
z_{a}(t)\end{array}\right]
\end{equation}
where $a\in[-\frac{1}{2},\infty)\setminus\{0\}$ is the Delaunay parameter, and $x_{a}$ and $z_{a}$ solve
\begin{equation}\label{eq:xz}
\left\{\begin{array}{l}
x''_{a}=(1+2\gamma_{a})x_{a}-2x_{a}^{3}\\
x_{a}(0)=1+a\\
x'_{a}(0)=0
\end{array}\right.\quad
\left\{\begin{array}{l}
z'_{a}=x_{a}^{2}-\gamma_{a}\\
z_{a}(0)=0
\end{array}\right.
\quad\text{where}\quad \gamma_{a}=a(1+a)\,.
\end{equation}
For every $a\in[-\frac{1}{2},\infty)\setminus\{0\}$ the mapping $X_{a}$ defines a parametric surface $\Sigma_{a}$ with mean curvature 1, axis of revolution on the $x_{3}$-axis, and neck-size $|a|$ (the neck-size is defined as the distance of $\Sigma_{a}$ from the symmetry axis). As $-\frac{1}{2}\le a<0$ the surfaces are embedded and are called \textit{unduloids}, as $a>0$ the surfaces have self-intersections and are called \textit{nodoids}. In fig.~1 one can see the support of the curve $t\mapsto(z_{a}(t),x_{a}(t))$ in the two cases. When $a\to 0$ the limit configuration of the surface $\Sigma_{a}$ corresponds to a row of infinitely many tangent unit spheres placed along the $x_{3}$-axis. Following \cite{MazzeoPacard2001}, the case $a\to 0$ is sometimes called the \textit{singular limit}.  
\bigskip


\noindent
\begin{tabular}{cc}
\begin{tikzpicture}[scale=1.15]
\draw[->](-1.2767,0)--(0.2+3*1.2767,0);
\draw[->](-1.2767,0)--(-1.2767,1.5);
\draw[dashed](-1.2767,0.1000)--(3*1.2767,0.1000);
\draw[dashed](-1.2767,0.9000)--(3*1.2767,0.9000);
\filldraw(-1.2767,0.1000)circle(0.4pt)node[left]{\small{$|a|$}};
\filldraw(-1.2767,0.9000)circle(0.4pt)node[left]{\small{$1\!-\!|a|$}};
\draw[thick,smooth,xshift=-1.2767cm](0.0000,0.9000)--(0.0280,0.8997)--(0.0429,0.8992)--(0.0728,0.8977)--(0.1025,0.8954)--(0.1320,0.8923)--(0.1612,0.8885)--(0.1756,0.8863)--(0.1900,0.8839)--(0.2184,0.8787)--(0.2324,0.8758)--(0.2463,0.8728)--(0.2738,0.8663)--(0.3007,0.8592)--(0.3270,0.8516)--(0.3527,0.8436)--(0.3778,0.8350)--(0.4022,0.8261)--(0.4259,0.8168)--(0.4490,0.8071)--(0.4713,0.7973)--(0.4930,0.7871)--(0.5140,0.7768)--(0.5343,0.7663)--(0.5539,0.7557)--(0.5821,0.7396)--(0.6088,0.7233)--(0.6341,0.7071)--(0.6580,0.6908)--(0.6879,0.6694)--(0.7090,0.6535)--(0.7353,0.6326)--(0.7597,0.6122)--(0.7823,0.5923)--(0.8085,0.5682)--(0.8370,0.5405)--(0.8667,0.5099)--(0.8932,0.4811)--(0.9233,0.4466)--(0.9525,0.4115)--(0.9781,0.3794)--(1.0008,0.3500)--(1.0232,0.3204)--(1.0434,0.2934)--(1.0663,0.2644)--(1.0870,0.2362)--(1.1089,0.2091)--(1.1258,0.1894)--(1.1439,0.1700)--(1.1620,0.1526)--(1.1775,0.1394)--(1.1876,0.1318)--(1.2039,0.1212)--(1.2193,0.1132)--(1.2351,0.1069)--(1.2497,0.1029)--(1.2629,0.1008)--(1.2712,0.1001)--(1.2767,0.1000);
\draw[thick,smooth,xshift=1.2767cm](-0.0000,0.9000)--(-0.0280,0.8997)--(-0.0429,0.8992)--(-0.0728,0.8977)--(-0.1025,0.8954)--(-0.1320,0.8923)--(-0.1612,0.8885)--(-0.1756,0.8863)--(-0.1900,0.8839)--(-0.2184,0.8787)--(-0.2324,0.8758)--(-0.2463,0.8728)--(-0.2738,0.8663)--(-0.3007,0.8592)--(-0.3270,0.8516)--(-0.3527,0.8436)--(-0.3778,0.8350)--(-0.4022,0.8261)--(-0.4259,0.8168)--(-0.4490,0.8071)--(-0.4713,0.7973)--(-0.4930,0.7871)--(-0.5140,0.7768)--(-0.5343,0.7663)--(-0.5539,0.7557)--(-0.5821,0.7396)--(-0.6088,0.7233)--(-0.6341,0.7071)--(-0.6580,0.6908)--(-0.6879,0.6694)--(-0.7090,0.6535)--(-0.7353,0.6326)--(-0.7597,0.6122)--(-0.7823,0.5923)--(-0.8085,0.5682)--(-0.8370,0.5405)--(-0.8667,0.5099)--(-0.8932,0.4811)--(-0.9233,0.4466)--(-0.9525,0.4115)--(-0.9781,0.3794)--(-1.0008,0.3500)--(-1.0232,0.3204)--(-1.0434,0.2934)--(-1.0663,0.2644)--(-1.0870,0.2362)--(-1.1089,0.2091)--(-1.1258,0.1894)--(-1.1439,0.1700)--(-1.1620,0.1526)--(-1.1775,0.1394)--(-1.1876,0.1318)--(-1.2039,0.1212)--(-1.2193,0.1132)--(-1.2351,0.1069)--(-1.2497,0.1029)--(-1.2629,0.1008)--(-1.2712,0.1001)--(-1.2767,0.1000);
\draw[thick,smooth,xshift=1.2767cm](0.0000,0.9000)--(0.0280,0.8997)--(0.0429,0.8992)--(0.0728,0.8977)--(0.1025,0.8954)--(0.1320,0.8923)--(0.1612,0.8885)--(0.1756,0.8863)--(0.1900,0.8839)--(0.2184,0.8787)--(0.2324,0.8758)--(0.2463,0.8728)--(0.2738,0.8663)--(0.3007,0.8592)--(0.3270,0.8516)--(0.3527,0.8436)--(0.3778,0.8350)--(0.4022,0.8261)--(0.4259,0.8168)--(0.4490,0.8071)--(0.4713,0.7973)--(0.4930,0.7871)--(0.5140,0.7768)--(0.5343,0.7663)--(0.5539,0.7557)--(0.5821,0.7396)--(0.6088,0.7233)--(0.6341,0.7071)--(0.6580,0.6908)--(0.6879,0.6694)--(0.7090,0.6535)--(0.7353,0.6326)--(0.7597,0.6122)--(0.7823,0.5923)--(0.8085,0.5682)--(0.8370,0.5405)--(0.8667,0.5099)--(0.8932,0.4811)--(0.9233,0.4466)--(0.9525,0.4115)--(0.9781,0.3794)--(1.0008,0.3500)--(1.0232,0.3204)--(1.0434,0.2934)--(1.0663,0.2644)--(1.0870,0.2362)--(1.1089,0.2091)--(1.1258,0.1894)--(1.1439,0.1700)--(1.1620,0.1526)--(1.1775,0.1394)--(1.1876,0.1318)--(1.2039,0.1212)--(1.2193,0.1132)--(1.2351,0.1069)--(1.2497,0.1029)--(1.2629,0.1008)--(1.2712,0.1001)--(1.2767,0.1000);
\draw[thick,smooth,xshift=3*1.2767cm](-0.0000,0.9000)--(-0.0280,0.8997)--(-0.0429,0.8992)--(-0.0728,0.8977)--(-0.1025,0.8954)--(-0.1320,0.8923)--(-0.1612,0.8885)--(-0.1756,0.8863)--(-0.1900,0.8839)--(-0.2184,0.8787)--(-0.2324,0.8758)--(-0.2463,0.8728)--(-0.2738,0.8663)--(-0.3007,0.8592)--(-0.3270,0.8516)--(-0.3527,0.8436)--(-0.3778,0.8350)--(-0.4022,0.8261)--(-0.4259,0.8168)--(-0.4490,0.8071)--(-0.4713,0.7973)--(-0.4930,0.7871)--(-0.5140,0.7768)--(-0.5343,0.7663)--(-0.5539,0.7557)--(-0.5821,0.7396)--(-0.6088,0.7233)--(-0.6341,0.7071)--(-0.6580,0.6908)--(-0.6879,0.6694)--(-0.7090,0.6535)--(-0.7353,0.6326)--(-0.7597,0.6122)--(-0.7823,0.5923)--(-0.8085,0.5682)--(-0.8370,0.5405)--(-0.8667,0.5099)--(-0.8932,0.4811)--(-0.9233,0.4466)--(-0.9525,0.4115)--(-0.9781,0.3794)--(-1.0008,0.3500)--(-1.0232,0.3204)--(-1.0434,0.2934)--(-1.0663,0.2644)--(-1.0870,0.2362)--(-1.1089,0.2091)--(-1.1258,0.1894)--(-1.1439,0.1700)--(-1.1620,0.1526)--(-1.1775,0.1394)--(-1.1876,0.1318)--(-1.2039,0.1212)--(-1.2193,0.1132)--(-1.2351,0.1069)--(-1.2497,0.1029)--(-1.2629,0.1008)--(-1.2712,0.1001)--(-1.2767,0.1000);
\end{tikzpicture}
&
\begin{tikzpicture}[scale=1.15]
\draw[->](-2*0.6403,0)--(0.2+6*0.6403,0);
\draw[white](-2*0.6403,0)--(-2*0.6403,-0.135);
\draw[->](-2*0.6403,0)--(-2*0.6403,1.55);
\draw[dashed](-2*0.6403,0.1727)--(6*0.6403,0.1727);
\draw[dashed](-2*0.6403,1.1727)--(6*0.6403,1.1727);
\filldraw(-2*0.6403,0.1727)circle(0.4pt)node[left]{\small{$a$}};
\filldraw(-2*0.6403,1.1727)circle(0.4pt)node[left]{\small{$1$+$a$}};
\draw[thick,smooth](0,1.1727)--(0.0650,1.1703)--(0.1292,1.1630)--(0.1918,1.1513)--(0.2402,1.1389)--(0.2983,1.1199)--(0.3531,1.0977)--(0.4043,1.0728)--(0.4519,1.0457)--(0.5119,1.0051)--(0.5646,0.9626)--(0.6103,0.9194)--(0.6443,0.8824)--(0.6739,0.8459)--(0.7111,0.7926)--(0.7436,0.7364)--(0.7725,0.6735)--(0.7954,0.6066)--(0.8125,0.5264)--(0.8178,0.4673)--(0.8155,0.3998)--(0.7985,0.3196)--(0.7779,0.2713)--(0.7454,0.2250)--(0.7047,0.1918)--(0.6632,0.1755)--(0.6403,0.1727);
\draw[thick,smooth](0,1.1727)--(-0.0650,1.1703)--(-0.1292,1.1630)--(-0.1918,1.1513)--(-0.2402,1.1389)--(-0.2983,1.1199)--(-0.3531,1.0977)--(-0.4043,1.0728)--(-0.4519,1.0457)--(-0.5119,1.0051)--(-0.5646,0.9626)--(-0.6103,0.9194)--(-0.6443,0.8824)--(-0.6739,0.8459)--(-0.7111,0.7926)--(-0.7436,0.7364)--(-0.7725,0.6735)--(-0.7954,0.6066)--(-0.8125,0.5264)--(-0.8178,0.4673)--(-0.8155,0.3998)--(-0.7985,0.3196)--(-0.7779,0.2713)--(-0.7454,0.2250)--(-0.7047,0.1918)--(-0.6632,0.1755)--(-0.6403,0.1727);
\draw[thick,smooth,xshift=-2*0.6403cm](0,1.1727)--(0.0650,1.1703)--(0.1292,1.1630)--(0.1918,1.1513)--(0.2402,1.1389)--(0.2983,1.1199)--(0.3531,1.0977)--(0.4043,1.0728)--(0.4519,1.0457)--(0.5119,1.0051)--(0.5646,0.9626)--(0.6103,0.9194)--(0.6443,0.8824)--(0.6739,0.8459)--(0.7111,0.7926)--(0.7436,0.7364)--(0.7725,0.6735)--(0.7954,0.6066)--(0.8125,0.5264)--(0.8178,0.4673)--(0.8155,0.3998)--(0.7985,0.3196)--(0.7779,0.2713)--(0.7454,0.2250)--(0.7047,0.1918)--(0.6632,0.1755)--(0.6403,0.1727);
\draw[thick,xshift=2*0.6403cm](0,1.1727)--(-0.0650,1.1703)--(-0.1292,1.1630)--(-0.1918,1.1513)--(-0.2402,1.1389)--(-0.2983,1.1199)--(-0.3531,1.0977)--(-0.4043,1.0728)--(-0.4519,1.0457)--(-0.5119,1.0051)--(-0.5646,0.9626)--(-0.6103,0.9194)--(-0.6443,0.8824)--(-0.6739,0.8459)--(-0.7111,0.7926)--(-0.7436,0.7364)--(-0.7725,0.6735)--(-0.7954,0.6066)--(-0.8125,0.5264)--(-0.8178,0.4673)--(-0.8155,0.3998)--(-0.7985,0.3196)--(-0.7779,0.2713)--(-0.7454,0.2250)--(-0.7047,0.1918)--(-0.6632,0.1755)--(-0.6403,0.1727);
\draw[thick,smooth,xshift=2*0.6403cm](0,1.1727)--(0.0650,1.1703)--(0.1292,1.1630)--(0.1918,1.1513)--(0.2402,1.1389)--(0.2983,1.1199)--(0.3531,1.0977)--(0.4043,1.0728)--(0.4519,1.0457)--(0.5119,1.0051)--(0.5646,0.9626)--(0.6103,0.9194)--(0.6443,0.8824)--(0.6739,0.8459)--(0.7111,0.7926)--(0.7436,0.7364)--(0.7725,0.6735)--(0.7954,0.6066)--(0.8125,0.5264)--(0.8178,0.4673)--(0.8155,0.3998)--(0.7985,0.3196)--(0.7779,0.2713)--(0.7454,0.2250)--(0.7047,0.1918)--(0.6632,0.1755)--(0.6403,0.1727);
\draw[thick,xshift=4*0.6403cm](0,1.1727)--(-0.0650,1.1703)--(-0.1292,1.1630)--(-0.1918,1.1513)--(-0.2402,1.1389)--(-0.2983,1.1199)--(-0.3531,1.0977)--(-0.4043,1.0728)--(-0.4519,1.0457)--(-0.5119,1.0051)--(-0.5646,0.9626)--(-0.6103,0.9194)--(-0.6443,0.8824)--(-0.6739,0.8459)--(-0.7111,0.7926)--(-0.7436,0.7364)--(-0.7725,0.6735)--(-0.7954,0.6066)--(-0.8125,0.5264)--(-0.8178,0.4673)--(-0.8155,0.3998)--(-0.7985,0.3196)--(-0.7779,0.2713)--(-0.7454,0.2250)--(-0.7047,0.1918)--(-0.6632,0.1755)--(-0.6403,0.1727);
\draw[thick,smooth,xshift=4*0.6403cm](0,1.1727)--(0.0650,1.1703)--(0.1292,1.1630)--(0.1918,1.1513)--(0.2402,1.1389)--(0.2983,1.1199)--(0.3531,1.0977)--(0.4043,1.0728)--(0.4519,1.0457)--(0.5119,1.0051)--(0.5646,0.9626)--(0.6103,0.9194)--(0.6443,0.8824)--(0.6739,0.8459)--(0.7111,0.7926)--(0.7436,0.7364)--(0.7725,0.6735)--(0.7954,0.6066)--(0.8125,0.5264)--(0.8178,0.4673)--(0.8155,0.3998)--(0.7985,0.3196)--(0.7779,0.2713)--(0.7454,0.2250)--(0.7047,0.1918)--(0.6632,0.1755)--(0.6403,0.1727);
\draw[thick,xshift=6*0.6403cm](0,1.1727)--(-0.0650,1.1703)--(-0.1292,1.1630)--(-0.1918,1.1513)--(-0.2402,1.1389)--(-0.2983,1.1199)--(-0.3531,1.0977)--(-0.4043,1.0728)--(-0.4519,1.0457)--(-0.5119,1.0051)--(-0.5646,0.9626)--(-0.6103,0.9194)--(-0.6443,0.8824)--(-0.6739,0.8459)--(-0.7111,0.7926)--(-0.7436,0.7364)--(-0.7725,0.6735)--(-0.7954,0.6066)--(-0.8125,0.5264)--(-0.8178,0.4673)--(-0.8155,0.3998)--(-0.7985,0.3196)--(-0.7779,0.2713)--(-0.7454,0.2250)--(-0.7047,0.1918)--(-0.6632,0.1755)--(-0.6403,0.1727);
\end{tikzpicture}
\\
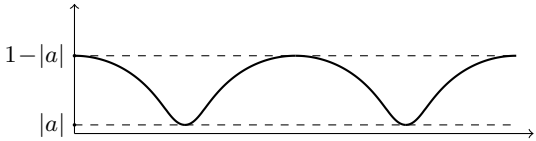
\small{Fig.~\!1-a. Profile curve of an unduloid ($-\frac{1}{2}\le a<0$) }~&~~
\small{Fig.~\!1-b. Profile curve of a nodoid ($a>0$)}
\end{tabular}
\bigskip

\noindent
We also point out that the function $x_{a}$ solving \eqref{eq:xz} is periodic of some period $2\tau_{a}>0$, whereas $z_{a}(t+2\tau_a)=2h_a+z_{a}(t)$ $(t\in\mathbb{R}$) for some $h_{a}>0$, with $\tau_{a}\to\infty$ and $h_{a}\to 1$, as $a\to 0$ (see Lemma \ref{L:tau-a-expansion}).

Then, fixing $a\in[-\frac{1}{2},\infty)\setminus\{0\}$ and $n\in\mathbb{N}$ sufficiently large, we define compact surfaces with genus one, setting
\begin{equation}
\label{eq:toroidal-nodoid}
X_{n,a}(t,\theta):=\left[\begin{array}{c}x_{a}(t)\cos\theta\\ \left(\frac{nh_{a}}{\pi}+x_{a}(t)\sin\theta\right)\cos\frac{\pi z_{a}(t)}{nh_{a}}\\[6pt] \left(\frac{nh_{a}}{\pi}+x_{a}(t)\sin\theta\right)\sin\frac{\pi z_{a}(t)}{nh_{a}}\end{array}\right].
\end{equation}
The map ${X}_{n,a}$ turns out to be doubly periodic with respect to the rectangle 
\[
K_{n,a}:=[-n\tau_{a},n\tau_{a}]\times[-\pi,\pi]
\] 
and provides a parametrization of a toroidal surface $\Sigma_{n,a}$ obtained by cutting a section of the Delaunay surface made by $n$ lobes, and bending it along a circumference contained in the plane $X\cdot\mathbf{e}_{1}=0$, with center at the origin and length $2nh_{a}$. The surface $\Sigma_{n,a}$ defines the so-called reference surface and it can be named \textit{Delaunay torus}.

We focus on a class of parametric surfaces which are normal graphs over $\Sigma_{n,a}$, that is, surfaces admitting parameterizations of the form
\[
X(t,\theta)=X_{n,a}(t,\theta)+\varphi(t,\theta)N_{n,a}(t,\theta)
\]
where 
\[
N_{n,a}:=\frac{\partial_{t}X_{n,a}\wedge\partial_{\theta}X_{n,a}}{|\partial_{t}X_{n,a}\wedge\partial_{\theta}X_{n,a}|}
\]
is the normal versor to $X_{n,a}$ and $\varphi\colon\R^{2}\to\R$ is a function of class $C^{2}$, doubly periodic with respect to the rectangle $K_{n,a}$. The mapping $\varphi$ rules the amplitude of the perturbation of the reference surface $\Sigma_{n,a}$ in the normal direction and has to be small in a suitable sense. Its smallness can vary in correspondence to the bulks or to the necks, according to a somehow accurate pointwise control. This is expressed in terms of the following weighted norm, involving also the derivatives, up to the second order:
\begin{equation}
\label{eq:weighted-C2-norm}
\|\varphi\|_{C^{2}(K_{n,a};x_{a}^{-1})}:=\sup_{(t,\theta)\in K_{n,a}}x_{a}(t)^{-1}\left(|\varphi(t,\theta)|+|\nabla\varphi(t,\theta)|+|D^{2}\varphi(t,\theta)|\right)\,.
\end{equation}

The weight $x_{a}^{-1}$ for the perturbation $\varphi$ is indeed necessary to ensure that $X_{n,a}+\varphi N_{n,a}$ is a parameterization of a regular surface, uniformly with respect to $a$ in a bounded set of the form $[-a_{0},a_{0}]\setminus\{0\}$ and $n\ge n_{0}$ (see Remarks \ref{R:sign-normal} and \ref{R:sign-normal-eps}). The same weight even for the first derivatives of $\varphi$ is the right one if one aims to control the distance between the Gauss map of $X_{n,a}+\varphi N_{n,a}$ and that one of $X_{n,a}$, uniformly with respect to the parameters $a\in [-a_{0},a_{0}]\setminus\{0\}$ and $n\ge n_{0}$ (see Remark \ref{R:weight-C1}). A similar control for the mean curvature necessarily involves also the second derivatives, and the weight $x_{a}^{-1}$ is enough for our purposes. In view of this, we will be allowed to write a first order expansion of the mean curvature of $X_{n,a}+\varphi N_{n,a}$ about the reference surface $X_{n,a}$ (see Sect.\,\ref{S:expansions-normal-graphs}). 

We also point out that controlling $\varphi$ with the weight $x_{a}^{-1}$ guarantees that, starting from unduloids, the embeddedness property is preserved. More precisely, when $-\frac{1}{2}<a<0$, for $n$ large enough also $\Sigma_{n,a}$ is embedded, as well as surfaces parameterized by $X_{n,a}+\varphi N_{n,a}$ if $\|x_{a}^{-1}\varphi\|_{C^{0}}$ is small enough (see \cite[Section 7]{CaldiroliMusso}, with easy adjustments). 
\medskip

Our main result can be stated as follows.

\begin{Theorem}
\label{mainteo}
Let $H\colon\mathbb{R}^{3}\to\mathbb{R}$ be a $C^{1}$ mapping satisfying \eqref{eq:H1} for some $A\in C^{1}(\mathbb{S}^{2})$, $\bta>0$ and $\nu>0$. In addition, assume that
\begin{gather}
\label{eq:H2}
A(\mathbf{e}_{2})\ne 0\\
\label{eq:H3}
\limsup_{|X|\to\infty}|X|^{\bta+1}|\nabla H(X)|<\infty\,.
\end{gather}
Then there cannot exist $\alpha\in(0,1]$, sequences $(\ak)_{k}\subset\mathbb{R}\setminus\{0\}$, $(\nk)_{k}\subset\mathbb{N}$ and $(\phik)_{k}\subset C^{2,\alpha}(\mathbb{R}^{2})$ such that:
\begin{itemize}
\item[(i)]
$\ak\to 0$ and $\nk \to \infty$, as $k\to \infty$;
\item[(ii)] $\phik=\phik(t,\theta)$ is doubly periodic of period $2\nk\tau_{\ak}$ with respect to $t$ and of period $2\pi$ with respect to $\theta$, and satisfies 
\begin{gather}
\label{eq:phi1}
\|\phik\|_{C^{2}(K_{\nk,\ak};x_{\ak}^{-1})}\le\frac{R}{n_{k}^{\bta}}\\
\label{eq:phi2}
[D^{2}\phik]_{\alpha,K_{\nk,\ak}}\le\frac{R}{n_{k}^{\bta}}
\end{gather}
for some constant $R>0$ independent of $k$;
\item[(iii)]
the surface $\widetilde{\Sigma}_{k}$ parameterized by $X_{\nk,\ak}+\phik N_{\nk,\ak}$ has mean curvature $H(X_{\nk,\ak}+\phik N_{\nk,\ak})$.
\end{itemize}
\end{Theorem}
In other words, we show that for $n$ large and $|a|$ small, a suitably small neighborhood of the Delaunay torus $\Sigma_{n,a}$ contains no parametric surface characterized as a normal graph over $\Sigma_{n,a}$ and with prescribed mean curvature $H$, where the smallness is made precise by the conditions \eqref{eq:phi1}--\eqref{eq:phi2}.
\medskip

Let us make some comments about the assumptions. The reason for which the direction of the versor $\mathbf{e}_{2}$ is preferred in condition \eqref{eq:H2} is due to the fact that the reference surfaces $X_{n,a}$ are defined with a lobe centered in the direction of $\mathbf{e}_{2}$. As a purely geometric condition independent of the parameterization we should require that $A$ has a definite sign on $\mathbb{S}^{2}$. The assumption \eqref{eq:H3} is consistent with \eqref{eq:H1} 
and essentially says that the decay rate in \eqref{eq:H1} holds in $C^{1}$. In fact, if $H(X)=1+A(X/|X|)|X|^{-\bta}$ with no extra term, then \eqref{eq:H3} is satisfied for free. 

The control in the condition \eqref{eq:phi1} takes account of the decay rate of $H$. The left-hand side in \eqref{eq:phi2} represents the H\"older constant of $D^{2}\phik$ in $K_{\nk,\ak}$ with exponent $\alpha$. 
\medskip

When $\bta>1$, Theorem \ref{mainteo} holds true even without assuming \eqref{eq:H2}. In fact, when $\bta>1$ the decay rates in the right-hand sides of \eqref{eq:phi1}--\eqref{eq:phi2} can be weakened, provided that a bound on the constant $R$ is assumed. More precisely, we have the following result:

\begin{Theorem}
\label{mainteo2}
Let $H\colon\mathbb{R}^{3}\to\mathbb{R}$ be a $C^{1}$ mapping satisfying \eqref{eq:H1} and \eqref{eq:H3} with $A\in C^{1}(\mathbb{S}^{2})$, $\bta>1$, and $\nu>0$. Then, there cannot exist $\alpha\in(0,1]$, sequences $(\ak)_{k}\subset\mathbb{R}\setminus\{0\}$, $(\nk)_{k}\subset\mathbb{N}$ and $(\phik)_{k}\subset C^{2,\alpha}(\mathbb{R}^{2})$ such that:
\begin{itemize}
\item[(i)]
$\ak\to 0$, $\nk \to \infty$, as $k\to \infty$;
\item[(ii)] $\phik=\phik(t,\theta)$ is doubly periodic of period $2\nk\tau_{\ak}$ with respect to $t$ and of period $2\pi$ with respect to $\theta$, and satisfies 
\begin{gather}
\label{eq:phi1+}
\|\phik\|_{C^{2}(K_{\nk,\ak};x_{\ak}^{-1})}\le\frac{R_{1}}{n_{k}}\\
\label{eq:phi2+}
[D^{2}\phik]_{\alpha,K_{\nk,\ak}}\le\frac{R}{n_{k}}
\end{gather}
for some constants $R_{1}\in\left(0,\frac{\pi}{2+\sqrt2}\right)$ and $R>0$, both independent of $k$;
\item[(iii)]
the surface $\widetilde{\Sigma}_{k}$ parameterized by $X_{\nk,\ak}+\phik N_{\nk,\ak}$ has mean curvature $H(X_{\nk,\ak}+\phik N_{\nk,\ak})$.
\end{itemize}
\end{Theorem}

Observe that when $\bta>1$ the conditions \eqref{eq:phi1}--\eqref{eq:phi2} imply \eqref{eq:phi1+}--\eqref{eq:phi2+} and then, in this case, Theorem \ref{mainteo} is a corollary of Theorem \ref{mainteo2}.
\medskip

The proof of Theorems \ref{mainteo} and \ref{mainteo2} is carried out by contradiction and relies on very careful estimates and expansions of the mean curvature operator evaluated at parametric surfaces of the form $X_{n,a}+\varphi N_{n,a}$ in the limit $n\to\infty$, uniformly with respect to $a$ in a bounded neighborhood of $0$, and with $\varphi$ small with respect to the norm \eqref{eq:weighted-C2-norm}. Then, we take the singular limit $a\to 0$ and we obtain a certain differential equation with no solution in some class of functions. 

The argument needs no particularly sophisticated technique and can be adapted to obtain similar non-existence results also for analogous prescribed almost constant mean curvature problems ruled by other classes of mean curvature functions. For example, $H(X)=1+H_{1}(\eps|X|)+o(\eps^{2})$, where $\eps$ is a smallness parameter, and $H_{1}$ vanishes at 1, which is also a nondegenerate critical point of $H_{1}$. 
\medskip

The very reason of our non-existence results can be explained, roughly speaking, as follows. We firstly observe that linearizing the mean curvature operator about any Delaunay surface yields a linear PDE solved by 
%
all the mappings deriving from the invariance groups. Among them, a key role is played by that one (the function $w_{a,0}^{+}(t)$ in Lemma \ref{L:Jacobi}) resulting from the invariance with respect to the Delaunay parameter (i.e., the neck-size, essentially). Such a mapping is not periodic, in fact neither bounded, and stays out from the space of periodic, admissible functions. In the singular limit, i.e., as $a\to 0$, one has that  $w_{a,0}^{+}(t)\to w_{0}(t)=-1+t\tanh t$ (see Lemma \ref{T:ker}). 

Moreover, we point out that the mean curvature of the Delaunay tori $\Sigma_{n,a}$ tends to 1 as $n\to\infty$, with an extra term which is proportional to the inverse of the distance from the origin (see Proposition \ref{prop:expmeancurv}). A balance with a prescribed mean curvature function of the form \eqref{eq:H1} can occur just for $\beta=1$. This can be detected by using the mapping $w_{0}$ as a test function. However, also when $\beta=1$ the contributions due to the first-order terms in the mean curvature operator and in the prescribed mean curvature function do not agree and again the function $w_{0}$ plays a key role. 
\medskip


The peculiarity of the Delaunay invariance and the possibility of taking the singular limit make the big difference with respect to the PDE problem on the nonlinear Schr\"odinger equation \eqref{NLS} (which indeed admits Delaunay-type solutions, but far from the singular limit, see \cite{Dancer, Malchiodi}), and also that one concerning curves with prescribed curvature and a very large number of spirals. 
%
%
%
\medskip

In fact, the problem of parametric surfaces with prescribed mean curvature, even in a perturbative setting, quite differs also with respect to that one of CMC surfaces and we suspect that the difficulty in applying Fredholm theory and perturbative techniques, like, e.g., in \cite{Kapouleas1990, Kapouleas1991, MazzeoPacard2001}, with the aim of getting existence results, is not due just to technical obstructions but hides more substantial features. This was already remarked, for different reasons, even on the $H$-bubble problem (see, e.g., \cite{CaldiroliMusina2004, CaldiroliMusina2006AdvDiffEq, CaldiroliMusina2006ARMA}). 
\medskip

In conclusion, our non-existence results, even if confined by some severe (but reasonable) restrictions, lead us to suspect that the Soap Bubble Theorem of Alexandrov recalled at the beginning has some rigidity in the class of almost constant mean curvature problems. 
\medskip

The paper is organized as follows: in Section 2 we recall some facts about Delaunay surfaces and on the corresponding Jacobi (i.e., the linearized CMC) operator. In Section 3 we move to Delaunay tori and establish some expansions on the mean curvature operator and on the corresponding linearized one about them. 
Before tackling the proof of our main results, we need other expansions both for the mean curvature operator and for the mapping $H$ evaluated at parametric surfaces which are normal graphs of Dealunay tori, with small variations with respect to the weighted norm \eqref{eq:weighted-C2-norm}. This is accomplished in Section 4. Finally, Section 5 contains the proof of Theorems \ref{mainteo} and \ref{mainteo2}. The paper is completed by two appendices: the first one collects some identities for orthogonal (not necessarily conformal) parameterizations. In the second one we write the first and second variation of the mean curvature operator with respect to normal variations for an arbitrary orthogonal parameterization. Even if the formula for the second variation is not used in the proof of the main results, it seems new and could be useful for further study.

\subsubsection*{Notation and abbreviations}
For reader's convenience, we collect here some notation and abbreviations used in this paper: 
\begin{itemize}
\item
$X\in\R^{3}$ and $\mathbf{e}_{i}$ denotes the $i$-th versor of the canonical basis of $\R^{3}$.  
\item
$\widehat X:=\frac{X}{|X|}$ and $\widecheck X:=X-(X\cdot\mathbf{e}_{1})\mathbf{e}_{1}$, for $X\in\R^3\setminus\{0\}$.
\item
When $Y$ is a vector or a matrix, $|Y|$ denotes its Euclidean norm. 
\item
$\R/_{\tau}=\R/\tau\mathbb{Z}$, where $\tau>0$, and a map $f\colon\R/_{2\tau}\times\R/_{2\pi}\to\R$ is a doubly periodic mapping with respect to the rectangle $[-\tau,\tau]\times[-\pi,\pi]$. 
\item
Given a map $F=F(t,\theta)\colon \R^2\to\R^d$, with $d\in\mathbb{N}$, $F_t$, $F_\theta$, $F_{tt}$, $F_{\theta\theta}$, $F_{t\theta}$ denote, respectively,  $\frac{\partial F}{\partial t}$, $\frac{\partial F}{\partial \theta}$, $\frac{\partial^2 F}{\partial t^2}$, $\frac{\partial^2 F}{\partial \theta^2}$, $\frac{\partial^2 F}{\partial t \partial \theta}$. Moreover, for real-valued mappings (i.e., when $d=1$), $\nabla F$ is the gradient and $D^{2}F$ the Hessian matrix. 
\item
Given $\alpha\in(0,1]$ and a map $F\colon K\subset\R^2\to\R^d$, with $d\in\mathbb{N}$, $[F]_{\alpha,K}$ denotes the H\"older constant of $F$ in $K$, with exponent $\alpha$. 
\item 
For $\sigma\in\mathbb{R}$, $R_{\sigma}$ denotes the rotation matrix of an angle $\sigma$ about the $x_{1}$-axis and $Q_{\sigma}=\frac{d}{d\sigma}R_{\sigma}$, i.e.,
\begin{equation}
\label{eq:rotation-matrix}
R_{\sigma}=\left[\begin{array}{ccc}1&0&0\\ 0&\cos\sigma&-\sin\sigma\\ 0&\sin\sigma&\cos\sigma\end{array}\right],\quad Q_{\sigma}=\left[\begin{array}{ccc}0&0&0\\ 0&-\sin\sigma&-\cos\sigma\\ 0&\cos\sigma&-\sin\sigma\end{array}\right].
\end{equation}
\item
$\Sigma_{a}$ is the Delaunay surface with parameterization \eqref{eq:Xa-vector}, being $x_{a}$ and $z_{a}$ solutions of \eqref{eq:xz}. The label $a$ is the Delaunay parameter and runs in $[-\frac{1}{2},\infty)\setminus\{0\}$. 
\item
$\tau_{a}$ is the half-period of the map $x_{a}$ appearing in the conformal representation \eqref{eq:Xa-vector} of $\Sigma_{a}$ and $h_{a}$ is the half-length of a lobe of $\Sigma_{a}$. 
\item
Given a regular surface $\Sigma$ in $\mathbb{R}^{3}$, with  parameterization $X=X(t,\theta)\colon\Omega\to\R^3$ and $\Omega\subset\mathbb{R}^{2}$ parameters domain, the normal vector is 
\[
{\normal}:=\frac{X_{t}\wedge X_{\theta}}{\left|X_{t}\wedge X_{\theta}\right|}~\!.
\]
The mean curvature of $\Sigma$ at a given point $X\in\Sigma$, denoted $\mathfrak{M}(X)$, is expressed in terms of the coefficients of the first and second fundamental form, in Gaussian notation,
\begin{equation}
\label{eq:Gaussian-notation}
\begin{array}{c}
{\E}=|X_{t}|^{2}~\!,\quad {\F}=X_{t}\cdot X_{\theta}~\!,\quad {\G}=|X_{\theta}|^{2}~\!,\\
{\L}=X_{tt}\cdot {\normal}~\!,\quad {\M}=X_{t\theta}\cdot {\normal}~\!,\quad {\N}=X_{\theta\theta}\cdot {\normal}~\!,
\end{array}
\end{equation}
as
\begin{equation}
\label{eq:mean-curvature-def}
\mathfrak{M}=\frac{{\E}{\N}-2{\F}{\M}+{\G}{\L}}{2({\E}{\G}-{\F}^{2})}~\!.
\end{equation}
\end{itemize}

\section{Delaunay surfaces}\label{S:Delaunay-surfaces}

In this Section we recall some properties of the Delaunay surfaces \cite{Delaunay}, defined in the Introduction by their parametric representations \eqref{eq:Xa-vector}--\eqref{eq:xz}. We will use just conformal parameterizations of them. This has advantages, both because allows us to address their description in a unified way, without distinction between unduloids and nodoids, and because the corresponding linearized problem turns out to be easier to handle. In addition, the coefficient $x_{a}$ in the conformal representation \eqref{eq:Xa-vector} of Delaunay surfaces constitute the natural weight in controlling perturbations in the normal direction. 

\subsection{Basic properties}

\begin{Lemma}
For every $a\in[-\frac{1}{2},\infty)\setminus\{0\}$  the surface $\Sigma_{a}$ with parameterization $X_{a}$ defined by \eqref{eq:Xa-vector}, with $x_{a},z_{a}\colon\mathbb{R}\to\mathbb{R}$ solving \eqref{eq:xz}, has constant mean curvature 1 and it is a surface of revolution about the $x_{3}$-axis. Moreover the parameterization $X_{a}$ is conformal, that is, satisfies \begin{equation}
\label{eq:X-conformal}
(X_{a})_{t}\cdot (X_{a})_{\theta}=0\,,\quad |(X_{a})_{t}|=|(X_{a})_{\theta}|\,.
\end{equation}
\end{Lemma}

\Proof
From \eqref{eq:xz} it follows that 
\begin{equation}
\label{eq:conformality}
x_{a}^{2}=(x'_{a})^{2}+(z'_{a})^{2}
\end{equation}
which turns out to be equivalent to \eqref{eq:X-conformal}. Moreover, one can compute the mean curvature of $X_{a}$ by \eqref{eq:Gaussian-notation}--\eqref{eq:mean-curvature-def} and, using \eqref{eq:xz} and \eqref{eq:conformality}, one obtains 
$\mathfrak{M}(X_{a})=(-x_{a}''z_{a}'+z_{a}''x_{a}'+x_{a}z_{a}')/(2x_{a}^{3})=1$. 
\qed

\begin{Remark}\label{R:Delaunay}
For $a\in[-\frac{1}{2},0)$, from \eqref{eq:xz} it follows that $z_{a}$ is strictly increasing in $\R$, hence $X_{a}$ is embedded. In this case $X_{a}$ is the parameterization of an {\rm unduloid}. In particular when $a=-\frac{1}{2}$ the surface parameterized by $X_{a}$ is a cylinder of radius $\frac{1}{2}$. For $a\in(0,\infty)$, from \eqref{eq:xz} it follows that $z'_{a}$ changes sign periodically. In this case the surface parametrized by $X_{a}$ is not embedded, it has self-intersections and is called {\rm nodoid}. When the parameter $a$ becomes large, the mutual intersection of distinct lobes becomes higher and higher. The parameter $a$ has a precise geometrical meaning: the neck-size of the surface $\Sigma_{a}$ is $|a|$, whereas the largest cross section is a circumference of radius $1+a$. When $a\to 0$ the mapping $X_{a}$ converges to a parameterization of a unit sphere centered at the origin, but geometrically the surface $\Sigma_{a}$ tends to a configuration given by a row of tangent unit spheres placed on the $x_{3}$-axis. 
\end{Remark}

In the next Lemma we collect some properties of $x_{a}$ and $z_{a}$ that will be useful in the sequel. Actually the following estimates \eqref{eq:xa-bound}--\eqref{eq:acca-a} are more accurate than what will need, but we state them in this form for the sake of completeness. 

\begin{Lemma}\label{L:tau-a-expansion}
For every $a\in[-\frac{1}{2},\infty)\setminus\{0\}$ one has that:
\begin{itemize}
\item[(i)] The mappings $x_{a}$ and $z'_{a}$ are even and periodic of period $2\tau_{a}>0$. Moreover $z_{a}$ is odd and there exists $h_{a}>0$ such that $z_{a}(t+2\tau_{a})=2h_{a}+z_{a}(t)$ for every $t\in\R$.
\item[(ii)]  $x'_{a}(t)<0$ for every $t\in(0,\tau_{a})$. Moreover
\begin{equation}
\label{xa-bounds}
 |a|\le x_{a}(t)\le 1+a\quad\forall t\in\R\,.
 \end{equation} 
 In addition, there exists $a_{0}\in(0,\frac{1}{2})$ such that for every $a\in[-a_{0},a_{0}]\setminus\{0\}$
\begin{equation}
\label{eq:xa-bound}
0<x_{a}(t)\le(1+|a|)\sqrt{\sech t}\quad\forall t\in[-\tau_{a},\tau_{a}]\,.
\end{equation}
\item[(iii)]  
The maps $x_{a}$ and $z_{a}$ depend in a smooth way on the parameter $a$. In particular
\[
x_{a}(t)\to\sech t\quad\text{and}\quad z_{a}(t)\to\tanh t\quad\text{in~~$C^{2}_{loc}(\R)$, as $a\to 0$\,.}
\]
\item[(iv)]  
The mapping $a\mapsto\tau_{a}$ is of class $C^{\infty}$ in $(-\frac{1}{2},\infty)\setminus\{0\}$ and 
\begin{equation}
\label{eq:tau-a}
\tau_{a}=-\log|a|+O(1)\quad\text{as~~$a\to 0$}\,.
\end{equation}
The mapping $a\mapsto h_{a}$ is continuous in $(-\frac{1}{2},\infty)$, of class $C^{\infty}$ in $(-\frac{1}{2},\infty)\setminus\{0\}$, and 
\begin{equation}
\label{eq:acca-a}
h_{a}=1+a\log|a|+O(a)\quad\text{as~~$a\to 0$}\,.
\end{equation}
 \end{itemize} 
\end{Lemma}

\Proof
Properties \textit{(i)}--\textit{(iii)} can be obtained by standard considerations about autonomous second order ode's and by the initial conditions in \eqref{eq:xz}. The proof of \eqref{eq:xa-bound} follows verbatim that one of \cite[Lemma 2.3]{CaldiroliMusso}, except for some obvious change in the notation ($a$ instead of $-a$). As far as concerns \textit{(iv)}, we observe that the function $x_{a}$ is strictly related to the \textit{Delta Amplitude function} with elliptic modulus $k^{2}$ which can be defined by
\begin{equation}
\label{eq:dn-Jacobi}
\mathrm{dn}(s,k)=\sqrt{1-k^{2}\sin^{2}\phi}~\!,\quad s=\int_{0}^{\phi}\frac{d\theta}{\sqrt{1-k^{2}\sin^{2}\theta}}
\end{equation}
and $k\in[0,1)$ (see \cite[Ch.~16]{METH}). Such function turns out to solve the Cauchy problem
\begin{equation}
\label{eq:CP-dn}
\left\{\begin{array}{l}y''=(2-k^{2})y-2y^{3}\\ y(0)=1\\ y'(0)=0~\!.\end{array}\right.
\end{equation}
(see \cite[Sect.~22]{DLMF}). 
Therefore, comparing \eqref{eq:CP-dn} with the Cauchy problem \eqref{eq:xz} for $x_{a}$, by uniqueness, we have that
\begin{equation}
\label{eq:xa=dn}
x_{a}(t)=(1+a)\mathrm{dn}((1+a)t,k_{a})\quad\text{where}\quad k_{a}=\sqrt{1-\frac{a^{2}}{(1+a)^{2}}}\,.
\end{equation}
It is known that the Delta Amplitude function is periodic in the $s$-variable, with period $2K(k)$, where
\[
K(k)=\int_{0}^{\frac{\pi}{2}}\frac{d\phi}{\sqrt{1-k^{2}\sin^{2}\phi}}
\]
is the complete elliptic integral of first type (see \cite[Ch.~16]{METH}). Then
\[
\tau_{a}=\frac{K(k_{a})}{1+a}\,.
\]
Moreover, integrating the equation for $z_{a}$ in \eqref{eq:xz} we obtain
\[
h_{a}=-\gamma_{a}\tau_{a}+\int_{0}^{\tau_{a}}x_{a}^{2}(t)\,dt\,.
\]
Using \eqref{eq:xa=dn} and making the change of variable $t\mapsto\phi$ through \eqref{eq:dn-Jacobi} we obtain
\[
\int_{0}^{\tau_{a}}x_{a}^{2}(t)\,dt=(1+a)\,E(k_{a})
\]
where
\[
E(k)=\int_{0}^{\frac{\pi}{2}}\sqrt{1-k^{2}\sin^{2}\phi}~\!d\phi
\]
is the complete elliptic integral of second type. Thus
\[
h_{a}=-\gamma_{a}\tau_{a}+(1+a)\,E(k_{a})=-a\,K(k_{a})+(1+a)\,E(k_{a})\,.
\]
The function $k\mapsto K(k)$ is $C^{\infty}$ in $[0,1)$ and singular at $k=1$, with a logarithmic behaviour:
\[
K(k)=-\log\sqrt{1-k^{2}}+O(1)\quad\text{as~~$k\to 1^{-}$}\,.
\]
The function $k\mapsto E(k)$ is $C^{1}$ in $[0,1]$ and $C^{\infty}$ in $[0,1)$. More precisely
\[
E(k)=1-\frac{1-k^{2}}{2}\log\sqrt{1-k^{2}}+(1-k^{2})\widetilde{E}(\sqrt{1-k^{2}})
\]
for some regular mapping $\widetilde{E}$ which is bounded, with bounded derivatives in $(0,1)$ (the properties of $K(k)$ and $E(k)$ used above can be found, e.g., in  \cite[Sect.~19]{DLMF}). Taking $k=k_{a}$ and observing that $a\to 0$ if and only if $k_{a}\to 1^{-}$, we obtain \eqref{eq:tau-a}--\eqref{eq:acca-a}.
\qed

Even if they will not be used in this work, for the sake of completeness, let us state (without proof) some expansions for the area and the volume integral of a lobe of a Delaunay surface $\Sigma_{a}$, defined respectively by
\[
\mathcal{A}(X_{a}):=\int_{K_{a}}|(X_{a})_{t}\wedge (X_{a})_{\theta}|\,dt\,d\theta\,,\quad \mathcal{V}(X_{a}):=\frac{1}{3}\int_{K_{a}}X_{a}\cdot (X_{a})_{t}\wedge (X_{a})_{\theta}\,dt\,d\theta
\]
where $K_{a}:=[-\tau_{a},\tau_{a}]\times[-\pi,\pi]$.

\begin{Lemma}
\label{L:area-volume} 
One has that
\begin{gather*}
\label{A-expansion}
\mathcal{A}(X_{a})=4\pi\left(1+a-\frac{a^{2}}{2}\log|a|+a^{2}\Phi_{0}(a)\right)
\\
\label{V-expansion}
\mathcal{V}(X_{a})=-\frac{4\pi}{3}\left(1+\frac{3a}{2}+a^{2}\Psi_{0}(a)\right)
\end{gather*}
where $\Phi_{0}$ and $\Psi_{0}$ are $C^{1}$ bounded functions with bounded derivative in $[-a_{0},a_{0}]\setminus\{0\}$ with $a_{0}\in\left(0,\frac{1}{2}\right)$ fixed.
\end{Lemma}
We refer to \cite{HMO,MO} for a proof in the case of unduloids but in fact the argument holds true for all Delaunay surfaces. 

\subsection{The Jacobi operator for Delaunay surfaces}
Fixing $a\in[-\frac{1}{2},\infty)\setminus\{0\}$, the Jacobi operator corresponding to the Delaunay surface $\Sigma_{a}$ is built by linearizing the mean curvature operator $\mathfrak{M}$ (defined in \eqref{eq:mean-curvature-def}) around the parameterization $X_{a}$ with respect to normal variations. More precisely, 
\[
\mathfrak{M}'_{a}(0)[\varphi]:=\frac{d}{ds}[\mathfrak{M}(X_{a}+s\varphi N_{a})]_{s=0}
\quad(\varphi\in C^{2}(\mathbb{R}^{2}))\quad\text{where}\quad N_{a}:=\frac{(X_{a})_{t}\wedge(X_{a})_{\theta}}{|(X_{a})_{t}\wedge(X_{a})_{\theta}|}\,.
\]
\begin{Remark}
\label{R:sign-normal}
Fixing $\varphi\in C^{2}(\mathbb{R}^{2})$ and $(t,\theta)\in\mathbb{R}^{2}$, one has that $[(X_{a}+s\varphi N_{a})_{t}\wedge(X_{a}+s\varphi N_{a})_{\theta}](t,\theta)\ne 0$ for $|s|$ small enough and then the above definition is well posed. In fact, using \eqref{eq:xz} and the identities \eqref{eq:propNetwedgeNeth}--\eqref{eq:propNetwedgeNeth-2}, one can compute
\begin{equation}
\label{sign-normal}
(X_{a}+\varphi N_{a})_{t}\wedge(X_{a}+\varphi N_{a})_{\theta}\cdot N_{a}=x_{a}^{2}\left[1-2\varphi+\frac{\varphi^{2}}{x_{a}^{2}}\left(x_{a}^{2}-\frac{\gamma_{a}^{2}}{x_{a}^{2}}\right)\right]\,.
\end{equation}
We note that the function $x_{a}^{2}-\frac{\gamma_{a}^{2}}{x_{a}^{2}}$ takes values in $[-1-2a,1+2a]$.
  However, fixing $a_{0}\in(0,\frac{1}{2})$ and $\varphi\in C^{2}(\mathbb{R}^{2})$, the coefficient $x_{a}^{-2}$ in front of $\varphi^{2}$ is not bounded uniformly with respect to $a\in[-a_{0},a_{0}]\setminus\{0\}$ (notice that $x_{a}(\tau_{a})=|a|$), hence it can produce large deviations, even when $\|\varphi\|_{C^{0}}$ is small. Therefore, if we aim to control the sign of \eqref{sign-normal} uniformly with respect to $a\in[-a_{0},a_{0}]\setminus\{0\}$ we need a bound on $\|x_{a}^{-1}\varphi\|_{C^{0}}$.
\end{Remark}

\begin{Lemma}\label{L:Jacobi}
One has that 
\[
\mathfrak{M}'_{a}(0)[\varphi]=\frac{1}{2x_{a}^{2}}\left[\Delta\varphi+2\left(x_{a}^{2}+\frac{\gamma_{a}^{2}}{x_{a}^{2}}\right)\varphi\right]
\] 
where $\Delta=\partial_{tt}+\partial_{\theta\theta}$ and $x_{a}$ and $\gamma_{a}$ are given by \eqref{eq:xz}.
\end{Lemma}
The proof of \ref{L:Jacobi} can be obtained by an easy adaptation of the computations made in \cite[Lemma 3.1]{CaldiroliMusso}. 


\begin{Lemma}
\label{Lemma-ker}
For every $a\in[-\frac{1}{2},\infty)\setminus\{0\}$ let
\[
w_{a,0}^{+}
=-\frac{z'_{a}}{x_{a}}\frac{\partial x_{a}}{\partial a}+\frac{x'_{a}}{x_{a}}\frac{\partial z_{a}}{\partial a}~\!,\quad w_{a,0}^{-}=\frac{x'_{a}}{x_{a}}~\!,\quad w_{a,1}^{+}=\frac{z'_{a}}{x_{a}}~\!,\quad w_{a,1}^{-}=x'_{a}+\frac{z_{a}z'_{a}}{x_{a}}
\]
and
\begin{equation}
\label{Wa}
\mathscr{W}_{a}:=\{w_{a,0}^{\pm}(t)\,,~~w_{a,1}^{\pm}(t)\cos\theta\,,~~w_{a,1}^{\pm}(t)\sin\theta\}\,.
\end{equation}
Then:
\begin{itemize}
\item[(i)]
Any linear combination of the functions in  $\mathscr{W}_{a}$ 
solves
\begin{equation}
\label{La=0}
\Delta\varphi+2\left(x_{a}^{2}+\frac{\gamma_{a}^{2}}{x_{a}^{2}}\right)\varphi=0\quad\text{on~~ }\mathbb{R}^{2}
\end{equation}
and is periodic with respect to $\theta$ with period $2\pi$.
\item[(ii)]
If $a\in\left(-\frac{1}{2},0\right)\cup\big(0,\frac{\sqrt{3}-1}{2}\big]$, then any solution of \eqref{La=0} which is doubly periodic with respect to a rectangle $[-\tau,\tau]\times[-\pi,\pi]$, with $\tau>0$, is a linear combination of the mappings $w_{a,0}^{-}$, $w_{a,1}^{+}\cos\theta$ and $w_{a,1}^{+}\sin\theta$ and $\tau=\tau_{a}$. 
\end{itemize}
\end{Lemma}

As explained in \cite{MazzeoPacard2001}, the functions in $\mathscr{W}_{a}\setminus\{w_{a,0}^{+}\}$ come from the invariance of the mean curvature with respect to Euclidean motions (translations and rotations). Instead, $w_{a,0}^{+}$ comes from the invariance with respect to the Delaunay parameter $a$. Notice that the mappings $w_{a,0}^{+}$ and $w_{a,1}^{-}$ are not periodic, if $a\in(-\frac{1}{2},\infty)\setminus\{0\}$, because of the presence of $\frac{\partial z_{a}}{\partial a}$ and $z_{a}$, respectively. 
\medskip

\Proof
Statement $(i)$ has been proved in \cite{CaldiroliMusso} and previously in \cite{MazzeoPacard2001}, with different notation. Let us discuss $(ii)$. Let $\varphi$ be a doubly periodic solution of \eqref{La=0}, with period $2\tau$ with respect to $t$ and period $2\pi$ with respect to $\theta$. We can express $\varphi$ by its Fourier series  
\begin{equation*}
\varphi(t,\theta)=\sum_{j\in\mathbb{Z}}\phi_{j}(t)\chi_{j}(\theta)
\end{equation*}
where
\[
\chi_{j}(\theta)=\left\{\begin{array}{ll}\frac{1}{\sqrt{\pi}}\cos(j\theta)&\text{for $j<0$}\\
\frac{1}{\sqrt{2\pi}}&\text{for $j=0$}\\
\frac{1}{\sqrt{\pi}}\sin(j\theta)&\text{for $j\ge 1$}\end{array}
\right.\quad\text{and}\quad\phi_{j}(t)=\int_{-\pi}^{\pi}\varphi(t,\theta)\chi_{j}(\theta)~\!d\theta~\!.
\]
Then, setting 
\[
p_{a}:=x_{a}^{2}+\frac{\gamma_{a}^{2}}{x_{a}^{2}}\,,
\]
and $K_{\tau}=[-\tau,\tau]\times[-\pi,\pi]$, since $\Delta\varphi(t,\theta)=\sum_{j\in\mathbb{Z}}\chi_{j}(\theta)[\phi_{j}''(t)-j^{2}\phi_{j}(t)]$, we infer that
\begin{equation}
\label{eq:intj}
0=\int_{K_{\tau}}[\Delta\varphi(t,\theta)+2 p_{a}(t)\varphi(t,\theta)]\chi_{j}(\theta)\,dt\,d\theta=-\int_{-\tau}^{\tau}\left[|\phi_{j}'|^{2}+\left(j^{2}-2 p_{a}(t)\right)|\phi_{j}|^{2}\right]\,dt\quad\forall j\in\mathbb{Z}\,.
\end{equation}
Recalling \eqref{xa-bounds}, we deduce that 
\[
\|p_{a}\|_{\infty}=\max_{s\in[|a|,1+a]}\left(s^{2}+\gamma_{a}^{2}s^{-2}\right)=a^{2}+(1+a)^{2}.
\] 
Therefore, since $a\in\left[-\frac{1}{2},0\right)\cup\big(0,\frac{\sqrt{3}-1}{2}\big]$ it holds that
\[
j^{2}-2 p_{a}(t)\ge 4-2\left(a^{2}+(1+a)^{2}\right)\ge 0\quad\forall |j|\ge 2~\!,
\]
and thus, by \eqref{eq:intj}, $\phi_{j}=0$ for all $|j|\ge 2$. For $j=0,\pm 1$ the equation 
\begin{equation}
\label{eq:phi-j}
\phi''=(j^{2}-2p_{a})\phi
\end{equation} 
admit exactly two independent solutions. In particular, $w_{a,j}^{\pm}$ solve \eqref{eq:phi-j} and are independent (see the proof of \cite[Lemma 4.3]{CaldiroliMusso}). Therefore the conclusion follows, taking account that only $w_{a,0}^{-}$ and $w_{a,1}^{+}$ are periodic.  
\qed

\begin{Lemma}
\label{T:ker}
One has that $w_{a,0}^{+}\to-1+t\tanh t$ in $C^{2}_{loc}(\R)$, as $a\to 0$. Moreover, fixing $a_{0}\in(0,\frac{1}{2})$, there exists $C>0$ such that $|w_{a,0}^{+}(t)|\le C(1+|t|)$ for every $t\in\R$ and for every $a\in [-a_{0},a_{0}]\setminus\{0\}$.
\end{Lemma}
The proof is the same as in \cite[Lemma 4.7]{CaldiroliMusso} with the already observed caution. In particular, since calculations in \cite{CaldiroliMusso} are made with $-a$ instead of $a$, the fuction $w_{a,0}^{+}$, obtained with a derivation with respect to $a$, turns out to have opposite sign to that one in \cite{CaldiroliMusso}.

\section{Delaunay tori}\label{S:Delaunay-tori}

Fix $a\in[-\frac{1}{2},\infty)\setminus\{0\}$ and $n\in\mathbb{N}$ sufficiently large. We denote by $\Sigma_{n,a}$ the surface obtained by cutting a section of the Delaunay surface $\Sigma_{a}$ made by $n$ periods, hence with length $2nh_{a}$, and by bending it along a circle placed on the plane $X\cdot\mathbf{e}_{1}=0$, centered at the origin, and with circumference $2nh_{a}$. The corresponding parametrization is given by the map ${X}_{n,a}$ defined in \eqref{eq:toroidal-nodoid}. 

By construction, the surface $\Sigma_{n,a}$ turns out to be topologically equivalent to a torus and we call it \textit{Delaunay torus}. 
The surface $\Sigma_{n,a}$ is symmetric with respect to rotations of an angle $\frac{2\pi}{n}$ about the $x_{1}$-axis, namely
\[
R_{\frac{2\pi}{n}}\Sigma_{n,a}=\Sigma_{n,a},
\]
where $R_{\sigma}$ denotes the rotation matrix defined in \eqref{eq:rotation-matrix}. It is also symmetric with respect to the plane $X\cdot\mathbf{e}_{1}=0$ and, in view of the initial conditions in \eqref{eq:xz}, also with respect to the plane $X\cdot\mathbf{e}_{3}=0$. 

The parameterization $X_{n,a}$ is $2n\tau_{a}$-periodic with respect to $t$ and $2\pi$-periodic with respect to $\theta$ and satisfies
\begin{equation}
\label{eq:discrete-symmetry}
{X}_{n,a}(t+2\tau_{a},\theta)=R_{\frac{2\pi}{n}}{X}_{n,a}(t,\theta)~\!.
\end{equation}

Instead of the parameter $n$, it is convenient to set
\[
\eps=\frac{\pi}{nh_{a}}
\]
and to write the parameterization $X_{n,a}$ in the form 
\begin{equation}
\label{eq:eps-toroidal-unduloid}
X_{\eps,a}(t,\theta)=\left[\begin{array}{c}x_{a}(t)\cos\theta\\ \left(\eps^{-1}+x_{a}(t)\sin\theta\right)\cos(\eps z_{a}(t))\\ \left(\eps^{-1}+x_{a}(t)\sin\theta\right)\sin(\eps z_{a}(t))\end{array}\right]
\end{equation}
having in mind that $\eps$ is a small and positive parameter. In fact, fixing $a_{0}\in(0,\frac{1}{2})$, the map $X_{\eps,a}$ defines a regular surface for every $a\in[-a_{0},a_{0}]\setminus\{0\}$ provided that $\eps$ is small enough, depending on $a_{0}$ (see \eqref{eq:normVwedgeW} and Remark \ref{rem:appendix1}). 

We also point out that considering the limit $\varepsilon\to 0$ is equivalent to $n\to\infty$, since $h_{a}$ is bounded and far from 0, uniformly with respect to $a\in[-a_{0},a_{0}]\setminus\{0\}$ (see Lemma \ref{L:tau-a-expansion}--(\textit{iv})). 

We denote $\Sigma_{\eps,a}$ the surface parameterized by $X_{\eps,a}$ and
\begin{equation}
\label{normale-eps-a}
N_{\eps,a}:=\frac{(X_{\eps,a})_{t}\wedge (X_{\eps,a})_{\theta}}{|(X_{\eps,a})_{t}\wedge (X_{\eps,a})_{\theta}|}
\end{equation}
the normal versor which, according with \eqref{eq:discrete-symmetry}, satisfies $N_{\eps,a}(t+2\tau_{a},\theta)=R_{\frac{2\pi}{n}}N_{\eps,a}(t,\theta)$. 

We observe that the Delaunay surface $\Sigma_{a}$ is the limit, up to a suitable natural translation, of $\Sigma_{\eps,a}$ in the limit $\eps\to 0$. More precisely
\begin{equation}
\label{limit-e2}
{X}_{\eps,a}-\eps^{-1}\mathbf{e}_{2}
\to
{X}_{a}\quad\text{as $\eps\to 0$, in $C^{1}_{loc}(\R\times[-\pi,\pi])$.}
\end{equation}

Notice that $X_{\eps,a}$ is just an orthogonal parameterization of $\Sigma_{\eps,a}$, i.e. $(X_{\eps,a})_{t}\cdot(X_{\eps,a})_{\theta}=0$, but it is not conformal. Indeed (see also \eqref{eq:X_tX_thetaVW})
\begin{equation}
\label{eq:no-conf}
\left|(X_{\eps,a})_{t}\right|^{2}-\left|(X_{\eps,a})_{\theta}\right|^{2}=
(z'_{a})^{2}\left[(1+\eps x_{a}\sin\theta)^{2}-1\right]\,.
\end{equation}
However $X_{\eps,a}$ becomes conformal in the limit $\eps\to 0$. In fact, by \eqref{eq:no-conf}, $\left|(X_{\eps,a})_{t}\right|^{2}-\left|(X_{\eps,a})_{\theta}\right|^{2}\to 0$ as $\eps\to 0$, uniformly on $\R\times[-\pi,\pi]$.

\subsection{An expansion for the mean curvature of Delaunay tori}
In this subsection we determine an expansion for the mean curvature $\mathfrak{M}(X_{\eps,a})$ of a Delaunay torus, where $X_{\eps,a}$ is the parametrization given by \eqref{eq:eps-toroidal-unduloid}. In order to obtain uniform estimates with respect to both the parameters $a$ and $\eps$, as it will be clear later (see Subsection~\ref{Ss:Jacobi}), it is convenient to study  $2x_a^2\mathfrak{M}(X_{\eps,a})$. The main result of this subsection is the following.

\begin{Proposition}\label{prop:expmeancurv}
Fixing $a_{0}\in(0,\frac{1}{2})$, there exists $\eps_0>0$ such that for all $a \in [-a_0,a_0]\setminus\{0\}$ and $\eps\in(0,\eps_0)$ one has
 \begin{equation}\label{eq:estmeanc}
2x_a^2\mathfrak{M}(X_{\eps,a})=2x_a^2+\eps\left(2x_a^3+2z_a^\prime x_a-4(z_a^\prime)^2 x_a-\frac{(z_a^\prime)^3}{x_a}-2\gamma_a\frac{(z_a^\prime)^2}{x_a}\right)\sin\theta + \eps^2 \sigma_{\eps,a},
\end{equation}
where $\sigma_{\eps,a}=\sigma_{\eps,a}(t,\theta)$ is doubly periodic with respect to the rectangle $K_{a}:=[-\tau_{a},\tau_{a}]\times[-\pi,\pi]$ and satisfies
\begin{equation}\label{eq:estremtermexpmeanc}
\|x_a^{-2}\sigma_{\eps,a}\|_{C^0(K_{a})}\leq C, 
\end{equation}
for some positive constant $C$ independent of $\eps$, $a$.
\end{Proposition}
 
For the proof of Proposition \ref{prop:expmeancurv} we need a preliminary technical result. We begin introducing some notation. Let
\begin{equation}\label{eq:defvectVW}
U_{\eps,a}:=\left[\begin{array}{c} x_a \cos\theta\\[10pt] \eps^{-1}+ x_a \sin\theta \\[10pt] 0\end{array}\right], \ \ \
V_{\eps,a}:=\left[\begin{array}{c} x_a^\prime \cos\theta\\[10pt] x_a^\prime \sin\theta \\[10pt] z_a^\prime(1+\eps x_a\sin\theta)\end{array}\right], \ \ \ W_{\eps,a}:=\left[\begin{array}{c} -x_a \sin\theta\\[10pt] x_a \cos\theta \\[10pt]0\end{array}\right],
\end{equation}
then in terms of the rotation matrix $R_\sigma$ (defined in \eqref{eq:rotation-matrix}) we can write $X_{\eps,a}$ as\
\begin{equation}\label{eq:eps-toroidal-unduloid_rotmatrix}
X_{\eps,a}=R_{\eps z_a} U_{\eps,a}\,,
\end{equation}
and we easily check that
\begin{equation}\label{eq:relXtXthetaR}
(X_{\eps,a})_t=R_{\eps z_a}V_{\eps,a},\ \ (X_{\eps,a})_\theta=R_{\eps z_a}W_{\eps,a}\,.
\end{equation}
Moreover in view of \eqref{eq:propRsigma} we get
\begin{equation}\label{eq:normRVW}
N_{\eps,a}=R_{\eps z_a} \frac{V_{\eps,a}\wedge W_{\eps,a}}{|V_{\eps,a}\wedge W_{\eps,a}|}.
\end{equation}
In particular we have
\begin{equation}
\label{eq:X_tX_thetaVW}
|(X_{\eps,a})_t|^2=|V_{\eps,a}|^2=x_a^2 + \eps\left( 2x_a (z_a^\prime)^2 \sin \theta\right) + \eps^2\left(z_a^\prime x_a\sin\theta\right)^2, \ \ |(X_{\eps,a})_\theta|^2=|W_{\eps,a}|^2=x_a^2\,,
\end{equation}
and thus, since $V_{\eps,a}\cdot W_{\eps,a}=0$, we infer that
\begin{equation}
\label{eq:normVwedgeW}
|V_{\eps,a}\wedge W_{\eps,a}|=|(X_{\eps,a})_t\wedge(X_{\eps,a})_\theta|=|(X_{\eps,a})_t| |(X_{\eps,a})_\theta|=x_a^2 \sqrt{1+2\eps \frac{(z_a^\prime)^2}{x_a}\sin\theta + \eps^2 (z_a^\prime)^2 \sin^2\theta}\,.
\end{equation}
Then, by direct computation, it holds that
\begin{equation}\label{eq:VwedgeW}
\frac{V_{\eps,a}\wedge W_{\eps,a}}{|V_{\eps,a}\wedge W_{\eps,a}|}=\frac{1}{\sqrt{1+2\eps \frac{(z_a^\prime)^2}{x_a}\sin\theta + \eps^2 (z_a^\prime)^2 \sin^2\theta}}\left(\left[\begin{array}{c} -\frac{z_a^\prime}{x_a} \cos\theta\\[10pt] -\frac{z_a^\prime}{x_a} \sin\theta \\[10pt] \frac{x_a^\prime}{x_a}\end{array}\right] + \eps \left[\begin{array}{c} - z_a^\prime \cos\theta \sin\theta\\[10pt] -z_a^\prime \sin^2\theta \\[10pt] 0 \end{array}\right]\right).
\end{equation}
\begin{Remark}
Due to the presence of the function $z_a$ in \eqref{eq:relXtXthetaR}, \eqref{eq:normRVW}  the vector-valued functions $\Xet$, $\Xeth$, $\Ne$ (as well as their partial derivatives) are not periodic with respect to the variable $t$. Nevertheless, since $z_{a}$ appears jut in the rotation matrix $R_{\eps z_{a}}$, their scalar products turn out to be doubly periodic, as will see in Lemmata \ref{lem:expsecderXN}--\ref{lem:expNsecX}.
\end{Remark}

\begin{Remark}\label{rem:appendix1}
By \eqref{eq:conformality} and \eqref{xa-bounds} the functions $x_a$, $x_a^\prime$, $z_a^\prime$, $\frac{z_a^\prime}{x_a}$, $\frac{x_a^\prime}{x_a}$ are uniformly bounded in $\R$ (and $2\tau_a$-periodic), whenever $a$ lies in a bounded subset of $[-\frac{1}{2},\infty)\setminus\{0\}$. The same happens for their higher order derivatives, taking account of \eqref{eq:xz}.
\end{Remark}
\begin{Remark}\label{rem:appendix2}
If $B\subset [-\frac{1}{2}, +\infty)\setminus\{0\}$ is a bounded subset  and $(f_a)_{a\in B}$ is a family of continuous functions satisfying $\|x_a^{-k}f_a\|_{C^0} \leq C$, for all $a \in B$, for some constants $C>0$, $k>0$ independent of $a$, then, for any $m\leq k$, $a\in B$  one has $\|x_a^{-m}f_a\|_{C^0} \leq C_1$, for some constant $C_1>0$ independent of $a$. In particular, one can take $C_1=C (1+\max_{a\in B}|a|)^{k-m}$.
\end{Remark}
Aiming to prove Proposition \ref{prop:expmeancurv}, the following estimates are useful.

\begin{Lemma}\label{lem:expsecderXN}
Fixing $a_{0}\in(0,\frac{1}{2})$, there exist $\eps_0>0$ and $C>0$ such that for all $a \in [-a_0,a_0]\setminus\{0\}$ and $\eps\in(0,\eps_0)$ one has
\begin{eqnarray}
&&\displaystyle (X_{\eps,a})_{tt}\cdot N_{\eps,a}=\displaystyle  x_a^2+\gamma_a + \eps \left(2x_a^3+x_az_a^\prime-2(z_a^\prime)^2x_a\right) \sin\theta + \eps^2 \sigma_{\eps,a}^{(1)}\label{eq:expXttN}\\[9pt]
&&\displaystyle (X_{\eps,a})_{\theta\theta}\cdot N_{\eps,a}=\displaystyle  z_a^\prime + \eps \left(x_az_a^\prime-\frac{(z_a^\prime)^3}{x_a}\right) \sin \theta + \eps^2  \sigma_{\eps,a}^{(2)}\label{eq:expXthetathetaN}\\[9pt]
&&\displaystyle (X_{\eps,a})_{t\theta}\cdot N_{\eps,a}=\displaystyle   \eps \left(x_a^\prime z_a^\prime\right) \cos\theta + \eps^2  \sigma_{\eps,a}^{(3)},\label{eq:expXtthetaN}
\end{eqnarray}
where $\sigma_{\eps,a}^{(i)}=\sigma_{\eps,a}^{(i)}(t,\theta)$ ($i=1,2,3$) is doubly periodic with respect to $K_{a}=[-\tau_{a},\tau_{a}]\times[-\pi,\pi]$ and satisfies
\begin{equation}\label{eq:stimarestoxaminus1}
\left\|x_a^{-2}{\sigma_{\eps,a}^{(i)}}\right\|_{C^0(K_{a})}\leq C.
\end{equation}
\end{Lemma}
\Proof
For every $a\in[-\frac{1}{2},\infty)\setminus\{0\}$ and $\eps>0$ small, from \eqref{eq:relXtXthetaR} we have 
\[
(X_{\eps,a})_{tt}=R_{\eps,a} (V_{\eps,a})_t + \eps z_a^\prime Q_{\eps z_a} V_{\eps,a}\,,
\]
where $Q_\sigma$ is the matrix defined in \eqref{eq:rotation-matrix}.
Then, in view of \eqref{eq:normRVW} and taking into account \eqref{eq:propRsigma}, \eqref{eq:proppartialRsigma} we deduce that
\[
(X_{\eps,a})_{tt}\cdot N_{\eps,a}= \underbrace{(V_{\eps,a})_t \cdot \frac{V_{\eps,a}\wedge W_{\eps,a}}{|V_{\eps,a}\wedge W_{\eps,a}|}}_{\textbf{(I)}}+\underbrace{\eps z_a^\prime  \frac{\textbf{e}_1\cdot [V_{\eps,a}\wedge(V_{\eps,a}\wedge W_{\eps,a})]}{|V_{\eps,a}\wedge W_{\eps,a}|}}_{\textbf{(II)}}\,.
\]

\textbf{Study of the term $(I)$}. Firstly, we notice that from the well known expansion of the function $y\mapsto (1+y)^{-1/2}$ near the origin one can write $\frac{1}{\sqrt{1+y}}=1-\frac{1}{2}y+g(y)$ in $[-\delta,\delta]$, for some $\delta>0$, with $g\in C^2([-\delta,\delta])$ such that $g(y)=O(y^2)$, $g^\prime(y)=O(y)$, $g^{\prime\prime}(y)=O(1)$.

Now, fixing $a_{0}\in(0,\frac12)$ and taking into account Remark \ref{rem:appendix1}, we can find $\eps_0>0$ such that for any $a\in [-a_0,a_0]\setminus\{0\}$ and $\eps \in (0,\eps_0)$ one has $|2\eps \frac{(z_a^\prime)^2}{x_a}\sin\theta + \eps^2 (z_a^\prime)^2 \sin^2\theta|\leq \frac{\delta}{2}$, for all $(t,\theta)\in\R^2$.
Hence, setting
\[
\eta_{\eps,a}(t,\theta):=\eps^{-2}g\left(2\eps \frac{(z_a^\prime)^2}{x_a}\sin\theta + \eps^2 (z_a^\prime)^2 \sin^2\theta\right)-\frac{1}{2}(z_a^\prime)^2 \sin^2\theta
\]
one has that $\eta_{\eps,a}=\eta_{\eps,a}(t,\theta)$ is doubly periodic and
\begin{equation}\label{eq:elemexpsqrt}
\frac{1}{\sqrt{1+2\eps \frac{(z_a^\prime)^2}{x_a}\sin\theta + \eps^2 (z_a^\prime)^2 \sin^2\theta}}=1-\eps \frac{(z_a^\prime)^2}{x_a}\sin\theta+\eps^2 \eta_{\eps,a}\,.
\end{equation}

Moreover, in view of Remark \ref{rem:appendix1} and exploiting \eqref{eq:xz} we infer that
$x_a^{-1}\frac{(z_a^\prime)^2}{x_a}$, $x_a^{-1}(z_a^\prime)^2$, $x_a^{-1}\left(\frac{(z_a^\prime)^2}{x_a}\right)^\prime$, $x_a^{-1}\left((z_a^\prime)^2\right)^\prime$, $x_a^{-1}\left(\frac{(z_a^\prime)^2}{x_a}\right)^{\prime\prime}$, $x_a^{-1}\left((z_a^\prime)^2\right)^{\prime\prime}$ are uniformly bounded. Therefore, taking into account the properties of $g$, up to choosing smaller $\eps_0>0$, for any $a\in [-a_0,a_0]\setminus\{0\}$ and $\eps \in (0,\eps_0)$ one has
\begin{equation}\label{eq:estremaindterm}
\|x_a^{-1}\eta_{\eps,a} \|_{C^0(K_{a})}+\|x_a^{-1}|\nabla \eta_{\eps,a}|\|_{C^0(K_{a})}+\|x_a^{-1}|D^2\eta_{\eps,a}| \|_{C^0(K_{a})}\leq C
\end{equation}
for some positive constant $C$ independent of $a$ and $\eps$.
 
At the end, from \eqref{eq:VwedgeW}, \eqref{eq:elemexpsqrt}, \eqref{eq:estremaindterm} we readily obtain
\begin{equation}\label{eq:VwedgeW_eps-toroidal}
\frac{V_{\eps,a}\wedge W_{\eps,a}}{|V_{\eps,a}\wedge W_{\eps,a}|}=\left[\begin{array}{c}-\frac{z^\prime_{a}}{x_a}\cos\theta\\[10pt] -\frac{z^\prime_{a}}{x_a} \sin\theta \\[10pt] \frac{x^\prime_{a}}{x_a}\end{array}\right] + \eps \left[\begin{array}{c} \left(- z^\prime_{a} + \frac{ (z^\prime_{a})^3}{ x^2_{a}}\right) \cos\theta \sin \theta \\[6pt]  \left(- z^\prime_{a} +  \frac{ (z^\prime_{a})^3}{ x^2_{a}}\right) \sin^2\theta\\[6pt] -\left(\frac{ z^\prime_{a}}{ x_{a}}\right)^2x^\prime_a \sin\theta\end{array}\right] + \eps^2 \Psi_{\eps,a}
\end{equation}
where $\Psi_{\eps,a}=\Psi_{\eps,a}(t,\theta)$ is a doubly periodic vector-valued function. Moreover, by definition and arguing as in \eqref{eq:estremaindterm}  we easily check that the components $(\Psi_{\eps,a})_j$, $j=1,2,3$ satisfy
\begin{equation}\label{eq:stimaremcompPsi}
\| x_a^{-1}(\Psi_{\eps,a})_j\|_{C^0(K_{a})} + \| x_a^{-1}|\nabla(\Psi_{\eps,a})_j|\|_{C^0(K_{a})}+ \| x_a^{-1}|D^2(\Psi_{\eps,a})_j|\|_{C^0(K_{a})}\leq C,
\end{equation}
for some positive constant $C$ independent of $\eps$, $a$. We notice that for the purposes of the present proof we need only a uniform bound on $x_a^{-1}(\Psi_{\eps,a})_j$, but the complete form of \eqref{eq:stimaremcompPsi} will be useful in the sequel (see the proofs of Lemma \ref{lem:expNtNthetasquared} and Lemma \ref{lem:expNsecX}).
\medskip

Now, since
\begin{equation}\label{eq:extVt}
 (V_{\eps,a})_t=\left[\begin{array}{c} x_a^{\prime\prime} \cos\theta\\[10pt] x_a^{\prime\prime} \sin\theta \\[10pt] z_a^{\prime\prime}\end{array}\right] + \eps \left[\begin{array}{c} 0\\[10pt] 0 \\[10pt] (z_a^{\prime\prime}x_a+z_a^\prime x_a^\prime)\sin\theta\end{array}\right]
 \end{equation}
then, exploiting \eqref{eq:VwedgeW_eps-toroidal}, taking into account \eqref{eq:xz} and using repeatedly  \eqref{eq:conformality}, we get that
\begin{equation}\label{eq:exptermI}
\begin{array}{lll}
&& \hspace{-12mm}\displaystyle (V_{\eps,a})_t\cdot \frac{V_{\eps,a}\wedge W_{\eps,a}}{|V_{\eps,a}\wedge W_{\eps,a}|}\\[10pt]
&=&\displaystyle\frac{-{x_a^{\prime\prime}z_a^\prime}+z_a^{\prime\prime}x_a^\prime}{x_a}+\eps\left(-x_a^{\prime\prime}{z_a^\prime} + \frac{x_a^{\prime\prime}(z_a^\prime)^3}{x_a^2}- \left(\frac{z_a^{\prime}}{x_a}\right)^2 x_a^\prime z_a^{\prime\prime} + z_a^{\prime\prime} x_a^\prime + \frac{z_a^\prime (x_a^\prime)^2}{x_a}\right)\sin\theta+\eps^2 \sigma^{(I)}_{\eps,a}\\[12pt]
&=&\displaystyle x_a^2+\gamma_a+\eps\left(-x_az_a^\prime + 2x_a^3 + \frac{(z_a^\prime)^3}{x_a}-2\frac{(z_a^\prime)^4}{x_a}- 2\frac{\left(z_a^{\prime}\right)^2(x_a^\prime)^2}{x_a}   + \frac{z_a^\prime (x_a^\prime)^2}{x_a}\right)\sin\theta+\eps^2\sigma^{(I)}_{\eps,a}\\[12pt]
&=&\displaystyle x_a^2+\gamma_a+\eps\left( 2x_a^3 - 2\left(z_a^{\prime}\right)^2x_a \right)\sin\theta+\eps^2 \sigma^{(I)}_{\eps,a},\\
\end{array}
\end{equation}
where the term $\sigma_{\eps,a}^{(I)}$ is doubly periodic. In addition, since $x_a^{-1}(V_{\eps,a})_t$ has uniformly bounded coefficients (see \eqref{eq:extVt} and Remark \ref{rem:appendix1}) and in view of \eqref{eq:stimaremcompPsi} we infer that $\sigma_{\eps,a}^{(I)}$  satisfies \eqref{eq:stimarestoxaminus1} (as already mentioned before, here we are using only the uniform $C^0$ bound on $x_a^{-1}(\Psi_{\eps,a})_j$, for $j=1,2,3$).\\

\textbf{Study of the term $(II)$}. We begin noticing that by \eqref{eq1:propcrossscal} one has
\[
V_{\eps,a}\wedge(V_{\eps,a}\wedge W_{\eps,a})=-|V_{\eps,a}|^{2} W_{\eps,a}\,.
\]
Then from \eqref{eq:defvectVW}, \eqref{eq:X_tX_thetaVW}, \eqref{eq:normVwedgeW}, \eqref{eq:elemexpsqrt} and Remark \ref{rem:appendix1}, we deduce that
\begin{equation}\label{eq:exptermII}
\eps z_a^\prime  \frac{\textbf{e}_1\cdot [V_{\eps,a}\wedge(V_{\eps,a}\wedge W_{\eps,a})]}{|V_{\eps,a}\wedge W_{\eps,a}|}=
%
\eps z_a^\prime x_a\sin\theta + \eps^2 \sigma_{\eps,a}^{(II)},
\end{equation}
where $\sigma_{\eps,a}^{(II)}$ is doubly periodic and satisfies the estimate \eqref{eq:stimarestoxaminus1}. Hence, from \eqref{eq:exptermI} and \eqref{eq:exptermII} we obtain  \eqref{eq:expXttN}.
\medskip

As far as concerns \eqref{eq:expXthetathetaN}, observing that $(X_{\eps,a})_{\theta\theta}\cdot N_{\eps,a}=(W_{\eps,a})_\theta \cdot \frac{V_{\eps,a}\wedge W_{\eps,a}}{|V_{\eps,a}\wedge W_{\eps,a}|}$, from \eqref{eq:VwedgeW_eps-toroidal}, taking into account Remark \ref{rem:appendix1}, and noticing that $x_a^{-1}(W_{\eps,a})_\theta$ has uniformly bounded coefficients, because
\[
(W_{\eps,a})_\theta= \left[\begin{array}{c}-{x_a}\cos\theta\\[10pt] -{x_a} \sin\theta \\[10pt] 0\end{array}\right]\,,
\] 
then we get that
\[
(X_{\eps,a})_{\theta\theta}\cdot N_{\eps,a}=z_a^\prime+\eps\left(z_a^\prime x_a -\frac{(z_a^\prime)^3}{x_a}\right)\sin\theta+\eps^2\sigma_{\eps,a}^{(2)}
\]
where $\sigma_{\eps,a}^{(2)}$  is doubly periodic and satisfies \eqref{eq:estremaindterm}.
\medskip

Finally, since $(X_{\eps,a})_{t\theta}=R_{\eps z_a} (V_{\eps,a})_\theta$  we get that $(X_{\eps,a})_{t\theta} \cdot N_{\eps,a}= (V_{\eps,a})_\theta \cdot \frac{V_{\eps,a}\wedge W_{\eps,a}}{|V_{\eps,a}\wedge W_{\eps,a}|}$. Hence, taking into account that
\[
(V_{\eps,a})_\theta= \left[\begin{array}{c}- x_a^\prime\sin\theta\\[3pt] x_a^\prime \cos\theta \\[3pt] 0\end{array}\right]+\eps\left[\begin{array}{c}0\\[3pt] 0 \\[3pt] z_a^\prime x_a \cos\theta\end{array}\right]
\]
and exploiting Remark \ref{rem:appendix1}, \eqref{eq:stimaremcompPsi}, we readily obtain \eqref{eq:expXtthetaN} and as in the previous cases the remainder term $\sigma_{\eps,a}^{(3)}$ is doubly periodic and satisfies \eqref{eq:estremaindterm}.
\qed

\noindent
\textit{Proof of Proposition \ref{prop:expmeancurv}.}
Fixing $a_0\in(0,\frac{1}{2})$, let $\eps_0>0$ be given by Lemma \ref{lem:expsecderXN}. From the definition of mean curvature (see \eqref{eq:mean-curvature-def}) and since $(X_{\eps,a})_t\cdot(X_{\eps,a})_\theta=0$, we have
\begin{equation}\label{eq:2xasqmeancurv}
2\mathfrak{M}(X_{\eps,a})=\left(\frac{(X_{\eps,a})_{tt}\cdot N_{\eps,a}}{|(X_{\eps,a})_t|^2}+\frac{(X_{\eps,a})_{\theta\theta}\cdot N_{\eps,a}}{|(X_{\eps,a})_\theta|^2}\right).
\end{equation}
In view \eqref{eq:X_tX_thetaVW}, \eqref{eq:expXttN} and arguing as in \eqref{eq:elemexpsqrt}, up to taking smaller $\eps_0$, if necessary, we get that
\begin{equation}\label{eq:meancurvpartI}
\begin{split}
\displaystyle x_a^2\frac{(X_{\eps,a})_{tt}\cdot N_{\eps,a}}{|(X_{\eps,a})_t|^2}&=\displaystyle\left(1-2\eps \frac{(z_a^\prime)^2}{x_a}\sin\theta+\eps^2 \eta_{\eps,a}^{(1)}\right) \left(x_a^2+\gamma_a + \eps \left(2x_a^3+x_az_a^\prime-2(z_a^\prime)^2x_a\right) \sin\theta + \eps^2 \sigma_{\eps,a}^{(1)}\right)\\
&\quad=\displaystyle x_a^2+\gamma_a + \eps \left(2x_a^3+x_az_a^\prime-4(z_a^\prime)^2x_a -2\frac{(z_a^\prime)^2}{x_a}\gamma_a\right) \sin\theta + \eps^2 \eta_{\eps,a}^{(2)}\,,
\end{split}
\end{equation}
where $\eta_{\eps,a}^{(1)}$ satisfies \eqref{eq:estremaindterm}, $\sigma_{\eps,a}^{(1)}$ satisfies \eqref{eq:stimarestoxaminus1} and thus, taking into account Remark \ref{rem:appendix1}, we easily check that $\eta_{\eps,a}^{(2)}$ satisfies \eqref{eq:estremtermexpmeanc}, for some positive constant $C$ independent of $\eps$ and $a$.

Similarly, from \eqref{eq:X_tX_thetaVW}, \eqref{eq:expXthetathetaN} we have
\begin{equation}\label{eq:meancurvpartII}
\begin{array}{lll}
\displaystyle x_a^2\frac{(X_{\eps,a})_{\theta\theta}\cdot N_{\eps,a}}{|(X_{\eps,a})_\theta|^2}&=&\displaystyle  z_a^\prime + \eps \left(x_az_a^\prime-\frac{(z_a^\prime)^3}{x_a}\right) \sin \theta + \eps^2  \sigma_{\eps,a}^{(2)}\,,
\end{array}
\end{equation}
and thus combining \eqref{eq:2xasqmeancurv}, \eqref{eq:meancurvpartI}, \eqref{eq:meancurvpartII} and using \eqref{eq:xz}, we obtain \eqref{eq:estmeanc}.
\qed


\subsection{The Jacobi operator for Delaunay tori}\label{Ss:Jacobi}

In this subsection we study the Jacobi operator corresponding to the Dealunay torus $\Sigma_{\eps,a}$, i.e., the linearization of the mean curvature operator $\mathfrak{M}$ around the parameterization $X_{\eps,a}$, given by \eqref{eq:eps-toroidal-unduloid}, with respect to normal variations. More precisely, 
\begin{equation}
\label{M-prime-eps-a}
\mathfrak{M}'_{\eps,a}(0)[\varphi]:=\frac{d}{ds}[\mathfrak{M}(X_{\eps,a}+s\varphi N_{\eps,a})]_{s=0}
\qquad(\varphi\in C^{2}(\mathbb{R}^{2}))
\end{equation}
with $N_{\eps,a}$ given by \eqref{normale-eps-a}. 
%
%

\begin{Remark}
\label{R:sign-normal-eps}
Fixing $a_{0}\in(0,\frac{1}{2})$ and taking $\eps_{0}>0$ such that $X_{\eps,a}$ is regular for all $a\in[-a_{0},a_{0}]\setminus\{0\}$ and $\eps\in(0,\eps_{0})$ for every $\varphi\in C^{2}(\mathbb{R}^{2})$ one can compute the right-hand side in \eqref{M-prime-eps-a} at every point $(t,\theta)\in\R^{2}$. However, as already seen in Remark \ref{R:sign-normal}, in order to guarantee the well-poseness of the normal vector on the whole surface defined by $X_{\eps,a}+\varphi N_{\eps,a}$ one needs a bound on $\|x_{a}^{-1}\varphi\|_{C^{0}}$ (see also the proof of Proposition \ref{prop:expmeancurvnormgraph}--(i)).
\end{Remark}

We look for an expansion of $\mathfrak{M}'_{\eps,a}(0)$ with uniform estimates with respect to the parameters $a$ and $\eps$. To this aim, it is convenient to introduce also the linear operator
\begin{equation}
\label{eq:L-eps-a}
\mathcal{L}_{\eps,a}:=2x_{a}^{2}\mathfrak{M}'_{\eps,a}(0)~\!.
\end{equation}
%
%
%
%
%
In the spirit of Proposition \ref{prop:expmeancurv}, we can decompose $\mathcal{L}_{\eps,a}$ as follows. 
%
%
\begin{Proposition}\label{prop:explinearized}
Fixing $a_{0}\in(0,\frac{1}{2})$, there exist $\eps_0>0$ and $C>0$ such that for all $a \in [-a_0,a_0]\setminus\{0\}$ and $\eps\in(0,\eps_0)$, for every $\varphi \in C^{2}(\mathbb{R}^{2})$ 
one has
\[
\mathcal{L}_{\eps,a} \varphi= \mathcal{L}_{a}\varphi + \eps \mathcal{L}_{a}^{(1)}\varphi+ \eps^2 \mathcal{L}_{\eps,a}^{(2)}\varphi
\]
where 
\begin{equation}
\label{normalizedJacop-a}
\mathcal{L}_{a}\varphi:=\Delta\varphi + 2\left(x_a^2+\frac{\gamma_a^2}{x_a^2}\right)\varphi
\end{equation} 
is the normalized Jacobi operator corresponding to $\Sigma_a$,
\begin{equation}\label{eq:expoperatorL1}
\begin{split}
&\displaystyle\mathcal{L}_{a}^{(1)}\varphi=\displaystyle \left(-2\frac{(z_a^\prime)^2}{x_a}\sin\theta\right)  \varphi_{tt} + \left[ \left(\frac{(z_a^\prime)^2x^\prime_a}{x_a^2}-4 z_a^\prime x_a^{\prime} \right)\sin\theta\right] \varphi_t+ \left(\frac{(z_a^\prime)^2}{x_a} \cos\theta\right) \varphi_{\theta}\\[3pt]
&\quad\displaystyle +2\left\{\left[\left(x_a+\frac{\gamma_a}{x_a}\right)\left(2x_a^2+z_a^\prime -\frac{(z_a^\prime)^3}{x_a}\right) - 4\frac{(z_a^\prime)^2}{x_a} \left(x_a+\frac{\gamma_a}{x_a}\right)^2+\frac{(x_a^\prime)^2}{x_a}\left(x_a-\frac{\gamma_a}{x_a}\right)^2\right]\sin\theta\right\} \varphi\,,
\end{split}
\end{equation}
\begin{equation}\label{eq:expoperatorL2}
\mathcal{L}_{\eps,a}^{(2)}\varphi=\alpha_{\eps,a}^{(1)}\varphi_{tt}  + \alpha_{\eps,a}^{(2)}\varphi_{t} + \alpha_{\eps,a}^{(3)}\varphi_{\theta} + \alpha_{\eps,a}^{(4)} \varphi,
\end{equation}
and the coefficients $\alpha_{\eps,a}^{(i)}=\alpha_{\eps,a}^{(i)}(t,\theta)$, $i=1,\ldots,4$ are doubly periodic with respect to the rectangle $K_{a}=[-\tau_{a},\tau_{a}]\times[-\pi,\pi]$ and satisfy
\begin{equation}\label{eq:unifboundexpLin}
\|x_{a}^{-1}\alpha_{\eps,a}^{(i)} \|_{C^0(K_{a})}  \leq C  \ \ \  i=1,\ldots,4.
\end{equation}
\end{Proposition}
The proof of Proposition \ref{prop:explinearized} needs some estimates which are stated in the next lemma. 

\begin{Lemma}\label{lem:expNtNthetasquared}
Fixing $a_{0}\in(0,\frac{1}{2})$, there exist $\eps_0>0$ and $C>0$ such that for all $a \in [-a_0,a_0]\setminus\{0\}$ and $\eps\in(0,\eps_0)$ one has
\begin{eqnarray}
&&\hspace{-13mm}|(N_{\eps,a})_t|^2=\left(x_a+\frac{\gamma_a}{x_a}\right)^2+ \eps \left[{2\left(x_a+\frac{\gamma_a}{x_a}\right)\left(2(x_a^\prime)^2 - 2 (z_a^\prime)^2 +\frac{(z_a^\prime)^3}{x_a^2}+z_a^\prime \right)} \right] \sin \theta + \eps^2 \psi^{(1)}_{\eps,a}\label{eq:norm_normal_t_eps-toroidal}\\[2pt]
&&\hspace{-13mm}|(N_{\eps,a})_\theta|^2=\left(x_a-\frac{\gamma_a}{x_a}\right)^2\ + \eps \left(2 \frac{ (z^\prime_{a})^2(x_a^\prime)^2}{ x_{a}^3} \sin \theta\right)  + \eps^2 \psi^{(2)}_{\eps,a}\label{eq:norm_normal_theta_eps-toroidal}
\end{eqnarray}
where $\psi^{(i)}_{\eps,a}=\psi^{(i)}_{\eps,a}(t,\theta)$ ($i=1,2$) is doubly periodic with respect to $K_{a}=[-\tau_{a},\tau_{a}]\times[-\pi,\pi]$ and satisfies
\begin{equation}\label{eq:estimaterempsi}
\|x_a^{-1}\psi^{(i)}_{\eps,a} \|_{C^0(K_{a})}  \leq C\,.
\end{equation}
\end{Lemma}
\Proof
Let us prove \eqref{eq:norm_normal_t_eps-toroidal}. Differentiating \eqref{eq:normRVW} with respect to $t$ one has 
\begin{equation}\label{eq:Nt}
(N_{\eps,a})_t=\eps z'_{a}Q_{\eps z_a}\frac{V_{\eps,a}\wedge W_{\eps,a}}{|V_{\eps,a}\wedge W_{\eps,a}|}+R_{\eps z_a}\left(\frac{V_{\eps,a}\wedge W_{\eps,a}}{|V_{\eps,a}\wedge W_{\eps,a}|}\right)_t
\end{equation}
where $Q_\sigma$ and $R_{\sigma}$ are the matrices defined in \eqref{eq:rotation-matrix}. Then,  applying \eqref{eq:proppartialRsigma}, we infer that
\begin{equation}\label{eq:normNtsquared}
|(N_{\eps,a})_t|^2=\eps^2(z_a^\prime)^2\frac{|\widecheck{(V_{\eps,a}\wedge W_{\eps,a}})|^2}{|V_{\eps,a}\wedge W_{\eps,a}|^2} + 2\eps z_a^\prime \frac{[\textbf{e}_1\cdot [(V_{\eps,a}\wedge W_{\eps,a})\wedge(V_{\eps,a}\wedge W_{\eps,a})_t]]}{|V_{\eps,a}\wedge W_{\eps,a}|^2}+  \left|\left(\frac{V_{\eps,a}\wedge W_{\eps,a}}{|V_{\eps,a}\wedge W_{\eps,a}|}\right)_t\right|^2.
\end{equation}
By definition and thanks to Remark \ref{rem:appendix1} we immediately recognize that the first term of the right-hand side of \eqref{eq:normNtsquared} is of the form $\eps^2\eta_{\eps,a}^{(1)}$,
for some doubly periodic function $\eta_{\eps,a}^{(1)}=\eta_{\eps,a}^{(1)}(t,\theta)$ and, fixing $a_{0}\in(0,\frac{1}{2})$, there exist $\eps_0>0$ and $C>0$ such that $\eta_{\eps,a}^{(1)}$ satisfies \eqref{eq:estimaterempsi} for all $a\in[-a_0,a_0]\setminus\{0\}$ and $\eps\in(0,\eps_0)$.
Moreover, since
\[
{V_{\eps,a}\wedge W_{\eps,a}}=\left[\begin{array}{c} -{z_a^\prime}{x_a} \cos\theta\\[6pt] -{z_a^\prime}{x_a} \sin\theta \\[6pt] {x_a^\prime}{x_a}\end{array}\right] + \eps \left[\begin{array}{c} - z_a^\prime x_a^2 \cos\theta \sin\theta\\[6pt] -z_a^\prime x_a^2\sin^2\theta \\[6pt] 0 \end{array}\right]\,,
\]
then, by direct computation and taking into account \eqref{eq:xz}, we infer that
\[
({V_{\eps,a}\wedge W_{\eps,a}})_t=\left[\begin{array}{c} -(2x_a^2+z_a^\prime)x_a^\prime \cos\theta\\[6pt] -(2x_a^2+z_a^\prime)x_a^\prime \sin\theta \\[6pt]x_a^2+(x_a^\prime)^2 - 2z_a^\prime x_a^2\end{array}\right] + \eps \left[\begin{array}{c} - 2x_a x_a^\prime(x_a^2+z_a^\prime) \cos\theta \sin\theta\\[6pt] -2x_a x_a^\prime(x_a^2+z_a^\prime) \sin^2\theta \\[6pt] 0 \end{array}\right]\,.
\]
Hence, taking into account \eqref{eq:conformality}, \eqref{eq:normVwedgeW} and \eqref{eq:elemexpsqrt}, we easily check that
\begin{equation}\label{eq:normNtsquaredpiece2}
2\eps z_a^\prime\frac{\textbf{e}_1\cdot [(V_{\eps,a}\wedge W_{\eps,a})\wedge(V_{\eps,a}\wedge W_{\eps,a})_t]}{|V_{\eps,a}\wedge W_{\eps,a}|^2}=2\eps {z_a^\prime}\left(x_a+\frac{\gamma_a}{x_a}\right) \sin\theta+\eps^2 \eta_{\eps,a}^{(2)}\,,
\end{equation}
where $\eta_{\eps,a}^{(2)}$ enjoys the same properties of $\eta_{\eps,a}^{(1)}$. Finally, differentiating \eqref{eq:VwedgeW_eps-toroidal} with respect to $t$, using \eqref{eq:xz}, \eqref{eq:conformality},  taking into account that
\[
\left(\frac{z_a^\prime}{x_a}\right)^\prime=x_a^\prime\left(1+\frac{\gamma_a}{x_a^2}\right), \ \ \left(\frac{x_a^\prime}{x_a}\right)^\prime=-z_a^\prime\left(1+\frac{\gamma_a}{x_a^2}\right), \ \ \left[\left(\frac{z_a^\prime}{x_a}\right)^2\right]^\prime=2\frac{z_a^\prime x_a^\prime}{x_a}\left(1+\frac{\gamma_a}{x_a^2}\right)\,,
\]
we deduce that
\begin{equation}\label{eq:VwedgeWnormalizedt}
\begin{split}
\displaystyle\left(\frac{V_{\eps,a}\wedge W_{\eps,a}}{|V_{\eps,a}\wedge W_{\eps,a}|}\right)_t&=\displaystyle\left[\begin{array}{c} -\left(1+\frac{\gamma_a}{x_a^2}\right)x_a^\prime \cos\theta\\[3pt] -\left(1+\frac{\gamma_a}{x_a^2}\right)x_a^\prime \sin\theta \\[3pt]-\left(1+\frac{\gamma_a}{x_a^2}\right)z_a^\prime\end{array}\right] + \eps \left[\begin{array}{c}  \frac{x_a^\prime}{x_a}\left(-2(x_a^\prime)^2 + 2(z_a^\prime)^2(1+\frac{\gamma_a}{x_a^2})\right) \cos\theta \sin\theta\\[3pt]  \frac{x_a^\prime}{x_a}\left(-2(x_a^\prime)^2 + 2(z_a^\prime)^2(1+\frac{\gamma_a}{x_a^2})\right) \sin^2\theta \\[3pt] \frac{z_a^\prime}{x_a}\left(+2(z_a^\prime)^2-2(x_a^\prime)^2 (1+\frac{\gamma_a}{x_a^2})-z_a^\prime\right)\sin\theta  \end{array}\right]\\
&\quad+\eps^2(\Psi_{\eps,a})_t
\end{split}
\end{equation}
and thus
\begin{equation}\label{eq:VwedgeWtsquarednorm}
\left|\left(\frac{V_{\eps,a}\wedge W_{\eps,a}}{|V_{\eps,a}\wedge W_{\eps,a}|}\right)_t\right|^2=\displaystyle \left(x_a+\frac{\gamma_a}{x_a}\right)^2+2\eps {\left(x_a+\frac{\gamma_a}{x_a}\right)\left(2(x_a^\prime)^2  - 2(z_a^\prime)^2 +\frac{(z_a^\prime)^3}{x_a^2}\right)} \sin \theta + \eps^2 \eta_{\eps,a}^{(3)}
\end{equation}
where $\eta_{\eps,a}^{(3)}$ enjoys the same properties of $\eta_{\eps,a}^{(1)}$, thanks to \eqref{eq:stimaremcompPsi}. Then, putting together \eqref{eq:normNtsquared}, \eqref{eq:normNtsquaredpiece2} and \eqref{eq:VwedgeWtsquarednorm}
we readily obtain \eqref{eq:norm_normal_t_eps-toroidal} with $\psi_{\eps,a}^{(1)}:=\eta_{\eps,a}^{(1)}+\eta_{\eps,a}^{(2)}+\eta_{\eps,a}^{(3)}$, which is doubly periodic and satisfies \eqref{eq:estimaterempsi}.
\medskip


As far as concerns \eqref{eq:norm_normal_theta_eps-toroidal}, we argue in a similar way. Indeed, since 
\begin{equation}\label{eq:Ntheta}
(N_{\eps,a})_\theta=R_{\eps z_a}\left(\frac{V_{\eps,a}\wedge W_{\eps,a}}{|V_{\eps,a}\wedge W_{\eps,a}|}\right)_\theta
\end{equation}
and since
\begin{equation}\label{eq:VwedgeWnormalizedtheta}
 \left(\frac{V_{\eps,a}\wedge W_{\eps,a}}{|V_{\eps,a}\wedge W_{\eps,a}|}\right)_\theta=\left[\begin{array}{c}\frac{z^\prime_{a}}{x_a}\sin\theta\\[8pt] -\frac{z^\prime_{a}}{x_a} \cos\theta \\[8pt] 0\end{array}\right] + \eps \left[\begin{array}{c} \left(- z^\prime_{a} + \frac{ (z^\prime_{a})^3}{ x^2_{a}}\right) \cos(2\theta) \\[4pt]  \left(- z^\prime_{a} +  \frac{ (z^\prime_{a})^3}{ x^2_{a}}\right) \sin(2\theta)\\[4pt] -\left(\frac{ z^\prime_{a}}{ x_{a}}\right)^2x^\prime_a \cos\theta\end{array}\right] + \eps^2 (\Psi_{\eps,a})_\theta,
 \end{equation}
then exploiting \eqref{eq:propRsigma} we easily conclude.
\qed

\noindent
\textit{Proof of Proposition \ref{prop:explinearized}.}
For every $a\in[-\frac{1}{2},\infty)\setminus\{0\}$ and $\eps>0$ small, 
fixing $\varphi\in C^{2}(\mathbb{R}^{2})$, according to Proposition \ref{Prop:firstvariationmeancurv} in Appendix, we have that
\begin{equation}\label{eq:linearized}
\begin{split}
&\displaystyle \frac{d}{ds}[\mathfrak{M}(X_{\eps,a}+s\varphi N_{\eps,a})]_{s=0}=\displaystyle  \left(\frac{1}{2|(X_{\eps,a})_t|^2}\right) \varphi_{tt} + \left(\frac{1}{2|(X_{\eps,a})_\theta|^2}\right) \varphi_{\theta\theta} \\
&\displaystyle + \left(\frac{1}{4|(X_{\eps,a})_t|^2|(X_{\eps,a})_\theta|^2} \frac{d|(X_{\eps,a})_\theta|^2}{dt}- \frac{1}{4|(X_{\eps,a})_t|^4} \frac{d|(X_{\eps,a})_t|^2}{dt}\right) \varphi_{t} \\
&\displaystyle+  \left(\frac{1}{4|(X_{\eps,a})_t|^2|(X_{\eps,a})_\theta|^2} \frac{d|(X_{\eps,a})_t|^2}{d\theta}- \frac{1}{4|(X_{\eps,a})_\theta|^4} \frac{d|(X_{\eps,a})_\theta|^2}{d\theta}\right) \varphi_{\theta}\\
&\displaystyle +  \left(\frac{2((X_{\eps,a})_{t\theta}\cdot N_{\eps,a})^2}{|(X_{\eps,a})_t|^2|(X_{\eps,a})_\theta|^2} - \frac{|(N_{\eps,a})_\theta|^2}{2|(X_{\eps,a})_\theta|^2} - \frac{|(N_{\eps,a})_t|^2}{2|(X_{\eps,a})_t|^2} + \frac{((X_{\eps,a})_{tt}\cdot N_{\eps,a})^2}{|(X_{\eps,a})_t|^4} + \frac{((X_{\eps,a})_{\theta\theta}\cdot N_{\eps,a})^2}{|(X_{\eps,a})_\theta|^4}  \right) \varphi.
\end{split}
\end{equation}
Let us analyze the coefficients of \eqref{eq:linearized}. Firstly, from \eqref{eq:X_tX_thetaVW} and arguing as in the proof of Lemma \ref{lem:expsecderXN}, fixing $a_{0}\in(0,\frac{1}{2})$, one can find $\eps_0>0$ such that for any $a\in[-a_0,a_0]\setminus\{0\}$ and $\eps\in(0,\eps_0)$ one has
 \begin{equation}\label{eq:expinvsqnormXt}
\frac{1}{|(X_{\eps,a})_t|^2} = x_a^{-2} \left(1-2\eps \frac{(z_a^\prime)^2}{x_a}\sin\theta+\eps^2 \eta_{\eps,a}^{(1)}\right),
\end{equation}
where $\eta_{\eps,a}^{(1)}$ satisfies \eqref{eq:estremaindterm}, for some $C>0$ independent of $a$ and $\eps$.

In addition, differentiating \eqref{eq:X_tX_thetaVW} and recalling that $z_a^{\prime\prime}=2x_a x_a^\prime$, we get that
\begin{equation}\label{eq:dersquarednormttheta}
\begin{split}
&\displaystyle \frac{d|(X_{\eps,a})_t|^2}{dt}= 2x_a x_a^\prime + \eps \left(2x_a^\prime (z_a^\prime)^2+8x_a^2 x_a^\prime z_a^\prime \right) \sin \theta + \eps^2 \left(4z_a^\prime x_a^3x_a^\prime+ 2 (z_a^\prime)^2x_a x_a^\prime \right) \sin^2\theta\,,
\\
&\displaystyle  \frac{d|(X_{\eps,a})_t|^2}{d\theta}= \displaystyle \eps \left(2x_a (z_a^\prime)^2\cos \theta\right) + \eps^2 \left(2\left(z_a^\prime x_a\right)^2\cos\theta \sin\theta\right)\,,
\\
&\displaystyle  
\frac{d|(X_{\eps,a})_\theta|^2}{dt}= \displaystyle  2x_a x_a^\prime\,,
\\  
&\displaystyle  
\frac{d|(X_{\eps,a})_\theta|^2}{d\theta}= 0\,.
\end{split}
\end{equation}


Putting together \eqref{M-prime-eps-a}--\eqref{eq:L-eps-a} and \eqref{eq:linearized}--\eqref{eq:dersquarednormttheta}, exploiting Lemma \ref{lem:expsecderXN}, Lemma \ref{lem:expNtNthetasquared} and taking into account \eqref{eq:xz}, 
for all $a\in[-a_0,a_0]\setminus\{0\}$ and $\eps\in(0,\eps_0)$ it holds
\[
\begin{split}
\mathfrak{L}_{\eps,a}\varphi=&\displaystyle  \left(1-2\eps \frac{(z_a^\prime)^2}{x_a}\sin\theta+\eps^2 \eta_{\eps,a}^{(1)}\right) \varphi_{tt} + \varphi_{\theta\theta}\\
&\displaystyle +  \left\{\left(1-2\eps \frac{(z_a^\prime)^2}{x_a}\sin\theta+\eps^2 \eta_{\eps,a}^{(1)}\right)\frac{x_a^\prime}{x_a}-\left(1-2\eps \frac{(z_a^\prime)^2}{x_a}\sin\theta+\eps^2 \eta_{\eps,a}^{(1)}\right)^2 \right.\\
&\displaystyle\ \ \ \ \ \ \ \cdot \left.\left[\frac{x_a^\prime}{x_a}+\eps\left( x_a^\prime \left(\frac{z_a^\prime}{x_a}\right)^2+4 z_a^\prime x_a^\prime\right)\sin\theta+\eps^2\left(2z_a^\prime x_a x_a^\prime + \left(\frac{z_a^\prime}{x_a}\right)^2 x_a x_a^\prime\right)\sin^2\theta\right]\right\} \varphi_{t}\\
&\displaystyle+ \left\{\left[\eps \frac{(z_a^\prime)^2}{x_a}\cos\theta + \eps^2 (z_a^\prime)^2 \cos\theta\sin\theta\right] \left(1-2\eps \frac{(z_a^\prime)^2}{x_a}\sin\theta+\eps^2 \eta_{\eps,a}^{(1)}\right)\right\}  \varphi_{\theta}\\
&\displaystyle + \left\{4\left(\eps \frac{x_a^\prime z_a^\prime}{x_a} \cos\theta+\eps^2\frac{\sigma_{\eps,a}^{(3)}}{x_a}\right)^2 \left(1-2\eps \frac{(z_a^\prime)^2}{x_a}\sin\theta+\eps^2 \eta_{\eps,a}^{(1)}\right)\right.\\
&\displaystyle\left.\ \ \ \ \ \ -\left[\left(x_a-\frac{\gamma_a}{x_a}\right)^2\ + \eps \left(2 \frac{ (z^\prime_{a})^2(x_a^\prime)^2}{ x_{a}^3} \sin \theta\right)  + \eps^2 \psi^{(2)}_{\eps,a}\right]\right.\\
& \displaystyle\left.\ \ \ \ \ \ - \left[\left(x_a+\frac{\gamma_a}{x_a}\right)^2+ \eps \left[{2\left(x_a+\frac{\gamma_a}{x_a}\right)\left(2(x_a^\prime)^2 - 2 (z_a^\prime)^2 +\frac{(z_a^\prime)^3}{x_a^2}+z_a^\prime \right)} \right] \sin \theta + \eps^2 \psi^{(1)}_{\eps,a}\right] \right.\\
& \displaystyle\left.\ \ \ \ \ \ \ \cdot\left(1-2\eps \frac{(z_a^\prime)^2}{x_a}\sin\theta+\eps^2 \eta_{\eps,a}^{(1)}\right)+2\left[ \left(x_a+\frac{\gamma_a}{x_a}\right) + \eps \left(2x_a^2+z_a^\prime-2(z_a^\prime)^2\right) \sin\theta + \eps^2\frac{ \sigma_{\eps,a}^{(1)}}{x_a}\right]^2 \cdot \right.\\
&\displaystyle\left. \ \ \ \ \ \  \cdot\left(1-2\eps \frac{(z_a^\prime)^2}{x_a}\sin\theta+\eps^2 \eta_{\eps,a}^{(1)}\right)^2 + 2 \left[\left(x_a-\frac{\gamma_a}{x_a}\right) + \eps \left(z_a^\prime-\frac{(z_a^\prime)^3}{x_a^2}\right) \sin \theta + \eps^2  \frac{\sigma_{\eps,a}^{(2)}}{x_a}\right]^2\right\} \varphi.
\end{split}
\]
Now, since $ \eta_{\eps,a}^{(j)}$, $\psi_{\eps,a}^{(j)}$  satisfy, respectively, \eqref{eq:estremaindterm},  \eqref{eq:estimaterempsi}, for $j=1,2$, and $\sigma_{\eps,a}^{(i)}$  satisfies \eqref{eq:stimarestoxaminus1}, for $i=1,2,3$, then, taking into account Remarks \ref{rem:appendix1}  and \ref{rem:appendix2}, \eqref{eq:xz}, \eqref{eq:conformality}, and noticing that
\[
\frac{(z^\prime_{a})^2(x_a^\prime)^2}{ x_{a}^3}=\left(x_a-\frac{\gamma_a}{x_a}\right)^2\frac{(x_a^\prime)^2}{x_a}\,,\qquad \left(x_a-\frac{\gamma_a}{x_a}\right)\left(z_a^\prime-\frac{(z_a^\prime)^3}{x_a^2}\right)=\left(x_a-\frac{\gamma_a}{x_a}\right)^2\frac{(x_a^\prime)^2}{x_a}\,,
\]
after some elementary computations we obtain
\[
\begin{split}
\mathfrak{L}_{\eps,a}\varphi=&\displaystyle  ~\Delta \varphi +  \left(-2\eps \frac{(z_a^\prime)^2}{x_a}\sin\theta\right) \varphi_{tt} + \eps^2\alpha_{\eps,a}^{(1)} \varphi_{tt} + \eps\left[\left( \frac{(z_a^\prime)^2}{x_a^{2}}x_a^\prime-4z_a^\prime x_a^\prime\right)\sin\theta\right] \varphi_{t} +\eps^2\alpha_{\eps,a}^{(2)} \varphi_{t}\\
&\displaystyle + \eps\left(\frac{(z_a^\prime)^2}{x_a}\cos\theta\right)\varphi_{\theta} +\eps^2\alpha_{\eps,a}^{(3)} \varphi_{\theta}+ \left[\left(x_a-\frac{\gamma_a}{x_a}\right)^2+\left(x_a+\frac{\gamma_a}{x_a}\right)^2 \right] \varphi\\
&\displaystyle  +2\eps \left\{\left[-\left(x_a-\frac{\gamma_a}{x_a}\right)^2\frac{(x_a^\prime)^2}{ x_{a}}  -\left(x_a+\frac{\gamma_a}{x_a}\right)\left(2(x_a^\prime)^2 - 2 (z_a^\prime)^2 +z_a^\prime+\frac{(z_a^\prime)^3}{x_a^2}\right) \right.\right.\\
&\displaystyle\left. \left.\ \ \ \ \ -4\left(x_a+\frac{\gamma_a}{x_a}\right)^2\frac{(z_a^\prime)^2}{x_a}+\left(x_a+\frac{\gamma_a}{x_a}\right)\left(4(x_a^\prime)^2  +2z_a^\prime \right) + 2\left(x_a-\frac{\gamma_a}{x_a}\right)^2\frac{(x_a^\prime)^2}{x_a}\right]\sin\theta\right\}\varphi\\
&\displaystyle +\eps^2\alpha_{\eps,a}^{(4)}\varphi,
\end{split}
\]
where $\alpha_{\eps,a}^{(k)}=\alpha_{\eps,a}^{(k)}(t,\theta)$ is doubly periodic and satisfies \eqref{eq:unifboundexpLin}, $k=1,\ldots,4$. Finally, after easy algebraic manipulations the assertion follows.
\qed

We complete this section with further estimates, of the same kind to those ones stated in Lemma \ref{lem:expsecderXN}, which will enter in the proof of Propositions \ref{prop:expmeancurvnormgraph} and \ref{prop:secvarexpansion}.

\begin{Lemma}\label{lem:expNsecX}
Fixing $a_{0}\in(0,\frac{1}{2})$, there exist $\eps_0>0$ and $C>0$ such that for all $a \in [-a_0,a_0]\setminus\{0\}$ and $\eps\in(0,\eps_0)$ one has
\begin{eqnarray}
&&(\Ne)_{t}\cdot (\Ne)_\theta= \eps \psi_{\eps,a}^{(1)}\,,\label{eq:0lem:expNsecX}\\[3pt]
&&(\Ne)_{tt}\cdot\Xet=-\frac{x_a^\prime z_a^\prime}{x_a} + \eps \psi_{\eps,a}^{(2)}\,, \ \ \ \ \ (\Ne)_{tt}\cdot\Xeth= \eps \psi_{\eps,a}^{(3)}\,,\label{eq:2lem:expNsecX}\\
&&(\Ne)_{\theta\theta}\cdot\Xet=\frac{x_a^\prime z_a^\prime}{x_a} + \eps \psi_{\eps,a}^{(4)}\,,\ \ \ \ \ \ \ (\Ne)_{\theta\theta}\cdot\Xeth= \eps \psi_{\eps,a}^{(5)}\label{eq:4lem:expNsecX}
\end{eqnarray}
where $\psi^{(i)}_{\eps,a}=\psi^{(i)}_{\eps,a}(t,\theta)$, for $i=1,\ldots,5$, is doubly periodic with respect to $K_{a}=[-\tau_{a},\tau_{a}]\times[-\pi,\pi]$ and satisfies
\begin{equation}\label{eq:estimaterempsi2}
\|x_a^{-1}\psi^{(i)}_{\eps,a} \|_{C^0(K_{a})}  \leq C\,.
\end{equation}
\end{Lemma}
\Proof
The proof is similar to that of Lemma \ref{lem:expNtNthetasquared}, so we will omit some details.
As far as concerns \eqref{eq:0lem:expNsecX}, in view of \eqref{eq:Nt}, \eqref{eq:Ntheta} and exploiting \eqref{eq:propRsigma}, \eqref{eq:proppartialRsigma} we get that
\[
(\Ne)_{t}\cdot (\Ne)_\theta=\eps z_a^\prime \textbf{e}_1 \cdot\left[\left(\frac{V_{\eps,a}\wedge W_{\eps,a}}{|V_{\eps,a}\wedge W_{\eps,a}|}\right)\wedge\left(\frac{V_{\eps,a}\wedge W_{\eps,a}}{|V_{\eps,a}\wedge W_{\eps,a}|}\right)_\theta\right] + \left(\frac{V_{\eps,a}\wedge W_{\eps,a}}{|V_{\eps,a}\wedge W_{\eps,a}|}\right)_t \cdot \left(\frac{V_{\eps,a}\wedge W_{\eps,a}}{|V_{\eps,a}\wedge W_{\eps,a}|}\right)_\theta
\]
and the assertion follows from \eqref{eq:VwedgeWnormalizedt} and \eqref{eq:VwedgeWnormalizedtheta}.
\medskip

For \eqref{eq:2lem:expNsecX}, differentiating \eqref{eq:Nt}, one has
\begin{eqnarray*}
(\Ne)_{tt}&=&\eps z_a^{\prime\prime} Q_{\eps z_a}\frac{V_{\eps,a}\wedge W_{\eps,a}}{|V_{\eps,a}\wedge W_{\eps,a}|} + (\eps z_a^{\prime})^2 R_{\eps z_a}\frac{V_{\eps,a}\wedge W_{\eps,a}}{|V_{\eps,a}\wedge W_{\eps,a}|}+ 2\eps z_a^{\prime} Q_{\eps z_a}\left(\frac{V_{\eps,a}\wedge W_{\eps,a}}{|V_{\eps,a}\wedge W_{\eps,a}|}\right)_t\\
&&+R_{\eps z_a}\left(\frac{V_{\eps,a}\wedge W_{\eps,a}}{|V_{\eps,a}\wedge W_{\eps,a}|}\right)_{tt}
\end{eqnarray*}
where $R_{\sigma}$ and $Q_{\sigma}$ are the matrices defined in \eqref{eq:rotation-matrix}. Then, exploiting \eqref{eq:relXtXthetaR}, taking into account \eqref{eq:propRsigma}, \eqref{eq:proppartialRsigma} and Remark \ref{rem:appendix1}, we easily infer that
\[
(\Ne)_{tt}\cdot \Xet=\left(\frac{V_{\eps,a}\wedge W_{\eps,a}}{|V_{\eps,a}\wedge W_{\eps,a}|}\right)_{tt}\cdot V_{\eps,a} + \eps \psi_{\eps,a}^{(2)}\,,\ \ \ (\Ne)_{tt}\cdot \Xeth=\left(\frac{V_{\eps,a}\wedge W_{\eps,a}}{|V_{\eps,a}\wedge W_{\eps,a}|}\right)_{tt}\cdot W_{\eps,a} + \eps \psi_{\eps,a}^{(3)}
\]
for some doubly periodic functions $\psi_{\eps,a}^{(2)}=\psi_{\eps,a}^{(2)}(t,\theta)$, $\psi_{\eps,a}^{(3)}=\psi_{\eps,a}^{(3)}(t,\theta)$ satisfying \eqref{eq:estimaterempsi2}.

Now, differentiating \eqref{eq:VwedgeWnormalizedt} with respect to $t$ one has
\[
\left(\frac{V_{\eps,a}\wedge W_{\eps,a}}{|V_{\eps,a}\wedge W_{\eps,a}|}\right)_{tt}=\left[\begin{array}{c} \left( 2 \gamma_a \frac{(x_a^\prime)^2}{x_a^3} -(1+\frac{\gamma_a}{x_a^2})x_a^{\prime\prime}\right) \cos\theta\\[6pt] 
\left(2 \gamma_a \frac{(x_a^\prime)^2}{x_a^3} -(1+\frac{\gamma_a}{x_a^2})x_a^{\prime\prime}\right) \sin\theta \\[6pt]
2 \gamma_a \frac{x_a^\prime z_a^\prime}{x_a^3} -(1+\frac{\gamma_a}{x_a^2})z_a^{\prime\prime}\end{array}\right] + \eps\Psi^{(1)}_{\eps,a}\,,
\]
where $\Psi_{\eps,a}^{(1)}$ is a doubly periodic vector-valued function and, in view of \eqref{eq:xz}, \eqref{eq:stimaremcompPsi} and thanks to Remark \ref{rem:appendix1}, each of its components satisfies \eqref{eq:estimaterempsi2}.

Finally, exploiting \eqref{eq:conformality} (which implies that $x_ax_a^\prime=x_a^\prime x_a^{\prime\prime}+z_a^\prime z_a^{\prime\prime}$) and \eqref{eq:defvectVW}, after some standard computations we obtain \eqref{eq:2lem:expNsecX}.
\medskip

The proof \eqref{eq:4lem:expNsecX} is analogous. We limit to observe that $(\Ne)_{\theta\theta}=R_{\eps z_a}\left(\frac{V_{\eps,a}\wedge W_{\eps,a}}{|V_{\eps,a}\wedge W_{\eps,a}|}\right)_{\theta\theta}$ and 
\[
\left(\frac{V_{\eps,a}\wedge W_{\eps,a}}{|V_{\eps,a}\wedge W_{\eps,a}|}\right)_{\theta\theta}=\left[\begin{array}{c}\frac{z^\prime_{a}}{x_a}\cos\theta\\[10pt] \frac{z^\prime_{a}}{x_a} \sin\theta \\[10pt] 0\end{array}\right] + \eps \left[\begin{array}{c} -2\left(- z^\prime_{a} + \frac{ (z^\prime_{a})^3}{ x^2_{a}}\right) \sin(2\theta) \\[6pt] 2 \left(- z^\prime_{a} +  \frac{ (z^\prime_{a})^3}{ x^2_{a}}\right) \cos(2\theta)\\[6pt] \left(\frac{ z^\prime_{a}}{ x_{a}}\right)^2x^\prime_a \sin\theta\end{array}\right] + \eps^2 (\Psi_{\eps,a})_{\theta\theta}\,,
\]
 as it follows by differentiating \eqref{eq:VwedgeWnormalizedtheta}. Then, from \eqref{eq:relXtXthetaR} and exploiting \eqref{eq:propRsigma} and taking into account \eqref{eq:stimaremcompPsi}, we easily conclude. 
\qed
\begin{Remark}
Since $\Ne \cdot \Net=\Ne\cdot \Neth=0$ 
 one has $\Net=A_{1} \widehat\Xet + A_{2} \widehat\Xeth$, for some $A_{1},A_2:\R^2\to\R$. Then, exploiting \eqref{eq:propNetwedgeNeth} one gets
\begin{equation}\label{eq:exprNtalt}
\Net= - \frac{\Xett \cdot \Ne }{ |\Xet|}\widehat \Xet - \frac{\Xetth \cdot \Ne}{ |\Xeth|}\widehat \Xeth,.
\end{equation}
Similarly one has
\begin{equation}\label{eq:exprNthetaalt}
\Neth= - \frac{\Xetth \cdot \Ne }{ |\Xet|}\widehat \Xet - \frac{\Xethth \cdot \Ne}{ |\Xeth|}\widehat \Xeth\,.
\end{equation}
In particular, using \eqref{eq:exprNtalt}, \eqref{eq:exprNthetaalt}, and Lemma \ref{lem:expsecderXN}, we can prove Lemma \ref{lem:expNtNthetasquared} and \eqref{eq:0lem:expNsecX} in an alternative way.
\end{Remark}

\section{Some expansions around normal graphs over Delaunay tori}
\label{S:expansions-normal-graphs}

Here we consider (generalized) approximate Delaunay tori, i.e., parametric surfaces of the form $X_{\eps,a}+\varphi N_{\eps,a}$, where $X_{\eps,a}$ and $N_{\eps,a}$ are given by  \eqref{eq:eps-toroidal-unduloid} and \eqref{normale-eps-a}, respectively, $\eps$ is a small positive parameter,  $\varphi\in C^{2}(\mathbb{R}/_{2n\tau_{a}}\times\mathbb{R}/_{2\pi})$ suitably small with respect to the weighted norm \eqref{eq:weighted-C2-norm} and 
\begin{equation}
\label{n-epsilon}
n=\frac{\pi}{\eps h_{a}}\,.
\end{equation}
We point out that, in view of \eqref{eq:acca-a}, taking $a$ in a bounded set, one has that $\eps\to 0$ if and only if $n\to\infty$. Moreover, when ${\pi}/(\eps h_{a})\in\mathbb{N}$ the function $X_{\eps,a}+\varphi N_{\eps,a}$ is a parameterization of a complete, compact surface of genus one. In fact, the results stated in this section hold true for all $\eps>0$ small enough.
\medskip
 
Fixing $a_{0}\in(0,\frac12)$, $\bta>0$ and $R>0$, for $a\in [-a_{0},a_{0}]\setminus\{0\}$ and $\eps>0$, let us set
\begin{equation}
\label{eq:U}
\mathscr{B}_{\eps,a}(\bta,R):=\{\varphi\in C^{2}(\mathbb{R}/_{2n\tau_{a}}\times\mathbb{R}/_{2\pi})~|~\|\varphi\|_{C^{2}(K_{n,a};x_{a}^{-1})}\le R\eps^{\bta}\}\,.
\end{equation}
The goal of this section is twofold: provide an expansion for $ \mathfrak{M}(X_{\eps,a}+\varphi N_{\eps,a})$ with respect to small $\eps$, uniformly with respect to $a\in[-a_{0},a_{0}]\setminus\{0\}$ and $\varphi\in \mathscr{B}_{\eps,a}(\bta,R)$. Moreover, fixing a function $H\colon\mathbb{R}^{3}\to\mathbb{R}$ satisfying the assumptions stated in Theorems \ref{mainteo} or \ref{mainteo2}, we aim to provide an expansion also for $H(X_{\eps,a}+\varphi N_{\eps,a})$, under the same conditions on $\eps$, $a$ and $\varphi$.
\medskip

As done in the previous sections, it is convenient to multiply by $2x_a^2$ and to write 
\begin{eqnarray}
&&\displaystyle 2x_a^2 \mathfrak{M}(X_{\varepsilon,a}+\varphi N_{\varepsilon,a}) =\displaystyle 2x_a^2 \mathfrak{M}(X_{\varepsilon,a})+\mathcal{L}_{\eps,a}\varphi + \mathfrak{B}_{\eps,a}^0(\varphi),\label{eq:firsttermdecom}\\[2pt]
&&\displaystyle 2x_a^2 H(X_{\varepsilon,a}+\varphi N_{\varepsilon,a})=\displaystyle 2x_a^2 H(X_{\varepsilon,a}) + 2x_a^2 (\nabla_X H (X_{\varepsilon,a})\cdot N_{\eps,a}) \varphi + \mathfrak{B}_{\eps,a}^1(\varphi), \label{eq:secondtermdecom}
\end{eqnarray}
where $\mathfrak{B}_{\eps,a}^0(\varphi)$ and $\mathfrak{B}_{\eps,a}^0(\varphi)$ are remainder terms and $\mathcal{L}_{\eps,a}$ is the normalized Jacobi operator (see \eqref{eq:L-eps-a}). We already provided expansions for $2x_a^2\mathfrak{M}(X_{\varepsilon,a})$ and $\mathcal{L}_{\eps,a}$ in Section \ref{S:Delaunay-tori}. Here we focus on the remaining terms of \eqref{eq:firsttermdecom} and \eqref{eq:secondtermdecom}.
\medskip

In the next result we prove that for all sufficiently small $\eps$, and for all $\varphi\in \mathscr{B}_{\eps,a}(\bta,R)$ the mapping $X_{\eps,a}+\varphi N_{\eps,a}$ defines a regular surface and we obtain the desired expansion for $\mathfrak{M}(\Xe+\varphi\Ne)$. As by-product, a uniform $C^0$-estimate for $\mathfrak{B}_{\eps,a}^0(\varphi)$, for $\varphi \in \mathscr{B}_{\eps,a}(\bta,R)$, will follow (see Corollary \ref{cor:stimaB0}).

\begin{Proposition}\label{prop:expmeancurvnormgraph}
Fixing $a_{0}\in(0,\frac{1}{2})$, for any $R>0$ and $\bta>0$ there exist $\eps_0>0$ such that for any $a\in[-a_0,a_0]\setminus\{0\}$, $\eps\in(0,\eps_0)$,  and for any $\varphi \in \mathscr{B}_{\eps,a}(\bta,R)$ 
\begin{itemize}
\item[(i)] $(X_{\eps,a}+\varphi N_{\eps,a})_{t}\wedge(X_{\eps,a}+\varphi N_{\eps,a})_{\theta}$ never vanishes,
\item[(ii)] there exists a function $h_{\eps,a}=h_{\eps,a}(t,\theta)\colon\mathbb{R}\times\mathbb{R}/_{2\pi}\to\mathbb{R}$ such that
\[
2x_a^2\mathfrak{M}(\Xe+\varphi\Ne)=2x_a^2\mathfrak{M}(\Xe) + \mathcal{L}_{\eps,a}\varphi + \eps^{2\bta}h_{\eps,a}\,,
\]
\begin{equation}\label{eq:boundremainderexpmean}
\|h_{\eps,a}\|_{C^0(\R\times[-\pi,\pi])}\leq C\,
\end{equation}
for some positive constant $C$ independent of $a$, $\eps$, $\varphi$. If in addition $\frac{\pi}{\eps h_{a}} \in \mathbb{N}$  then $h_{\eps,a}$ is doubly periodic with respect to $K_{n,a}=[-n\tau_{a},n\tau_{a}]\times[-\pi,\pi]$, where $n=\frac{\pi}{\eps h_{a}}$.
\end{itemize}
\end{Proposition}

\Proof
{\it (i)} Fixing $R>0$ and $\bta>0$, observe that by definition, for any given $a\in[-\frac{1}{2},\infty)\setminus\{0\}$, $\eps>0$  and for any $\varphi \in \mathscr{B}_{\eps,a}(\bta,R)$ one has
\begin{equation}\label{eq:boundweighvarphi}
|\varphi(t,\theta)|\leq R \eps^{\bta} x_a(t) \ \ \forall (t,\theta)\in \R^2,
\end{equation}
and the same bound holds for the first and second partial derivatives of $\varphi$.
 
On the other hand, by a straightforward computation (see \eqref{eq1foglio} in the appendix) and recalling \eqref{eq:2xasqmeancurv}, for any $\varphi \in C^2(\R^2,\R)$ one has 
 \begin{equation}\label{eq:XtwedgeXthetavarphi2}
\begin{split}
&\displaystyle\left[(X_{\eps,a}+\varphi \Ne)_{t}\wedge (\Xe +\varphi \Ne)_{\theta}\right]\cdot \Ne\\[3pt] 
&= \displaystyle |(\Xe)_t| |(\Xe)_\theta| \left\{ 1-  2\mathfrak{M}(\Xe)\varphi +\left[\frac{(\Xe)_{tt} \cdot \Ne}{|(\Xe)_t|^2}\frac{ (\Xe)_{\theta\theta} \cdot \Ne}{|(\Xe)_\theta|^2}-\left(\frac{(\Xe)_{t \theta} \cdot \Ne}{|(\Xe)_t| |(\Xe)_\theta|}\right)^2\right]\varphi^2\right\}.
\end{split}
\end{equation}
Now, thanks to Proposition \ref{prop:expmeancurv} and taking into account \eqref{eq:boundweighvarphi} and Remark \ref{rem:appendix1}, fixing $a_{0}\in(0,\frac{1}{2})$, one can find $\eps_0>0$ and $C>0$ such that for any  $a\in[-a_0,a_0]\setminus\{0\}$, $\eps\in(0,\eps_0)$, $\varphi \in \mathscr{B}_{\eps,a}(\bta,R)$
\begin{equation}\label{eq:estimateMXea}
|\mathfrak{M}(\Xe)\varphi|\leq C \eps^{\bta} \ \ \forall (t,\theta)\in\R^2.
\end{equation}
Similarly, from Lemma \ref{lem:expsecderXN}, \eqref{eq:expinvsqnormXt}, \eqref{eq:boundweighvarphi} and Remark \ref{rem:appendix1}, up to taking a smaller $\eps_0$ (whose choice depends only on $R$, $\bta$ and $a_{0}$), for any $a\in[-a_0,a_0]\setminus\{0\}$, $\eps\in(0,\eps_0)$,  $\varphi \in \mathscr{B}_{\eps,a}(\bta,R)$ one has
\begin{equation}\label{eq:estimatesecondpiece}
\left|\left[\frac{(\Xe)_{tt} \cdot \Ne}{|(\Xe)_t|^2}\frac{ (\Xe)_{\theta\theta} \cdot \Ne}{|(\Xe)_\theta|^2}-\left(\frac{(\Xe)_{t \theta} \cdot \Ne}{|(\Xe)_t| |(\Xe)_\theta|}\right)^2\right] \varphi^2\right|\leq C \eps^{2\bta} \ \ \forall (t,\theta)\in\R^2,
\end{equation}
for some constant $C>0$ independent of $a$, $\eps$ and $\varphi$.
\medskip

Therefore, putting together \eqref{eq:XtwedgeXthetavarphi2}--\eqref{eq:estimatesecondpiece} and recalling \eqref{eq:normVwedgeW} we find $\eps_0>0$ such that for any $a\in[-a_0,a_0]\setminus\{0\}$, $\eps\in(0,\eps_0)$, $\varphi \in \mathscr{B}_{\eps,a}(\bta,R)$ one has 
\[
[((X_{\eps}+\varphi N_{\eps,a})_t\wedge (\Xe+\varphi\Ne)_\theta)\cdot \Ne](t,\theta)\neq 0 \ \ \forall(t,\theta)\in\R^2\,,
\]
which readily implies {\it (i)}.
\medskip

{\it (ii)} Let us fix $R>0$, $\bta>0$ and $a_{0}\in(0,\frac{1}{2})$, and let $\eps_0>0$  be given according to {\it (i)}. Then for any $a\in[-a_0,a_0]\setminus\{0\}$, $\eps\in(0,\eps_0)$ and for any $\varphi\in\mathscr{B}_{\eps,a}(\bta,R)$ one has that  $\mathfrak{M}(\Xe+\varphi\Ne)$ is well defined and
\begin{equation}\label{eq:defmeancurvphi}
\mathfrak{M}(X_{\eps,a}+\varphi N_{\eps,a})=\frac{{\E}_{\varphi}{\N}_{\varphi}-2{\F}_{\varphi}{\M}_{\varphi}+{\G}_{\varphi}{\L}_{\varphi}}{2({\E}_{\varphi}{\G}_{\varphi}-{\F}_{\varphi}^{2})}
\end{equation}
where
\[
\begin{array}{lll}
{\E}_{\varphi}=|(X_{\eps,a}+\varphi{\normal}_{\eps,a})_{t}|^{2}~\!,
&
{\F}_{\varphi}=(X_{\eps,a}+\varphi{\normal}_{\eps,a})_{t}\cdot (X_{\eps,a}+\varphi{\normal}_{\eps,a})_{\theta}~\!,
&{\G}_{\varphi}=|(X_{\eps,a}+\varphi{\normal}_{\eps,a})_{\theta}|^{2}~\!,
\\
{\L}_{\varphi}=(X_{\eps,a}+\varphi{\normal}_{\eps,a})_{tt}\cdot {\normal}_{\varphi}~\!,
&{\M}_{\varphi}=(X_{\eps,a}+\varphi{\normal}_{\eps,a})_{t\theta}\cdot {\normal}_{\varphi}~\!,
&{\N}_{\varphi}=(X_{\eps,a}+\varphi{\normal}_{\eps,a})_{\theta\theta}\cdot {\normal}_{\varphi}~\!,
\end{array}
\]
\[
{\normal}_{\varphi}=\frac{(X_{\eps,a}+\varphi{\normal}_{\eps,a})_{t}\wedge (X_{\eps,a}+\varphi{\normal}_{\eps,a})_{\theta}}{\left|(X_{\eps,a}+\varphi{\normal}_{\eps,a})_{t}\wedge (X_{\eps,a}+\varphi{\normal}_{\eps,a})_{\theta}\right|}~\!.
\]

By direct computation and exploiting \eqref{identitN} one has
\begin{equation}\label{eq:EFG}
\begin{split}
\displaystyle {\E}_{\varphi}&= \displaystyle |\Xet|^2 - 2(\Xett\cdot \Ne) \varphi + \varphi_t^2+\varphi^2|\Net|^2\,,\\[3pt]
\displaystyle {\F}_{\varphi}&=\displaystyle -2 ((\Xe)_{t\theta}\cdot \Ne)\varphi + \varphi_t\varphi_\theta + (\Net\cdot\Neth) \varphi^2\,,\\[3pt]
\displaystyle{\G}_{\varphi}&= \displaystyle |\Xeth|^2 - 2(\Xethth\cdot \Ne) \varphi + \varphi_\theta^2+\varphi^2|\Neth|^2\,.
\end{split}
\end{equation}
For the remaining coefficients, we first compute an expansion of $N_\varphi$.

 Setting $\widehat\Xet:=\frac{\Xet}{|\Xet|}$, $\widehat\Xeth:=\frac{\Xeth}{|\Xeth|}$ then, after a straightforward computation (see \eqref{eq1foglio} in the appendix), one has
 \begin{small}
\begin{equation}\label{eq:XtvarwedgeXthetavartang}
\begin{split}
 &\displaystyle [(X_{\eps,a}+\varphi{\normal}_{\eps,a})_{t}\wedge (X_{\eps,a}+\varphi{\normal}_{\eps,a})_{\theta}]\cdot \widehat\Xet\\[3pt]
 &\quad=\displaystyle |\Xet||\Xeth|\left(-\frac{\varphi_t}{|\Xet|}+\frac{\Xethth\cdot\Ne}{|\Xeth|}\frac{\varphi_t\varphi}{|\Xet||\Xeth|} -\frac{\Xetth\cdot\Ne}{|\Xeth|}\frac{\varphi_\theta\varphi}{|\Xet||\Xeth|}\right),\\[6pt]
&\displaystyle  [(X_{\eps,a}+\varphi{\normal}_{\eps,a})_{t}\wedge (X_{\eps,a}+\varphi{\normal}_{\eps,a})_{\theta}]\cdot \widehat\Xeth\\[3pt]
&\quad=\displaystyle |\Xet||\Xeth|\left(-\frac{\varphi_\theta}{|\Xeth|}+{\frac{\Xett\cdot\Ne}{|\Xet|}\frac{\varphi_\theta\varphi}{|\Xet||\Xeth|}}-\frac{\Xetth\cdot\Ne}{|\Xet|}\frac{\varphi_t \varphi}{|\Xet||\Xeth|}\right).
 \end{split}
 \end{equation} 
 \end{small}
 Then, exploiting the bound \eqref{eq:boundweighvarphi} for $\varphi$ and its derivatives, recalling \eqref{eq:XtwedgeXthetavarphi2}, Lemma \ref{lem:expsecderXN}, taking into account Remark \ref{rem:appendix1},  Remark \ref{rem:appendix2}, we infer that, up to choosing a smaller $\eps_0>0$, one has
\[
|(X_{\eps,a}+\varphi{\normal}_{\eps,a})_{t}\wedge (X_{\eps,a}+\varphi{\normal}_{\eps,a})_{\theta}|=|\Xet||\Xeth| \sqrt{1-4\mathfrak{M}(\Xe)\varphi+\eps^{2\bta} h^{(0)}_{\eps,a}}\,,
\]  
where $h_{\eps,a}^{(0)}= h_{\eps,a}^{(0)}(t,\theta)\colon\mathbb{R}\times\mathbb{R}/_{2\pi}\to\mathbb{R}$ and satisfies \eqref{eq:boundremainderexpmean}, for some constant $C>0$ independent of $a$, $\eps$, $\varphi$. Moreover, if $\frac{\pi}{\eps h_{a}} \in \mathbb{N}$, then, since $\varphi\colon \mathbb{R}/_{2n\tau_{a}}\times\mathbb{R}/_{2\pi}\to \R$, where $n=\frac{\pi}{\eps h_{a}}$, and the functions appearing in Lemma \ref{lem:expsecderXN} are doubly periodic with respect to $K_a=[-\tau_{a},\tau_{a}]\times[-\pi,\pi]$ we easily deduce that $h_{\eps,a}^{(0)}$ is doubly periodic with respect to $K_{n,a}=[-n\tau_{a},n\tau_{a}]\times[-\pi,\pi]$. 

Arguing in a similar way, from \eqref{eq:XtwedgeXthetavarphi2}, \eqref{eq:XtvarwedgeXthetavartang} one has
\begin{equation}\label{eq:normvarphiexpansion}
\begin{split}
N_\varphi&=\frac{1}{\sqrt{1-4\mathfrak{M}(\Xe)\varphi+\eps^{2\bta} h^{(0)}_{\eps,a}}}\left[ \left(-\frac{\varphi_t}{|\Xet|}+\eps^{2\bta}h_{\eps,a}^{(1)}\right)\widehat\Xet\right.
\\
&\qquad\qquad\left.
+ \left(-\frac{\varphi_\theta}{|\Xeth|}+\eps^{2\bta}h_{\eps,a}^{(2)}\right)\widehat\Xeth + \left(1-2\mathfrak{M}(\Xe)\varphi+\eps^{2\bta}h_{\eps,a}^{(3)}\right)\Ne\right]\,,
\end{split}
\end{equation}
for some functions $h_{\eps,a}^{(j)}=h_{\eps,a}^{(j)}(t,\theta)\colon\mathbb{R}\times\mathbb{R}/_{2\pi}\to\mathbb{R}$, $j=1,2,3$, satisfying \eqref{eq:boundremainderexpmean}, and as before if $\frac{\pi}{\eps h_{a}}\in \mathbb{N}$ they are doubly periodic in $K_{n,a}$.
\medskip

Setting for brevity 
\begin{equation}\label{eq:defDepsa}
D_{\eps,a}:=\sqrt{1-4\mathfrak{M}(\Xe)\varphi+\eps^{2\bta} h^{(0)}_{\eps,a}}\,,
\end{equation}
then by definition and \eqref{eq:normvarphiexpansion}, we have
\begin{equation}\label{eq:expLphi}
\begin{array}{lll}
\displaystyle {\L}_{\varphi}&=&\displaystyle \left(\Xett+\varphi_{tt} \Ne +2\varphi_t \Net + \varphi \Nett\right)\cdot N_\varphi\\[6pt] 
&=&\displaystyle \frac{1}{D_{\eps,a}}\left( - \frac{\Xett\cdot \Xet}{|\Xet|^2}\varphi_t  -  \frac{\Xett\cdot \Xeth}{|\Xeth|^2} \varphi_\theta +(\Xett\cdot \Ne) -2\mathfrak{M}(\Xe)(\Xett\cdot \Ne)\varphi   \right.\\[22pt]
&&\displaystyle \ \ \ \left. +\eps^{2\bta}  h_{\eps,a}^{(1)}  \frac{\Xett\cdot \Xet}{|\Xet|}  + \eps^{2\bta}  h_{\eps,a}^{(2)}  \frac{\Xett\cdot \Xeth}{|\Xeth|}  + \eps^{2\bta}  h_{\eps,a}^{(3)} (\Xett\cdot \Ne)\right.\\[12pt]
&&\displaystyle\left. \ \ \ + \varphi_{tt}   -2\mathfrak{M}(\Xe)\varphi_{tt}\varphi+ \eps^{2\bta}  h_{\eps,a}^{(3)}\varphi_{tt}  -2  \frac{\Net\cdot \Xet}{|\Xet|^2}\varphi_t^2- 2 \frac{\Net\cdot \Xeth}{|\Xeth|^2}  \varphi_t \varphi_\theta  \right.\\[12pt]
&&\displaystyle\ \ \ \left.  +2\eps^{2\bta}  h_{\eps,a}^{(1)}  \frac{\Net\cdot \Xet}{|\Xet|} \varphi_t +2 \eps^{2\bta}  h_{\eps,a}^{(2)} \frac{\Net\cdot \Xeth}{|\Xeth|} \varphi_t -  \frac{\Nett\cdot \Xet}{|\Xet|^2} \varphi_t\varphi \right.\\[12pt]
&&\displaystyle \ \ \ \left. -  \frac{\Nett\cdot \Xeth}{|\Xeth|^2} \varphi_\theta \varphi + ((\Ne)_{tt}\cdot \Ne)\varphi -2\mathfrak{M}(\Xe)(\Nett\cdot \Ne)\varphi^2\right.\\[12pt]
&&\displaystyle\left. \ \ \   +\eps^{2\bta}  h_{\eps,a}^{(1)}  \frac{\Nett\cdot \Xet}{|\Xet|} \varphi + \eps^{2\bta}  h_{\eps,a}^{(2)} \frac{\Nett\cdot \Xeth}{|\Xeth|} \varphi+\eps^{2\bta}  h_{\eps,a}^{(3)} (\Nett\cdot \Ne)\varphi  \right)
\end{array}
\end{equation}
Now, regrouping the terms, using the identities \eqref{identitN}--\eqref{idDerSecX2}, taking into account \eqref{eq:norm_normal_t_eps-toroidal}, \eqref{eq:dersquarednormttheta}, \eqref{eq:2lem:expNsecX}, \eqref{eq:boundweighvarphi}, Proposition \ref{prop:expmeancurv}, Remarks \ref{rem:appendix1} and \ref{rem:appendix2}, we can rewrite \eqref{eq:expLphi} as
\begin{equation}\label{eq:expLphi_2}
\begin{split}
\displaystyle {\L}_{\varphi}&=\displaystyle \frac{1}{D_{\eps,a}}\left(\Xett\cdot \Ne + \varphi_{tt} - \frac{1}{2|\Xet|^2}\frac{d|\Xet|^2}{dt}\varphi_t  + \frac{1}{2|\Xeth|^2}\frac{d|\Xet|^2}{d\theta} \varphi_\theta  \right.\\[6pt]
&\quad\displaystyle \ \ \ \ \ \ \ \ \left. -2\mathfrak{M}(\Xe)[\Xett\cdot \Ne]\varphi   - |\Net|^2 \varphi +\eps^{2\bta}  h_{\eps,a}^{(4)} \right),
\end{split}
\end{equation}
for some function $h_{\eps,a}^{(4)}=h_{\eps,a}^{(4)}(t,\theta)\colon\mathbb{R}\times\mathbb{R}/_{2\pi}\to\mathbb{R}$ satisfying \eqref{eq:boundremainderexpmean}, doubly periodic in $K_{n,a}$ if $\frac{\pi}{\eps h_{a}}\in \mathbb{N}$.

By the same argument we easily check that
\begin{equation}\label{eq:expcalNphi}
\begin{split}
\displaystyle {\N}_{\varphi}&=\displaystyle \frac{1}{D_{\eps,a}}\left(\Xethth\cdot \Ne + \varphi_{\theta\theta} + \frac{1}{2|\Xet|^2}\frac{d|\Xeth|^2}{dt}\varphi_t  - \frac{1}{2|\Xeth|^2}\frac{d|\Xeth|^2}{d\theta} \varphi_\theta  \right.\\[6pt]
&\quad\displaystyle \ \ \ \ \ \ \ \ \left. -2\mathfrak{M}(\Xe)[\Xethth\cdot \Ne]\varphi   - |\Neth|^2 \varphi +\eps^{2\bta}  h_{\eps,a}^{(5)} \right)
\end{split}
\end{equation}
and
\begin{equation}\label{eq:expcalMphi}
\begin{split}
\displaystyle {\M}_{\varphi}&=\displaystyle \frac{1}{D_{\eps,a}}\left((\Xe)_{t\theta}\cdot \Ne + \varphi_{t\theta} - \frac{1}{2|\Xet|^2}\frac{d|\Xet|^2}{d\theta}\varphi_t  - \frac{1}{2|\Xeth|^2}\frac{d|\Xeth|^2}{dt} \varphi_\theta  \right.\\[6pt]
&\quad\displaystyle \ \ \ \ \ \ \ \ \left. -2\mathfrak{M}(\Xe)[(\Xe)_{t\theta}\cdot \Ne   - \Net\cdot\Neth] \varphi +\eps^{2\bta}  h_{\eps,a}^{(6)} \right)\,,
\end{split}
\end{equation}
for some functions $h_{\eps,a}^{(5)}, h_{\eps,a}^{(6)}\colon\mathbb{R}\times\mathbb{R}/_{2\pi}\to\mathbb{R}$ satisfying \eqref{eq:boundremainderexpmean}, doubly periodic in $K_{n,a}$ if $\frac{\pi}{\eps h_{a}}\in \mathbb{N}$.
\medskip

At the end from \eqref{eq:EFG} and \eqref{eq:expLphi_2}--\eqref{eq:expcalMphi}, and taking into account \eqref{eq:boundweighvarphi}, after some elementary computations one has
\begin{equation}\label{eq:EN2FMGL}
\begin{split}
&\displaystyle \displaystyle \E_\varphi \N_\varphi -2\F_\varphi \M_\varphi+\G_\varphi {\L}_{\varphi}=\displaystyle \frac{1}{D_{\eps,a}}\left[(\Xett\cdot \Ne)|\Xeth|^2 + (\Xethth\cdot \Ne)|\Xet|^2\right.\\
&\quad\displaystyle \left. + \varphi_{tt} |\Xeth|^2  + \varphi_{\theta\theta} |\Xet|^2+ \frac{1}{2}\frac{d|\Xeth|^2}{dt} \varphi_t  - \frac{|\Xeth|^2}{2|\Xet|^2}\frac{d|\Xet|^2}{dt}\varphi_t   \right.\\[6pt]
&\displaystyle\quad\left. + \frac{1}{2}\frac{d|\Xet|^2}{d\theta} \varphi_\theta  - \frac{|\Xet|^2}{2|\Xeth|^2}\frac{d|\Xeth|^2}{d\theta}\varphi_\theta - 4 (\Xett\cdot \Ne) (\Xethth\cdot \Ne) \varphi\right.\\[6pt]
&\displaystyle\quad\left. + 4 ((\Xe)_{t\theta}\cdot\Ne)^2 \varphi  -2\mathfrak{M}(\Xe)\left((\Xett\cdot \Ne) |\Xeth|^2+(\Xethth\cdot \Ne) |\Xet|^2\right)\varphi \right.\\[6pt]
&\displaystyle\quad\left.  - \left(|\Net|^2 |\Xeth|^2 + |\Neth|^2 |\Xet|^2\right)\varphi +\eps^{2\bta}  \sigma_{\eps,a} \right],
\end{split}
\end{equation}
for some function $\sigma_{\eps,a}=\sigma_{\eps,a}(t,\theta)\colon\mathbb{R}\times\mathbb{R}/_{2\pi}\to\mathbb{R}$ satisfying 
\begin{equation}\label{eq:stimaremaindersigma}
\|x_a^{-2}\sigma_{\eps,a}\|_{C^0(\R\times[-\pi,\pi])}\leq C,
\end{equation}
where $C$ is a positive constant independent of $a$, $\eps$, $\varphi$, and again, if $\frac{\pi}{\eps h_{a}}\in \mathbb{N}$, then we easily check that $\sigma_{\eps,a}$ is doubly periodic with respect to the rectangle $K_{n,a}$.
\medskip

On the other hand, from \eqref{eq:EFG}, Lemma \ref{lem:expsecderXN}--Lemma \ref{lem:expNsecX}, exploiting \eqref{eq:boundweighvarphi} and taking into account Remark \ref{rem:appendix1}, Remark \ref{rem:appendix2}, \eqref{eq:2xasqmeancurv}, we get that
\begin{equation}\label{eq:EG-F2}
\begin{array}{lll}
&&\displaystyle \displaystyle \E_\varphi \G_\varphi-2\F_\varphi^2\\[6pt]
&=& \displaystyle |\Xet|^2|\Xeth|^2\left[ 1 -2\left(\frac{\Xett\cdot\Ne}{|\Xet|^2}+\frac{\Xethth\cdot\Ne}{|\Xeth|^2}\right) \varphi + \frac{\varphi_t^2}{|\Xet|^2} + \frac{\varphi_\theta^2}{|\Xeth|^2} \right.\\[12pt]
&&\displaystyle \left. \ +\left(\frac{|\Net|^2}{|\Xet|^2}+\frac{|\Neth|^2}{|\Xeth|^2}\right) \varphi^2 + \frac{(\Xett\cdot\Ne)(\Xethth\cdot\Ne)-((\Xe)_{t\theta}\cdot \Ne)^2}{|\Xet|^2|\Xeth|^2}\varphi^2 \right.\\[12pt]
  &&\displaystyle \left.  \  -2\frac{\Xett\cdot\Ne}{|\Xet|^2|\Xeth|^2}\varphi_\theta^2\varphi -2\frac{\Xethth\cdot\Ne}{|\Xet|^2|\Xeth|^2}\varphi_t^2\varphi  + 4  \frac{((\Xe)_{t\theta}\cdot \Ne)(\Net\cdot\Neth)}{|\Xet|^2|\Xeth|^2} \varphi^3  \right.\\[12pt]
 &&\displaystyle \left. \  -2\frac{(\Xett\cdot\Ne)|\Neth|^2+(\Xethth\cdot\Ne)|\Net|^2}{|\Xet|^2|\Xeth|^2}\varphi^3 + 4  \frac{(\Xe)_{t\theta}\cdot \Ne}{|\Xet|^2|\Xeth|^2} \varphi_t\varphi_\theta \varphi \right.\\[12pt]
  &&\displaystyle \left. \  +\frac{|\Net|^2}{|\Xet|^2|\Xeth|^2}\varphi_\theta^2 \varphi^2 +\frac{|\Neth|^2}{|\Xet|^2|\Xeth|^2}\varphi_t^2 \varphi^2 - 2  \frac{\Net\cdot\Neth}{|\Xet|^2|\Xeth|^2} \varphi_t\varphi_\theta\varphi^2 \right.\\[12pt]
  &&\displaystyle \left.  \ +  \frac{|\Net|^2|\Neth|^2-2(\Net\cdot\Neth)^2}{|\Xet|^2|\Xeth|^2} \varphi^4 \right]\\[12pt]
&=&\displaystyle |\Xet|^2|\Xeth|^2\left(1-4\mathfrak{M}(\Xe)\varphi +\eps^{2\bta}  h_{\eps,a}^{(7)} \right),
\end{array}
\end{equation}
for some function $h_{\eps,a}^{(7)}= h_{\eps,a}^{(7)}(t,\theta)\colon\mathbb{R}\times\mathbb{R}/_{2\pi}\to\mathbb{R}$ satisfying \eqref{eq:boundremainderexpmean}, doubly periodic in $K_{n,a}$ when $\frac{\pi}{\eps h_{a}}\in \mathbb{N}$.\\

From \eqref{eq:defmeancurvphi}, \eqref{eq:defDepsa}, \eqref{eq:EN2FMGL}, \eqref{eq:EG-F2}, and taking into account \eqref{eq:stimaremaindersigma} we obtain that
\begin{equation}\label{eq:expMXvarphiprovv}
\begin{array}{lll}
&&\displaystyle2x_a^2\mathfrak{M}(\Xe+\varphi\Ne)\\[12pt]
&=&\displaystyle\frac{2x_a^2}{\sqrt{1-12\mathfrak{M}(\Xe)\varphi +\eps^{2\bta}  h_{\eps,a}^{(8)}}}\left\{\frac{1}{2} \left[ \frac{(\Xett\cdot \Ne)}{|\Xet|^2} + \frac{(\Xethth\cdot \Ne)}{|\Xeth|^2} \right]\right.\\[16pt]
&&\displaystyle \left.  + \frac{\varphi_{tt}}{2|\Xet|^2}  + \frac{\varphi_{\theta\theta}}{ 2|\Xeth|^2} + \left(  \frac{1}{4|\Xet|^2|\Xeth|^2}\frac{d|\Xeth|^2}{dt}   - \frac{1}{4|\Xet|^4}\frac{d|\Xet|^2}{dt}\right)\varphi_t \right.\\[12pt]
&&\displaystyle \left.  + \left(\frac{1}{4|\Xet|^2|\Xeth|^2}\frac{d|\Xet|^2}{d\theta}   - \frac{1}{4|\Xeth|^4}\frac{d|\Xeth|^2}{d\theta}\right)\varphi_\theta \right.\\[12pt]
&&\displaystyle \left.   + \left[- 2\frac{(\Xett\cdot \Ne) (\Xethth\cdot \Ne)}{|\Xet|^2|\Xeth|^2}+2\frac{((X_{\eps,a})_{t\theta}\cdot N_{\eps,a})^2}{|(X_{\eps,a})_t|^2|(X_{\eps,a})_\theta|^2} - \frac{|(N_{\eps,a})_\theta|^2}{2|(X_{\eps,a})_\theta|^2} - \frac{|(N_{\eps,a})_t|^2}{2|(X_{\eps,a})_t|^2} \right.\right.\\[12pt]
&&\displaystyle\left.\left. - \mathfrak{M}(\Xe)\left(\frac{\Xett\cdot\Ne}{|\Xet|^2}+\frac{\Xethth\cdot\Ne}{|\Xeth|^2}\right)\right]\varphi \right\} + \eps^{2\bta}\frac{1}{D_{\eps,a}}\frac{2x_a^2\sigma_{\eps,a}}{|\Xet|^2|\Xeth|^2},
 \end{array}
\end{equation}
for some function $h_{\eps,a}^{(8)}= h_{\eps,a}^{(8)}(t,\theta)\colon\mathbb{R}\times\mathbb{R}/_{2\pi}\to\mathbb{R}$ satisfying \eqref{eq:boundremainderexpmean}, doubly periodic in $K_{n,a}$ if $\frac{\pi}{\eps h_{a}}\in \mathbb{N}$. Moreover, thanks to \eqref{eq:X_tX_thetaVW},  \eqref{eq:stimaremaindersigma} and since $D_{\eps,a}$ is uniformly bounded one has that
\begin{equation}\label{eq:stimahordtermremterm}
\left\|\frac{1}{D_{\eps,a}}\frac{2x_a^2\sigma_{\eps,a}}{|\Xet|^2|\Xeth|^2}\right\|_{C^0(\R\times[-\pi,\pi])}\leq C,
\end{equation}
for some positive constant $C$ independent of $a$, $\eps$ and $\varphi$.
\medskip

Now, exploiting \eqref{eq:boundweighvarphi}, taking into account Proposition \ref{prop:expmeancurv} and arguing as in the proof of \eqref{eq:elemexpsqrt} one has
\begin{equation}\label{eq:expsqrt12M}
\frac{1}{\sqrt{1-12\,\mathfrak{M}(\Xe)\varphi +\eps^{2\bta}  h_{\eps,a}^{(8)}}}= 1 + 6\mathfrak{M}(\Xe)\varphi +\eps^{2\bta}  h_{\eps,a}^{(9)},
\end{equation}
 for some function $h_{\eps,a}^{(9)}= h_{\eps,a}^{(9)}(t,\theta)\colon\mathbb{R}\times\mathbb{R}/_{2\pi}\to\mathbb{R}$ satisfying \eqref{eq:boundremainderexpmean}, doubly periodic in $K_{n,a}$ when $\frac{\pi}{\eps h_{a}}\in \mathbb{N}$.\\ 
 
 Finally, from \eqref{eq:expMXvarphiprovv}--\eqref{eq:expsqrt12M}, using repetitively \eqref{eq:2xasqmeancurv}, exploiting \eqref{eq:boundweighvarphi}, and recalling \eqref{M-prime-eps-a}, \eqref{eq:L-eps-a}, \eqref{eq:linearized}, then after some elementary computations one has
\[
\begin{split}
\displaystyle2x_a^2\mathfrak{M}(\Xe+\varphi\Ne)&=\displaystyle 2x_a^2 \mathfrak{M}(\Xe)+ 2x_a^2 \left[\frac{\varphi_{tt}}{2|\Xet|^2} + \frac{\varphi_{\theta\theta}}{ 2|\Xeth|^2} \right.\\[6pt]
&\quad\displaystyle \left. + \left(  \frac{1}{4|\Xet|^2|\Xeth|^2}\frac{d|\Xeth|^2}{dt}   - \frac{1}{4|\Xet|^4}\frac{d|\Xet|^2}{dt}\right)\varphi_t\right.\\[6pt]
&\quad\displaystyle \left. + \left(\frac{1}{4|\Xet|^2|\Xeth|^2}\frac{d|\Xet|^2}{d\theta}   - \frac{1}{4|\Xeth|^4}\frac{d|\Xeth|^2}{d\theta}\right)\varphi_\theta \right.\\[6pt]
&\quad\displaystyle \left.   + \left((2\mathfrak{M}(\Xe))^2- 2\frac{(\Xett\cdot \Ne) (\Xethth\cdot \Ne)}{|\Xet|^2|\Xeth|^2} \right.\right.\\[6pt]
&\quad\displaystyle \left.  \left. \ \ +2\frac{((X_{\eps,a})_{t\theta}\cdot N_{\eps,a})^2}{|(X_{\eps,a})_t|^2|(X_{\eps,a})_\theta|^2} - \frac{|(N_{\eps,a})_\theta|^2}{2|(X_{\eps,a})_\theta|^2} - \frac{|(N_{\eps,a})_t|^2}{2|(X_{\eps,a})_t|^2}\right)\varphi  \right]+ \eps^{2\bta} h_{\eps,a}^{(10)}
\end{split}
\]
\[
\begin{split}
\phantom{\displaystyle2x_a^2\mathfrak{M}(\Xe+\varphi\Ne)}&=\displaystyle 2x_a^2 \mathfrak{M}(\Xe)+2x_a^2 \left[\frac{\varphi_{tt}}{2|\Xet|^2} + \frac{\varphi_{\theta\theta}}{ 2|\Xeth|^2} \right.\\[6pt]
&\quad\displaystyle \left. + \left(  \frac{1}{4|\Xet|^2|\Xeth|^2}\frac{d|\Xeth|^2}{dt}   - \frac{1}{4|\Xet|^4}\frac{d|\Xet|^2}{dt}\right)\varphi_t\right.\\[6pt]
&\quad\displaystyle \left. + \left(\frac{1}{4|\Xet|^2|\Xeth|^2}\frac{d|\Xet|^2}{d\theta}   - \frac{1}{4|\Xeth|^4}\frac{d|\Xeth|^2}{d\theta}\right)\varphi_\theta \right.\\[6pt]
&\quad\displaystyle \left.   + \left(\frac{(\Xett\cdot \Ne)^2}{|\Xet|^4}+\frac{(\Xethth\cdot \Ne)^2}{|\Xeth|^4} \right.\right.\\[6pt]
&\quad\displaystyle \left.  \left. \ \ +2\frac{((X_{\eps,a})_{t\theta}\cdot N_{\eps,a})^2}{|(X_{\eps,a})_t|^2|(X_{\eps,a})_\theta|^2} - \frac{|(N_{\eps,a})_\theta|^2}{2|(X_{\eps,a})_\theta|^2} - \frac{|(N_{\eps,a})_t|^2}{2|(X_{\eps,a})_t|^2}\right)\varphi  \right]+ \eps^{2\bta} h_{\eps,a}^{(10)}\\[6pt]
&=\displaystyle 2x_a^2 \mathfrak{M}(\Xe)+\mathcal{L}_{\eps,a}\varphi + \eps^{2\bta} h_{\eps,a}^{(10)}
 \end{split}
\]
for some function $h_{\eps,a}^{(10)}= h_{\eps,a}^{(10)}(t,\theta)\colon\mathbb{R}\times\mathbb{R}/_{2\pi}\to\mathbb{R}$ satisfying \eqref{eq:boundremainderexpmean}, doubly periodic in $K_{n,a}$ if $\frac{\pi}{\eps h_{a}}\in \mathbb{N}$. This ends the proof of {\it (ii)}.
\qed
As in immediate consequence of Proposition \ref{prop:expmeancurvnormgraph}-(ii) we obtain:
\begin{Corollary}\label{cor:stimaB0}
Fixing $a_{0}\in(0,\frac{1}{2})$, for every $R>0$ and $\bta>0$ there exists $\eps_0>0$ such that for any $a\in[-a_0,a_0]\setminus\{0\}$, $\eps\in(0,\eps_0)$,  $\varphi \in \mathscr{B}_{\eps,a}(\bta,R)$ one has
\[
\|\mathfrak{B}_{\eps,a}^0(\varphi)\|_{C^0(\R\times [-\pi,\pi])} \leq C \eps^{2\bta}\,,
\]
for some positive constant $C$ independent of $a$, $\eps$,  $\varphi$. If in addition $\frac{\pi}{\eps h_a} \in \mathbb{N}$ then $\mathfrak{B}_{\eps,a}^0(\varphi)$ is doubly periodic with respect to the rectangle $K_{n,a}=[-n\tau_a,n\tau_a]\times[-\pi,\pi]$, where $n=\frac{\pi}{\eps h_a}$.
\end{Corollary}

In the last part of this section we study the right-hand side of \eqref{eq:secondtermdecom} and we discuss the following result.

\begin{Proposition}\label{prop:expHXvarphiepsa}
Let $R>0$ and let $H\colon\R^3\to \R$ be a function of class $C^{1}$ satisfying \eqref{eq:H1} and \eqref{eq:H3} for some $A\colon\mathbb{S}^{2}\to\mathbb{R}$, $\bta>0$ and $\nu>0$. Then, fixing $a_0\in(0,\frac{1}{2})$,  there exist $\eps_0>0$ and $C>0$ such that for any  $a\in[-a_0,a_0]\setminus\{0\}$, $\eps\in(0,\eps_0)$, and for any $\varphi\in\mathscr{B}_{\eps,a}(\bta,R)$ one has 
\begin{equation}\label{eq:prop:expHXvarphiepsa}
2x_a^2 H(X_{\varepsilon,a}+\varphi N_{\varepsilon,a})=2x_a^2+  2x_{a}^{2}\eps^\bta A(\widehat{X}_{\eps,a}) +  \eps^{\bta+\widetilde\nu} \xi_{\eps,a}
\end{equation}
for some function $\xi_{\eps,a}=\xi_{\eps,a}(t,\theta)\colon\mathbb{R}\times\mathbb{R}/_{2\pi}\to\mathbb{R}$ satisfying
\begin{equation}\label{eqtesi2:prop:expHXepsa}
\|\xi_{\eps,a}\|_{C^0(\R \times [-\pi,\pi])} \leq C,
\end{equation}
and $\widetilde\nu=\min\{1,\nu\}$. If in addition $\frac{\pi}{\eps h_a}\in \mathbb{N}$ then $\xi_{\eps,a}$ is doubly periodic with respect to the rectangle $K_{n,a}=[-n\tau_{a},n\tau_{a}]\times[-\pi,\pi]$, where $n=\frac{\pi}{\eps h_a}$.
\end{Proposition}

We firstly state a couple of preliminary lemmata. 

\begin{Lemma}\label{lem:expHXepsa}
Let $H\colon\R^3\to \R$ be a function satisfying \eqref{eq:H1} for some continuous function $A\colon\mathbb{S}^{2}\to\mathbb{R}$, $\bta>0$ and $\nu>0$. Then, fixing $a_0\in(0,\frac{1}{2})$, there exist $\eps_0>0$ and $C>0$ such that for any  $a\in[-a_0,a_0]\setminus\{0\}$ and $\eps\in(0,\eps_0)$ one has
\begin{equation}\label{eqtesi:lem:expHXepsa}
2x_a^2H(X_{\eps,a})=2x_a^2+  2x_a^2\eps^\bta A(\widehat{X}_{\eps,a}) +  \eps^{\bta+\widetilde\nu} \rho_{\eps,a}\quad\text{where}\quad\widehat{X}_{\eps,a}=\frac{{X}_{\eps,a}}{|{X}_{\eps,a}|}\,,
\end{equation}
$\widetilde\nu=\min\{1,\nu\}$ and $\rho_{\eps,a}=\rho_{\eps,a}(t,\theta)\colon\mathbb{R}\times\mathbb{R}/_{2\pi}\to\mathbb{R}$ satisfies
\begin{equation}\label{eqtesi2:lem:expHXepsa}
\|x_a^{-2}\rho_{\eps,a}\|_{C^0(\R\times[-\pi,\pi])} \leq C\,.
\end{equation}
If in addition $\frac{\pi}{\eps h_a} \in \mathbb{N}$ then $\rho_{\eps,a}$ is doubly-periodic with respect to $K_{n,a}=[-n\tau_{a},n\tau_{a}]\times[-\pi,\pi]$, where $n=\frac{\pi}{\eps h_a}$.
\end{Lemma}
\Proof
By \eqref{eq:H1}, for any $X\in\R^{3}$ with $|X|$ large enough we can write
\[
H(X)=1+\frac{A(\widehat{X})}{|X|^{\bta}}+\frac{H_{1}(X)}{|X|^{\bta+\nu}}
\]
where $\widehat{X}=X/|X|$ and $H_{1}$ is bounded. 
By direct computation (see \eqref{eq:eps-toroidal-unduloid} or \eqref{eq:eps-toroidal-unduloid_rotmatrix}) we have 
\begin{equation}\label{eq:moduloXepsa}
|X_{\eps,a}|=\eps^{-1} \sqrt{1+2\eps x_a\sin\theta + \eps^2 x_a^2}.
\end{equation}
Then, fixing $a_{0}\in(0,\frac12)$ and arguing as in the proof of \eqref{eq:elemexpsqrt}, taking into account Remark \ref{rem:appendix1}, we can find $\eps_0>0$ such that for all $a\in[-a_0,a_0]\setminus\{0\}$, $\eps\in(0,\eps_0)$
\begin{equation}\label{eq2:lem:expHXepsa}
\frac{1}{|X_{\eps,a}|^{\bta}}=
\eps^{\bta}\left(1+\eps x_{a}^{-2}\rho^{(1)}_{\eps,a}\right)\,,\quad
\frac{1}{|X_{\eps,a}|^{\bta+\nu}}=\eps^{\bta+\nu}x_{a}^{-2}\rho^{(2)}_{\eps,a}
\end{equation}
for some doubly-periodic functions $\rho_{\eps,a}^{(i)}=\rho_{\eps,a}^{(i)}(t,\theta)$ with respect to the rectangle $K_a=[-\tau_a,\tau_a]\times[-\pi,\pi]$, $(i=1,2)$, satisfying 
\begin{equation}\label{eq2:lem:expHXepsabis}
\|x_a^{-2}\rho^{(i)}_{\eps,a}\|_{C^0(K_{a})} \leq C
\end{equation}
for some positive constant $C$ independent of $a$, $\eps$. Then, by \eqref{eq:moduloXepsa}--\eqref{eq2:lem:expHXepsa}, after multiplication by $2x_{a}^{2}$ we obtain
\[
2x_{a}^{2}H(X_{\eps,a})=2x_{a}^{2}+2x_{a}^{2}\eps^{\bta}A(\widehat{X}_{\eps,a})+
2\eps^{\bta+1}\rho^{(1)}_{\eps,a}A(\widehat{X}_{\eps,a})+2\eps^{\bta+\nu}\rho^{(2)}_{\eps,a}
\]
which readily implies \eqref{eqtesi:lem:expHXepsa}, and taking into account \eqref{eq2:lem:expHXepsabis} we obtain \eqref{eqtesi2:lem:expHXepsa}. Moreover, when $\frac{\pi}{\eps h_a} \in \mathbb{N}$, then we easily check that $\rho_{\eps,a}$ is doubly-periodic with respect to $K_{n,a}=[-n\tau_{a},n\tau_{a}]\times[-\pi,\pi]$, with $n=\frac{\pi}{\eps h_a}$.
\qed

\begin{Lemma}\label{lem:expnablaHXepsa}
Let $H\colon\R^3\to \R$ be a function of class $C^{1}$ satisfying \eqref{eq:H3} for some $\bta>0$ and let $R>0$. Then, fixing $a_0\in(0,\frac{1}{2})$, there exist $\eps_0>0$ and $C>0$ such that for any  $a\in[-a_0,a_0]\setminus\{0\}$, $\eps\in(0,\eps_0)$, for any $\varphi\in\mathscr{B}_{\eps,a}(\bta,R)$ one has
\begin{itemize}
\item[(i)] $\|2x_a^2 (\nabla H (X_{\varepsilon,a})\cdot N_{\eps,a}) \varphi\|_{C^0(\R\times[-\pi, \pi])}  \leq C \eps^{2\bta+1}$;
\item[(ii)] $\|\mathfrak{B}_{\eps,a}^1(\varphi)\|_{C^0(\R\times[-\pi, \pi])}  \leq C \eps^{2\bta+1}$.
\end{itemize}
If in addition $\frac{\pi}{\eps h_a} \in \mathbb{N}$ then $2x_a^2 (\nabla H (X_{\varepsilon,a})\cdot N_{\eps,a}) \varphi$, $\mathfrak{B}_{\eps,a}^1(\varphi)$ are doubly periodic with respect to $K_{n,a}=[-n\tau_{a},n\tau_{a}]\times[-\pi,\pi]$, where $n=\frac{\pi}{\eps h_a}$.
\end{Lemma}
\Proof
{\it (i)} By \eqref{eq:H3} there exists $C>0$ such that 
\[
|\nabla H (X)|\le\frac{C}{|X|^{\bta+1}}
\]
for every $X\in\mathbb{R}^{3}$ with $|X|$ large enough. Then, fixing $a_{0}\in(0,\frac12)$ and arguing as in the proof of Lemma \ref{lem:expHXepsa}, we can find $\eps_0>0$ and $C>0$ such that for all $a\in[-a_0,a_0]\setminus\{0\}$ and $\eps\in(0,\eps_0)$
\[
|\nabla H (X_{\varepsilon,a})\cdot N_{\eps,a}|\le C_{1}\eps^{\bta+1}
\]
for some constant $C_{1}$ independent of $a$ and $\eps$. Then, thanks to \eqref{eq:boundweighvarphi}, {\it (i)} follows. 
\medskip

{\it (ii)} By definition \eqref{eq:secondtermdecom}, we can write $\mathfrak{B}_{\eps,a}^1(\varphi)$  as
\begin{equation}\label{eq:writeBepsa1}
\mathfrak{B}_{\eps,a}^1(\varphi)=\mathfrak{J}_{\eps,a}(\varphi)\varphi
\end{equation}
where 
\[
\mathfrak{J}_{\eps,a}(\varphi):=\int_0^1\left(\nabla H(\Xe+r\varphi\Ne)-\nabla H(\Xe)\right)\cdot \Ne\, dr\,.
\]
Then, fixing $a_{0}\in(0,\frac12)$,  using again \eqref{eq:H3} and arguing as in the proof of \textit{(i)}, we find $\eps_{0}>0$ and $C>0$ such that for all $a\in[-a_0,a_0]\setminus\{0\}$, $\eps\in(0,\eps_0)$ and $\varphi\in \mathscr{B}_{\eps,a}(\bta,R)$ 
\begin{equation}\label{eq:stimaunif3}
|\mathfrak{J}_{\eps,a}(\varphi)|\leq \sup_{r\in[0,1]}|\nabla  H(\Xe+r\varphi\Ne)|+|\nabla H(\Xe)| \leq\displaystyle C \eps^{\bta+1} \ \ \forall (t,\theta)\in\R^2.
\end{equation}
Finally, joining \eqref{eq:writeBepsa1}, \eqref{eq:stimaunif3} and exploiting \eqref{eq:boundweighvarphi} we readily deduce {\it (ii)}. Moreover, if $\frac{\pi}{\eps h_a}\in \mathbb{N}$, by construction and \eqref{eq:writeBepsa1}, we immediately check that $2x_a^2 (\nabla H (X_{\varepsilon,a})\cdot N_{\eps,a}) \varphi$ and $\mathfrak{B}_{\eps,a}^1(\varphi)$ are doubly periodic in $K_{n,a}=[-n\tau_{a},n\tau_{a}]\times[-\pi,\pi]$, with $n=\frac{\pi}{\eps h_a}$.
\qed

\noindent
{\it Proof of Proposition \ref{prop:expHXvarphiepsa}.}
Recalling \eqref{eq:secondtermdecom} and using Lemmata \ref{lem:expHXepsa} and \ref{lem:expnablaHXepsa}, fixing $a_0\in(0,\frac{1}{2})$, there exist $\eps_0>0$ such that for every $a\in[-a_0,a_0]\setminus\{0\}$, $\eps \in (0,\eps_0)$ and $\varphi\in\mathscr{B}_{\eps,a}(\bta,R)$ one has
\begin{equation}\label{eq1:propexpHvarphi}
2x_a^2 H(X_{\varepsilon,a}+\varphi N_{\varepsilon,a})=\displaystyle 2x_a^2+ 2x_a^2  \eps^\bta A(\widehat{X}_{\eps,a}) +  \eps^{\bta+\widetilde\nu} \rho_{\eps,a} + 2x_a^2 (\nabla H (X_{\varepsilon,a})\cdot N_{\eps,a}) \varphi + \mathfrak{B}_{\eps,a}^1(\varphi)
\end{equation}
where $\rho_{\eps,a}=\rho_{\eps,a}(t,\theta)\colon\mathbb{R}\times\mathbb{R}/_{2\pi}\to\mathbb{R}$ satisfies \eqref{eqtesi2:lem:expHXepsa}, for some positive constant $C$ independent of $a$, $\eps$, $\varphi$. Then from \eqref{eq1:propexpHvarphi} we infer that \eqref{eq:prop:expHXvarphiepsa} holds with
\[
\xi_{\eps,a}:=\rho_{\eps,a} + \eps^{-\bta-\widetilde\nu}\left[ 2x_a^2 (\nabla H (X_{\varepsilon,a})\cdot N_{\eps,a})\varphi +\mathfrak{B}_{\eps,a}^1(\varphi)\right]\,.
\]
Then, taking into account \eqref{eqtesi2:lem:expHXepsa}, Remark \ref{rem:appendix2} and the estimates {\it (i)} and  {\it (ii)} of Lemma \ref{lem:expnablaHXepsa}, we easily check that $\xi_{\eps,a}$ satisfies \eqref{eqtesi2:prop:expHXepsa} with a constant $C>0$ independent of $a$, $\eps$, $\varphi$.

 Finally, if $\frac{\pi}{\eps h_a}\in \mathbb{N}$ then $\xi_{\eps,a}$ is doubly periodic  with respect to $K_{n,a}=[-n\tau_{a},n\tau_{a}]\times[-\pi,\pi]$, where $n=\frac{\pi}{\eps h_a}$, because $\rho_{\eps,a}$, $2x_a^2 (\nabla H (X_{\varepsilon,a})\cdot N_{\eps,a}) \varphi$ and $\mathfrak{B}_{\eps,a}^1(\varphi)$ are so (see Lemmata \ref{lem:expHXepsa} and \ref{lem:expnablaHXepsa}).
\qed

\section{Proof of Theorems \ref{mainteo} and \ref{mainteo2}}
In this section we prove Theorems \ref{mainteo} and \ref{mainteo2}. To this aim we firstly need to state some facts occurring in the singular limit, that is, as $a\to 0$. In this limit, taking account of Lemma \ref{L:tau-a-expansion} $(iii)$, the normalized Jacobi operator $\mathcal{L}_{a}$ introduced in \eqref{normalizedJacop-a} becomes
\begin{equation}\label{eq:defL0}
\mathcal{L}_0\varphi:=\Delta\varphi+2(\sech t)^2\varphi\,.
\end{equation}
Moreover, taking the limit $a\to 0$ in the mappings belonging to the set $\mathscr{W}_{a}$ introduced in \eqref{Wa}, we obtain the set
\begin{equation}
\label{W0}
\mathscr{W}_{0}:=\{-1+t\tanh t\,,~~-\tanh t\,,~~\sech t\cos\theta\,,~~\sech t\sin\theta\}\,,
\end{equation}
having omitted the limits of $w_{a,1}^{-}\cos\theta$ and $w_{a,1}^{-}\sin\theta$, which are the null function. 
Moreover, in view of Lemma \ref{L:tau-a-expansion} again, in particular of the fact that $x_{a}(t)\to\sech t$ and $\tau_{a}\to\infty$ as $a\to 0$ (see \eqref{eq:tau-a}), the bound with respect to the norm \eqref{eq:weighted-C2-norm} becomes a condition of the form
\begin{equation}\label{eq:Lemmamainteo}
\sup_{(t,\theta)\in\mathbb{R}\times[-\pi,\pi]}(\cosh t)\left[|\varphi(t,\theta)|+|\nabla\varphi (t,\theta)|+|D^2\varphi (t,\theta)|\right]<\infty\,.
\end{equation}
We have that:
\begin{Lemma}\label{lem:techmainteo}
For every $\varphi \in C^2(\mathbb{R}\times\mathbb{R}/_{2\pi})$ satisfying \eqref{eq:Lemmamainteo} one has that
\begin{equation}\label{eq:ClaimCase1}
 \int_{\R\times[-\pi,\pi]}w\,\mathcal{L}_0\varphi\,dt\,d\theta=0\,.
\end{equation}
where $\mathcal{L}_0$ is the operator defined in \eqref{eq:defL0} and $w$ is any mapping in the class $\mathscr{W}_{0}$ defined in \eqref{W0}.
\end{Lemma}

\Proof
We firstly observe that by definition of $\mathcal{L}_0$ and in view of the assumption \eqref{eq:Lemmamainteo}, the integral in \eqref{eq:ClaimCase1} is well-defined and finite. Then 
\[
\begin{split}
\int_{\R\times[-\pi,\pi]}w\,\mathcal{L}_0\varphi\,dt\,d\theta
&=\lim_{M\to\infty}\int_{[-M,M]\times[-\pi,\pi]}w\left[\varphi_{tt}+\varphi_{\theta\theta}+2(\sech t)^{2}\varphi\right]\,dt\,d\theta\,.
\end{split}
\]
After a double integration by parts, we have
\[
\begin{split}
&\int_{[-M,M]\times[-\pi,\pi]}w\varphi_{\theta\theta}\,dt\,d\theta=\int_{[-M,M]\times[-\pi,\pi]}w_{\theta\theta}\varphi\,dt\,d\theta\\
&\int_{[-M,M]\times[-\pi,\pi]}w\varphi_{tt}\,dt\,d\theta=\int_{[-M,M]\times[-\pi,\pi]}w_{tt}\varphi\,dt\,d\theta+I(M)-I(-M)
\end{split}
\]
where
\[
I(\pm M)=\int_{-\pi}^{\pi}\left[w(\pm M,\theta)\varphi_{t}(\pm M,\theta)-w_{t}(\pm M,\theta)\varphi(\pm M,\theta)\right]\,d\theta\,.
\]
Then
\[
\int_{[-M,M]\times[-\pi,\pi]}w\left[\varphi_{tt}+\varphi_{\theta\theta}+2(\sech t)^{2}\varphi\right]\,dt\,d\theta=\int_{[-M,M]\times[-\pi,\pi]}\varphi\,\mathcal{L}_{0}w\,dt\,d\theta+I(M)-I(-M)
\]
and \eqref{eq:ClaimCase1} follows, since $\mathcal{L}_{0}w=0$ and, by the shape of $w$ and by \eqref{eq:Lemmamainteo}, $I(\pm M)\to 0$ as $M\to\infty$.
\qed


\noindent
\textit{Proof. of Theorem \ref{mainteo}.}
Let $H\colon\R^3\to \R$ be as in the statement of Theorem \ref{mainteo}. Assume by contradiction that the thesis is false. This means that there exist sequences $(\ak)_k\subset \R\setminus\{0\}$ with $\ak\to 0$, $(\nk)_k\subset \mathbb{N}$ with $\nk\to \infty$, and $(\phik)_k \subset C^{2,\alpha}(\R^2,\R)$ where for any $k\in \mathbb{N}$ the function $\phik$ is doubly periodic with respect to the rectangle $K_{n_{k},a_{k}}=[-n_{k}\tau_{a_{k}},n_{k}\tau_{a_{k}}]\times[-\pi,\pi]$, satisfies
\[
\|\phik\|_{C^{2}(K_{\nk,\ak};x_{\ak}^{-1})}\le\frac{R}{n_{k}^{\bta}}
\]
for some constant $R>0$ 
independent of $k$, and
\begin{equation}\label{eq:equationmeancurvsequence}
 \mathfrak{M}(X_{\ek,\ak}+\phik N_{\ek,\ak})=H(X_{\ek,\ak}+\phik N_{\ek,\ak})
\end{equation}
being
\[
\ek=\frac{\pi}{\nk h_{\ak}}\,.
\]
Recalling the definition \eqref{eq:U} and taking into account that, as $a_{k}\to 0$, then $h_{a_{k}}\to 1$ (see \eqref{eq:acca-a}), every $\phik$ belongs to $\mathscr{B}_{\ek,\ak}(\bta,\overline{R})$, for some $\overline{R}>0$ independent of $k$. Since $\ek\to 0$ and $\ak\to 0$, we can find $\bar k\in \mathbb{N}$ such that for all $k\geq \bar k$ the expansions provided by Propositions \ref{prop:expmeancurv}, \ref{prop:explinearized}, \ref{prop:expmeancurvnormgraph} and \ref{prop:expHXvarphiepsa} hold true.
Hence, setting for brevity $\sigma_{k}=:\sigma_{\ek,\ak}$, $h_{k}:=h_{\ek,\ak}$, $\xi_{k}:=\xi_{\ek,\ak}$ and
\begin{equation}\label{eq:defgan}
g_{k}:=2x_{\ak}^3+2z_{\ak}^\prime x_{\ak}-4(z_{\ak}^\prime)^2 x_{\ak}-\frac{(z_{\ak}^\prime)^3}{x_{\ak}}-2\gamma_{\ak}\frac{(z_{\ak}^\prime)^2}{x_{\ak}}\,,
\end{equation}
for all $k\geq\bar k$ we have
\begin{equation}\label{eq:decompequationmeancurvsequence1}
 2x_{\ak}^2\mathfrak{M}(X_{\ek,\ak}+\phik N_{\ek,\ak})=2x_{\ak}^2+\ek g_k\sin\theta + \ek^2 \sigma_k + \mathcal{L}_{\ak}\phik + \ek \mathcal{L}_{\ak}^{(1)}\phik+ \ek^2 \mathcal{L}_{\ek,\ak}^{(2)}\phik + \ek^{2\bta}h_k
\end{equation}
and
\begin{equation}\label{eq:decompequationmeancurvsequence2}
2x_{\ak}^2H(X_{\ek,\ak}+\phik N_{\ek,\ak})=2x_{\ak}^2+ 2x_{\ak}^2\ek^\bta A(\widehat{X}_{\eps_{k},a_{k}}) +  \ek^{\bta+\widetilde\nu} \xi_k\quad\text{where~~}\widehat{X}_{\eps_{k},a_{k}}=\frac{X_{\eps_{k},a_{k}}}{|X_{\eps_{k},a_{k}}|}\,,
\end{equation}
$\mathcal{L}_{\ak}$ is the normalized Jacobi operator introduced in \eqref{normalizedJacop-a}, $ \mathcal{L}_{\ak}^{(1)}$ and $ \mathcal{L}_{\ek,\ak}^{(2)}$ are the operators defined in \eqref{eq:expoperatorL1}, \eqref{eq:expoperatorL2}, and the mappings $\sigma_k$, $h_k$, $\xi_k$ are doubly periodic functions with respect to the rectangle $K_{n_{k},a_{k}}$ and  satisfy the uniform estimate 
\begin{equation}\label{eq:unifboundremterms}
\|\sigma_k\|_{C^0(\R^{2})} + \|h_k\|_{C^0(\R^{2})} + \|\xi_k\|_{C^0(\R^{2})}\leq C
\end{equation}
for some constant $C>0$ independent of $k$ (as it follows from \eqref{eq:estremtermexpmeanc} and Remark \ref{rem:appendix2}, \eqref{eq:boundremainderexpmean}, \eqref{eqtesi2:prop:expHXepsa}).
\medskip

According to the values of $\bta$ we distinguish two cases and we prove they all lead a contradiction.
\medskip

\textbf{Case 1:} Assume that $\bta\in(0,1)$.  Setting for $k\geq \bar k$
\begin{equation}
\label{phi-tilde-k}
\widetilde{\varphi}_{k}(t,\theta):=\frac{1}{\ek^\bta}\phik(t,\theta) \ \ \ (t,\theta)\in \R^2\,,
\end{equation}
we obtain a new sequence $(\widetilde{\varphi}_{k})_k\subset C^{2,\alpha}(\R^2)$ of doubly periodic functions (with respect to the rectangle $K_{n_{k},a_{k}}$), such that for all $k\geq \bar k$ one has 
\begin{equation}\label{eq:proofmainteo1}
|\widetilde{\varphi}_{k}(t,\theta)|+|\nabla\widetilde{\varphi}_{k}(t,\theta)|+|D^2\widetilde{\varphi}_{k}(t,\theta)|\le \overline{R} x_{\ak}(t),
\quad\forall(t,\theta)\in\R^2,
\end{equation}
and, by the hypothesis \eqref{eq:phi2},
\[
[D^{2}\widetilde{\varphi}_{k}]_{\alpha,K_{\nk,\ak}}\le \overline{R}\,.
\]
In particular, in view of Remark \ref{rem:appendix1}, and taking into account that $\ak\to 0$, there exists a positive constant $C_1$ independent of $k$ such that
\begin{equation}\label{eq:proofmainteo2}
\|\widetilde{\varphi}_{k}\|_{C^{2,\alpha}(\R^2)} \leq C_1 \ \ \ \forall k\geq \bar k\,.
\end{equation}
Now, observe that in terms of $\widetilde{\varphi}_{k}$ we can rewrite  \eqref{eq:decompequationmeancurvsequence1} as 
\begin{equation}\label{eq:decompequationmeancurvsequence1bis}
\begin{split}
\displaystyle 2x_{\ak}^2&\mathfrak{M}(X_{\ek,\ak}+\phik N_{\ek,\ak})\\
&=\displaystyle 2x_{\ak}^2+\ek g_k\sin\theta +\ek^{2}\sigma_k + \ek^\bta \mathcal{L}_{\ak}\widetilde{\varphi}_{k} + \ek^{1+\bta} \mathcal{L}_{\ak}^{(1)}\widetilde{\varphi}_{k}+ \ek^{2+\bta} \mathcal{L}_{\ek,\ak}^{(2)}\widetilde{\varphi}_{k} +\ek^{2\bta}h_k\\
&=\displaystyle 2x_{\ak}^2+\ek g_k\sin\theta + \ek^\bta \mathcal{L}_{\ak}\widetilde{\varphi}_{k} + \ek^{2\bta} \zeta_k,
\end{split}
\end{equation}
where we have set  
\[
\zeta_k:=\ek^{2-2\bta}\sigma_k+\ek^{1-\bta} \mathcal{L}_{\ak}^{(1)}\widetilde{\varphi}_{k}+ \ek^{2-\bta} \mathcal{L}_{\ek,\ak}^{(2)}\widetilde{\varphi}_{k} +h_k\,.
\] 
In particular $\zeta_k$ is doubly periodic and thanks to \eqref{eq:expoperatorL1}--\eqref{eq:unifboundexpLin}, \eqref{eq:unifboundremterms}, \eqref{eq:proofmainteo2},  taking into account  that $\bta \in (0,1)$, $\ak\to 0$, $\ek\to 0$ and Remark \ref{rem:appendix1}, we see that
\[
\|\zeta_k\|_{C^0(\R^2)}\leq C_2\quad \forall k\geq \bar k
\]
for some positive constant $C_2$ independent of $k$.
\medskip

At the end, from \eqref{eq:equationmeancurvsequence}, \eqref{eq:decompequationmeancurvsequence2},  \eqref{eq:decompequationmeancurvsequence1bis}, we obtain the following crucial identity
\begin{equation}\label{eq:fundamentalrelation}
\ek g_{\ak}\sin\theta + \ek^\bta \mathcal{L}_{\ak}\widetilde{\varphi}_{k} + \ek^{2\bta} \zeta_k=2x_{\ak}^2 \ek^\bta A(\widehat{X}_{\ek,\ak}) +  \ek^{\bta+\widetilde\nu} \xi_k\,.
\end{equation}
Let us observe that, since
\[
|\widehat{X}_{\eps,a}-\mathbf{e}_{2}|\le
\frac{\eps|X_{\eps,a}-\eps^{-1}\mathbf{e}_{2}|+\left(1-\eps|X_{\eps,a}|\right)}{\eps|X_{\eps,a}|}\,,
\]
by \eqref{limit-e2} and \eqref{eq:moduloXepsa}, $\widehat{X}_{\eps,a}\to\mathbf{e}_{2}$ as $\eps\to 0$ and then, by the continuity of $A$,
\begin{equation}
\label{Ae2}
A(\widehat{X}_{\ek,\ak})\to A(\mathbf{e}_{2})\,.
\end{equation}
Moreover, from \eqref{eq:proofmainteo2} and by a standard compactness result (see \cite[Lemma 6.36]{GT}), up to a subsequence (still indexed by $k$ for brevity) we get that $\widetilde{\varphi}_{k} \to \widetilde{\varphi}$ in $C^2_{loc}(\R^2)$ for some function $\tilde \varphi\in C^2(\R\times\R/_{2\pi})$. In addition, from \eqref{eq:proofmainteo1} and since $x_{\ak}(t)\to \sech t$ (see Lemma \ref{L:tau-a-expansion}-(iii)), we infer that
\begin{equation}\label{eq:proofmainteo3}
|\widetilde{\varphi}(t,\theta)|+|\nabla\widetilde{\varphi} (t,\theta)|+|D^2\widetilde{\varphi} (t,\theta)|\le \overline{R}\, \sech t,
\quad\forall(t,\theta)\in\R^2.
\end{equation}
Hence, fixing $(t,\theta)\in\R^2$, dividing each side of \eqref{eq:fundamentalrelation} by $\ek^\bta$ and passing to the limit as $k\to\infty$, taking into account of \eqref{Ae2}, that $\bta\in(0,1)$, $x_{\ak}(t)\to \sech(t)$, ${\gamma_{\ak}}/{x_{\ak}(t)}\to 0$, and that $\|g_k\|_{C^0(\R)}\leq C_3$, for some positive constant $C_3$ independent of $k$ (as it follows by Remark \ref{rem:appendix1}), we get that
\begin{equation}
\label{eq:fundamentalrelation2}
\Delta\tilde\varphi+2(\sech t)^{2}\tilde\varphi= 2 A(\mathbf{e}_{2})(\sech t)^{2} \,.
\end{equation}

Then, multiplying each side of \eqref{eq:fundamentalrelation2} by the function $1-t\tanh t$, and integrating over $\R\times [-\pi,\pi]$ (which makes sense in view of \eqref{eq:proofmainteo3}), we have
 \begin{equation}\label{eqint:mainteoproof}
\int_{\R\times[-\pi,\pi]}(1-t\tanh t)(\Delta\tilde\varphi+2(\sech t)^{2}\tilde\varphi)\, dt\,d\theta= 2\pi A(\mathbf{e}_{2}) \int_{-\infty}^{\infty}(\sech t)^{2} (1-t\tanh t) \, dt\,.
\end{equation}
Now, by Lemma \ref{lem:techmainteo} the left-hand side of \eqref{eqint:mainteoproof} is zero. Therefore, observing that
\[
\int_{-\infty}^{\infty}(\sech t)^{2} (1-t\tanh t) \, dt= \frac{1}{2} \left[t(\sech t)^{2}+ \tanh t\right]\big|_{-\infty}^{\infty}=1\,,
\]
one obtains $A(\mathbf{e}_{2})=0$, contrary to the hypothesis \eqref{eq:H2}.
\medskip

\textbf{Case 2:} Assume that $\bta=1$. As in case 1, we introduce the functions $\widetilde{\varphi}_{k}$ defined by \eqref{phi-tilde-k} with $\beta=1$. Each $\widetilde{\varphi}_{k}$ turns out to solve
\begin{equation}\label{eq:fundamentalrelation22}
\ek g_{\ak}\sin\theta + \ek \mathcal{L}_{\ak}\widetilde{\varphi}_{k} + \ek^{2} \zeta_k= 2x_{\ak}^2\ek A(\widehat{X}_{\ek,\ak}) +  \ek^{1+\widetilde\nu} \xi_k
\end{equation}
where $\zeta_k$, $\xi_k$ are some doubly periodic functions satisfying
\[
\| \zeta_k\|_{C^0(\R^{2})}+\|\xi_k\|_{C^0(\R^{2})}\leq C
\]
for some constant $C$ independent of $k$. Moreover \eqref{eq:proofmainteo1} holds true again, and then, up to a subsequence, $\widetilde{\varphi}_{k}\to \widetilde\varphi$ in $C_{loc}^2(\R^2)$ for some $\widetilde\varphi\in C^{2}(\R\times\R/_{2\pi})$ satisfying \eqref{eq:proofmainteo3}. In order to pass to the limit in \eqref{eq:fundamentalrelation22}, we observe that, in view of Lemma \ref{L:tau-a-expansion}-(iii), \eqref{eq:xz} and since $\gamma_{\ak}\to 0$, one has $x_{\ak}\to \sech t$, $z_{\ak}^\prime\to \sech^2t$, and thus we deduce that (see \eqref{eq:defgan})
\begin{equation}\label{eq:limitgn}
 g_{\ak}(t) \to 4(\sech t)^3-5(\sech t)^5=:g_{0}(t)\,.
 \end{equation}
Then, dividing \eqref{eq:fundamentalrelation22} by $\ek$, and passing to the limit as $k\to\infty$, we obtain that $\tilde\varphi\in C^{2}(\R\times\R/_{2\pi})$ satisfies 
\begin{equation}\label{eq:fundamentalrelationCase2}
g_{0}(t)\, \sin\theta +\Delta\widetilde\varphi+2(\sech t)^{2}\widetilde\varphi= 2 A(\mathbf{e}_{2})\,(\sech t)^{2}\,.
\end{equation}
Finally, multiplying \eqref{eq:fundamentalrelationCase2} by $w(t)=1-t\tanh t$, integrating over $\R\times[-\pi,\pi]$ and observing that
\[
\int_{\R\times[-\pi,\pi]} (1-t\tanh t)\,g_{0}(t)\, \sin\theta \, dt\,d\theta=0\,,
\]
then, thanks to  Lemma \ref{lem:techmainteo}, we obtain a contradiction as in \textbf{Case 1}.

It remains the case $\beta>1$ but in this situation Theorem \ref{mainteo} is a corollary of Theorem \ref{mainteo2}. 
\qed

\begin{Remark} $(i)$
We observe that in the proof of Theorem \ref{mainteo}, when $\bta\in(0,1]$, if we test the limit equation \eqref{eq:fundamentalrelation2} (in case $0<\beta<1$) or \eqref{eq:fundamentalrelationCase2} (in case $\beta=1$) with any other function in the class $\mathscr{W}_{0}$ different from $-1+t\tanh t$ (see \eqref{W0}), then we get no contradiction, because both sides of the resulting equation vanish. \smallskip

$(ii)$ Let us consider the case $\bta>1$. If we repeat the argument developed in the previous proof, we find that each $\widetilde{\varphi}_{k}$, defined by \eqref{phi-tilde-k}, solves
\begin{equation}
\label{beta-maggiore-di-uno}
\ek g_{\ak}\sin\theta + \ek^{2}\sigma_k+\ek^\bta \mathcal{L}_{\ak}\widetilde{\varphi}_{k} + \ek^{1+\bta} \mathcal{L}_{\ak}^{(1)}\widetilde{\varphi}_{k}+ \ek^{2+\bta} \mathcal{L}_{\ek,\ak}^{(2)}\widetilde{\varphi}_{k} +\ek^{2\bta}h_k  = 2x_{\ak}^2 \ek^\bta A(\widehat{X}_{\ek,\ak}) +  \ek^{\bta+\widetilde\nu} \xi_k\,.
\end{equation}
Moreover, since the sequence $(\widetilde{\varphi}_{k})_{k}$ turns out to be bounded in $C^{2}(\R^{2})$, and $\bta>1$, dividing \eqref{beta-maggiore-di-uno}
 by $\ek$ and taking into account  \eqref{eq:limitgn}, we arrive to find $g_{0}(t)\,\sin\theta\equiv 0$, as a limit equation. This is obviously a contradiction. However this argument suggests that, in case $\bta>1$, solutions of the form $X_{n,a}+\varphi N_{n,a}$ might exist by taking perturbations $\varphi$ whose weighted norm \eqref{eq:weighted-C2-norm} satisfy a bound of order $n^{-1}$ (i.e., $\varphi\in\mathscr{B}_{\eps,a}(1,R)$), rather than $n^{-\beta}$ (i.e., $\varphi\in\mathscr{B}_{\eps,a}(\bta,R)$). In fact, also this case is excluded, according to Theorem \ref{mainteo2}, provided that $R$ is not too large. 
\end{Remark}

\noindent
{\it Proof of Theorem \ref{mainteo2}.} We argue by contradiction and we repeat the same steps as in the proof of Theorem \ref{mainteo}. In the first part, the only difference concerns the normalization condition, which now is accomplished by setting
\[
\widetilde{\varphi}_{k}(t,\theta):=\frac{1}{\ek}\phik(t,\theta), \ \ \ (t,\theta)\in \R^2\,.
\]
In this way we obtain a new sequence $(\widetilde{\varphi}_{k})_k\subset C^{2,\alpha}(\R^2)$ of doubly periodic functions (with respect to the rectangle $K_{n_{k},a_{k}}$), satisfying for all $k$ large enough 
\[
|\widetilde{\varphi}_{k}(t,\theta)|+|\nabla\widetilde{\varphi}_{k}(t,\theta)|+|D^2\widetilde{\varphi}_{k}(t,\theta)|\le \frac{R_{1}h_{\ak}}{\pi} x_{\ak}(t)
\quad\forall(t,\theta)\in\R^2\,.
\]
Moreover each $\widetilde{\varphi}_{k}$ solves the equation
\[
\ek g_{\ak}\sin\theta + \ek \mathcal{L}_{\ak}\widetilde{\varphi}_{k} + \ek^{2} \zeta_k=2x_{\ak}^2 \ek^\bta A(\widehat{X}_{\ek,\ak}) +  \ek^{\bta+\widetilde\nu} \xi_k
\]
where $\zeta_k$ and $\xi_k$ are doubly periodic functions uniformly bounded on $\R^{2}$ with a constant independent of $k$. Passing to the limit $k\to\infty$, we find a function $\widetilde\varphi\in C^{2}(\R\times\R/_{2\pi})$ such that
\begin{equation}
\label{equation-tildephi-beta>1}
\left[4(\sech t)^3-5(\sech t)^5\right] \sin\theta +\Delta\widetilde\varphi+2(\sech t)^{2}\widetilde\varphi=0\quad\text{on~~}\R^{2}
\end{equation}
and
\begin{equation}
\label{bound-tildephi-beta>1}
|\widetilde{\varphi}(t,\theta)|+|\nabla\widetilde{\varphi}(t,\theta)|+|D^2\widetilde{\varphi}(t,\theta)|\le \frac{R_{1}}{\pi} \sech t
\quad\forall(t,\theta)\in\R^2\,,
\end{equation}
because $h_{a_{k}}\to 1$. Evaluating \eqref{equation-tildephi-beta>1}--\eqref{bound-tildephi-beta>1} at $(t,\theta)=(0,\frac{\pi}{2})$, and using the estimate $|\Delta\widetilde\varphi|\le\sqrt{2}\,|D^{2}\widetilde\varphi|$, we obtain
\[
1=|\Delta\widetilde\varphi(0,\tfrac{\pi}{2})+2\widetilde\varphi(0,\tfrac{\pi}{2})|
\le\sqrt{2}\,\frac{R_{1}}{\pi}+2\,\frac{R_{1}}{\pi}
\]
contrary to the bound $R_{1}<\frac{\pi}{2+\sqrt 2}$. 
\qed

\begin{Remark}\label{R:stime-uniformi}
In \cite{CaldiroliMusso} a delicate step is achieving uniform estimates for the linearized operator $\mathcal{L}_{a}$ in some space of periodic functions (see \cite[Theorem 4.6]{CaldiroliMusso}). In view of the previous proofs, we believe that such estimates, which are essential in the argument developed in that paper, in fact cannot be obtained because of the presence of the function $w_{a,0}^{+}$. Indeed, this function does not belong to the kernel of $\mathcal{L}_{a}$, because it is not periodic, but its limit as $a\to 0$ cannot be neglected and it represents the true obstacle in getting certain estimates and actually forbids them.
\end{Remark}
 
\appendix
\renewcommand{\Xe}{X}
\renewcommand{\Xet}{X_t}
\renewcommand{\Xett}{X_{tt}}
\renewcommand{\Xeth}{X_\theta}
\renewcommand{\Xetth}{X_{t\theta}}
\renewcommand{\Xethth}{X_{\theta\theta}}
\renewcommand{\Ne}{N}
\renewcommand{\Net}{N_t}
\renewcommand{\Nett}{N_{tt}}
\renewcommand{\Neth}{N_\theta}
\renewcommand{\Netth}{N_{t\theta}}
\renewcommand{\Nethth}{N_{\theta\theta}}
\section{Basic identities for orthogonal parameterizations}
In this appendix we write some identities for orthogonal parameterizations which are  repeatedly used all through the paper. Firstly, let us recall some elementary vectorial identities. For any $V_1,V_2, V_3 \in \R^3$ one has
\begin{gather}
V_1\wedge (V_2\wedge V_3)=V_2(V_1\cdot V_3)-V_3(V_1\cdot V_2)\,,\label{eq1:propcrossscal}\\[6pt]
V_1\cdot (V_2\wedge V_3)=V_2\cdot(V_3\wedge V_1) =V_3\cdot(V_1\wedge V_2)\,.\label{eq2:propcrossscal}
\end{gather}
Moreover
\begin{equation}\label{eq:propRsigma}
R_{\sigma} V_1\wedge R_{\sigma}V_2= R_\sigma (V_1\wedge V_2),\ \ R_\sigma V_1\cdot R_\sigma V_2=V_1\cdot V_2,\ \text{for any $\sigma\in\R$, $V_1,V_2 \in \R^3$}
\end{equation}
where $R_{\sigma}$ denotes the rotation matrix defined in \eqref{eq:rotation-matrix}. In addition, 
\begin{equation}\label{eq:proppartialRsigma}
Q_\sigma V_1\cdot R_\sigma V_2=\textbf{e}_1\cdot (V_1\wedge V_2) ,\ \ Q_\sigma V_1 \cdot Q_\sigma V_2=\widecheck V_1\cdot \widecheck V_2,\ \text{for any $\sigma\in\R$, $V_1,V_2 \in \R^3$},
\end{equation}
where $\widecheck V \in \R^3$ denotes the vector whose first component is zero and whose remaining two components coincide with those of $V$, that is $\widecheck V:=V-(V\cdot\mathbf{e}_{1})\mathbf{e}_{1}$, and $Q_{\sigma}$ is the matrix obtained by differentiating $R_\sigma$ with respect to $\sigma$ (see \eqref{eq:rotation-matrix} for its explicit form). 

Let  $\Omega\subset\R^2$ be a domain and let $X\colon\Omega\to\R^3$ be a parametrization of a regular surface $\Sigma$. We denote $(t,\theta)$ the pair of parameters in $\Omega$, with the same notation used in the rest of the paper. Assume that $X$ is an orthogonal parametrization, i.e. 
\begin{equation}
\label{eq:XtXth=0}
\Xet\cdot\Xeth=0\,.
\end{equation} 
Then, the unit normal is given by
\begin{equation}\label{eq:defnormal}
\Ne=\frac{\Xet\wedge\Xeth}{|\Xet||\Xeth|}.
\end{equation}
Differentiating the identity $\Ne\cdot \Ne=1$ one has
\begin{equation}\label{eq:prop1NtNtheta}
\Net \cdot \Ne =0,\ \ \Neth \cdot \Ne =0\,.
\end{equation}
Moreover, differentiating \eqref{eq:prop1NtNtheta} and the relations $\Xet \cdot \Ne=0$, $\Xeth \cdot \Ne=0$ one gets 
\begin{equation}\label{identitN}
\begin{array}{ll}
\Nett\cdot \Ne=-|\Net|^2\,
&
\Xet \cdot \Neth=-\Xe_{t\theta}\cdot \Ne=\Xeth \cdot \Net\,,
\\[6pt]
\Nethth \cdot \Ne= -|\Neth|^2\,,
&\Xet \cdot \Net =-\Xett \cdot \Ne,
\\[6pt]
\Ne_{t\theta}\cdot \Ne=-\Net\cdot\Neth\,,\qquad
&\Xeth\cdot \Neth=-\Xethth\cdot \Ne.
\end{array}
\end{equation}
Concerning the second partial derivatives of $\Xe$, by direct computation it is easy to check that
\begin{equation}\label{idDerSecX1}
\begin{array}{ll}
\displaystyle \Xett \cdot \Xet=\displaystyle \frac{1}{2}\frac{d |\Xet|^2}{dt}\,,\qquad
& \displaystyle  \Xe_{t\theta}\cdot \Xet =\displaystyle \frac{1}{2}\frac{d |\Xet|^2}{d\theta},\\[6pt]
  \displaystyle   \Xe_{t\theta}\cdot \Xeth=\displaystyle \frac{1}{2}\frac{d |\Xeth|^2}{dt}\,,\quad
  &\displaystyle  \Xethth\cdot \Xeth \displaystyle \frac{1}{2}\frac{d |\Xeth|^2}{d\theta}\,.
\end{array}
\end{equation}
Moreover, since
\[
\Xett \cdot \Xeth = \frac{d}{dt} \left[\Xet\cdot \Xeth\right] - \Xe_{t\theta}\cdot \Xet\,,\qquad\Xethth \cdot \Xet = \frac{d}{d\theta} \left[\Xet\cdot \Xeth\right] - \Xe_{t\theta}\cdot \Xeth,,
\]
then, exploiting \eqref{eq:XtXth=0} and \eqref{idDerSecX1} we get that
\begin{equation}\label{idDerSecX2}
\Xett \cdot \Xeth=- \frac{1}{2}\frac{d |\Xet|^2}{d\theta}\,,\qquad\Xethth \cdot \Xet =- \frac{1}{2}\frac{d |\Xeth|^2}{dt}.
\end{equation}
From \eqref{eq:defnormal} and exploiting \eqref{eq1:propcrossscal} one gets
\begin{equation}\label{eq:propXetwedgeXeth}
\begin{array}{l}
\displaystyle \Xet \wedge \Ne=\displaystyle  -\frac{|\Xet|}{|\Xeth|} \Xeth= -{|\Xet|} \widehat\Xeth\quad\text{where}\quad \widehat\Xeth=\frac{\Xeth}{|\Xeth|}\,,\\[12pt]
\displaystyle  \Ne \wedge \Xeth=\displaystyle  -\frac{|\Xeth|}{|\Xet|} \Xet=-{|\Xeth|} \widehat\Xet\quad\text{where}\quad \widehat\Xet=\frac{\Xet}{|\Xet|}\,.
\end{array}
\end{equation}
Then, decomposing vectors with respect to the basis $\{\widehat\Xet,\widehat\Xeth,\Ne\}$ and taking into account \eqref{eq:XtXth=0}--\eqref{identitN}, we deduce that
\begin{equation}\label{eq:propNetwedgeNeth}
\begin{array}{ll}
\displaystyle \Xet \wedge \Neth =\displaystyle  -\frac{|\Xet| }{|\Xeth|} (\Xethth \cdot \Ne) \Ne\,,\quad
&\displaystyle \Ne \wedge \Neth =\displaystyle \frac{(\Xethth \cdot \Ne) }{ |\Xeth|}\widehat \Xet - \frac{(\Xetth \cdot \Ne) }{|\Xet| } \widehat\Xeth\,,\\[6pt]
\displaystyle \Net \wedge \Xeth =\displaystyle  -\frac{|\Xeth| }{|\Xet|} (\Xett \cdot \Ne) \Ne\,,&
\displaystyle \Net \wedge \Ne =\displaystyle - \frac{(\Xetth \cdot \Ne) }{|\Xeth| } \widehat\Xet + \frac{(\Xett \cdot \Ne) }{ |\Xet|}\widehat \Xeth\,.
\end{array}
\end{equation}
Finally, using again \eqref{eq:prop1NtNtheta}, \eqref{identitN}, \eqref{eq2:propcrossscal} and also \eqref{eq:propNetwedgeNeth} and decomposing with respect to the basis $\{\widehat\Xet,\widehat\Xeth,\Ne\}$ we obtain
\begin{equation}
\label{eq:propNetwedgeNeth-2}
\displaystyle\Net \wedge \Neth =\displaystyle \frac{(\Xett \cdot \Ne) (\Xethth \cdot \Ne)-(\Xe_{t \theta} \cdot \Ne)^2}{|\Xet| |\Xeth|}  \Ne\,.
\end{equation}


\section{The first and second variations of the mean curvature operator}
In this section we compute the first and second variation of the mean curvature operator with respect to normal variations of orthogonal parametrizations of smooth regular surfaces in $\R^3$.


\begin{Proposition}
\label{Prop:firstvariationmeancurv}
Let  $\Omega\subset\R^2$ be a domain and let $X=X(t,\theta)\colon\Omega\to\R^3$ be an orthogonal parametrization of a smooth regular surface with unit normal $N$. For any given $\varphi \in C^2(\Omega)$, and for any fixed $(t,\theta)\in \Omega$ one has that
\begin{equation}
\label{eq:formulaProp:firstvariationmeancurv}
\begin{split}
&\displaystyle \frac{d}{ds}[\mathfrak{M}(\Xe+s\varphi \Ne)]_{s=0}=\displaystyle  \left(\frac{1}{2|\Xe_t|^2}\right) \varphi_{tt} + \left(\frac{1}{2|\Xe_\theta|^2}\right) \varphi_{\theta\theta} \\[6pt] 
&\qquad\displaystyle + \left(\frac{1}{4|\Xe_t|^2|\Xeth|^2} \frac{d|\Xeth|^2}{dt}- \frac{1}{4|\Xet|^4} \frac{d|\Xet|^2}{dt}\right) \varphi_{t}+  \left(\frac{1}{4|\Xet|^2|\Xeth|^2} \frac{d|\Xet|^2}{d\theta}- \frac{1}{4|\Xeth|^4} \frac{d|\Xeth|^2}{d\theta}\right) \varphi_{\theta}\\[6pt]
&\qquad\displaystyle +  \left(\frac{2(\Xetth\cdot \Ne)^2}{|\Xet|^2|\Xeth|^2} - \frac{|\Neth|^2}{2|\Xeth|^2} - \frac{|\Net|^2}{2|\Xet|^2} + \frac{(\Xett\cdot \Ne)^2}{|\Xet|^4} + \frac{(\Xethth\cdot \Ne)^2}{|\Xeth|^4}  \right) \varphi\,.
\end{split}
\end{equation}
\end{Proposition}

\begin{Proposition}\label{Prop:secvariationmeancurv}
Under the same assumptions of Proposition \ref{Prop:firstvariationmeancurv} one has
\begin{equation}
\label{eq:formulaProp:secvariationmeancurv}
\begin{split}
&\displaystyle \frac{d^2}{ds^2}\left[\mathfrak{M}(\Xe+s\varphi \Ne)\right]_{s=0}= \displaystyle  \varphi_t^2 \left(\frac{1}{2} \frac{\Xett \cdot \Ne}{|\Xet|^4}-\frac{1}{2}\frac{\Xethth\cdot \Ne}{|\Xet|^2 |\Xeth|^2}  \right)  +\varphi_\theta^2 \left(-\frac{1}{2}\frac{\Xett \cdot \Ne}{|\Xet|^2 |\Xeth|^2}+ \frac{1}{2}\frac{\Xethth \cdot \Ne}{|\Xeth|^4}\right) \\[6pt]
&\quad \displaystyle  +\varphi^2 \left( -6\frac{(\Net\cdot \Neth)(\Xetth\cdot \Ne)}{|\Xet|^2 |\Xeth|^2}
+12 \frac{(\Xett\cdot \Ne)(\Xetth\cdot \Ne)^2}{|\Xet|^4 |\Xeth|^2}\right.\\[6pt]
&\quad \displaystyle \left.\ \ \ \ \ \  \ +12 \frac{(\Xethth\cdot \Ne)(\Xetth\cdot \Ne)^2}{|\Xet|^2 |\Xeth|^4}  -3 \frac{(\Xett \cdot \Ne)|\Net|^2}{|\Xet|^4} -3\frac{(\Xethth \cdot \Ne)|\Neth|^2}{ |\Xeth|^4}\right.\\[6pt]
&\quad  \displaystyle \left.\ \ \ \ \ \  \ +4 \frac{(\Xett \cdot \Ne)^3}{|\Xet|^6} +4 \frac{(\Xethth \cdot \Ne)^3}{|\Xeth|^6}   \right) +\varphi_t\varphi_\theta \left(2\frac{\Xetth\cdot \Ne}{|\Xet|^2 |\Xeth|^2} \right)\\[6pt]
&\quad  \displaystyle  +\varphi\varphi_t \left(- \frac{\Nethth\cdot \Xet}{|\Xet|^2 |\Xeth|^2} - \frac{\Nett\cdot \Xet}{|\Xet|^4}  -\frac{3}{2} \frac{\Xett\cdot \Ne}{|\Xet|^6}\frac{d |\Xet|^2}{dt}+\frac{\Xethth\cdot \Ne}{|\Xet|^2|\Xeth|^4}\frac{d |\Xeth|^2}{dt}  +\frac{1}{2} \frac{\Xett\cdot \Ne}{|\Xet|^4 |\Xeth|^2}\frac{d|\Xeth|^2}{dt} \right.\\[6pt]
&\quad \displaystyle \ \ \ \ \ \ \ \ \left. - \frac{1}{2} \frac{(\Xetth\cdot N)}{|\Xet|^2|\Xeth|^4} \frac{d|\Xeth|^2}{d\theta}- \frac{3}{2} \frac{\Xetth\cdot \Ne}{|\Xet|^4 |\Xeth|^2} \frac{d|\Xet|^2}{d\theta} \right)\\[6pt]
&\quad  \displaystyle  +\varphi\varphi_\theta\left(- \frac{\Nett\cdot \Xeth}{|\Xet|^2 |\Xeth|^2} - \frac{\Nethth\cdot \Xeth}{|\Xeth|^4}  -\frac{3}{2} \frac{\Xethth\cdot \Ne}{|\Xeth|^6}\frac{d |\Xeth|^2}{d\theta} +\frac{\Xett\cdot \Ne}{|\Xet|^4|\Xeth|^2}\frac{d |\Xet|^2}{d\theta} +\frac{1}{2} \frac{\Xethth\cdot \Ne}{|\Xet|^2 |\Xeth|^4}\frac{d|\Xet|^2}{d\theta}\right.\\[6pt]
&\quad  \displaystyle  \ \ \ \ \ \ \ \ \left.\ \ - \frac{1}{2} \frac{\Xetth\cdot \Ne}{|\Xet|^4|\Xeth|^2} \frac{d|\Xet|^2}{dt}- \frac{3}{2} \frac{\Xetth\cdot \Ne}{|\Xet|^2 |\Xeth|^4} \frac{d|\Xeth|^2}{dt} \right)\\[6pt]
&\quad  \displaystyle  +\varphi\varphi_{\theta\theta} \left(2\frac{\Xethth\cdot \Ne}{|\Xeth|^4} \right) +\varphi\varphi_{tt} \left(2\frac{\Xett\cdot \Ne}{|\Xet|^4}\right) + \varphi\varphi_{t\theta} \left( 4\frac{\Xetth\cdot \Ne}{|\Xet^2 |\Xeth|^2}\right)\,.
\end{split}
\end{equation}

\end{Proposition}
For the proofs of Proposition \ref{Prop:firstvariationmeancurv}, Proposition \ref{Prop:secvariationmeancurv} we need a preliminary technical result.

\begin{Lemma}\label{lem:techNphi}
Under the same assumptions of Proposition \ref{Prop:firstvariationmeancurv}, letting
\begin{equation}
\label{N-varphi-s}
N_{\varphi}(s)=\frac{(X+s\varphi{\normal})_{t}\wedge (X+s\varphi{\normal})_{\theta}}{|(X+s\varphi{\normal})_{t}\wedge (X+s\varphi{\normal})_{\theta}|}\,,
\end{equation}
one has that
\begin{equation}\label{eq5foglio}
\left[\displaystyle \frac{d}{ds}N_{\varphi}(s)\right]_{s=0}=\displaystyle -\varphi_\theta \frac{\widehat \Xeth}{|\Xeth|}  - \varphi_t \frac{ \widehat \Xet}{|\Xet|}
\end{equation}
and
\begin{equation}\label{eq6foglio}
\begin{split}
\left[\displaystyle \frac{d^{2}}{ds^{2}}N_{\varphi}(s)\right]_{s=0}&=\displaystyle -2\left(\varphi_t \varphi  \frac{\Xett\cdot \Ne}{|\Xet|^2} +  \varphi_\theta \varphi  \frac{\Xetth\cdot \Ne}{|\Xeth|^2} \right)\frac{\widehat \Xet}{|\Xet|} \\[6pt]
&\quad  -2 \left(\varphi_\theta \varphi  \frac{\Xethth\cdot \Ne}{|\Xeth|^2} + \varphi_t \varphi  \frac{\Xetth\cdot \Ne}{|\Xet|^2} \right)\frac{\widehat \Xeth}{|\Xeth|} \displaystyle -\left(\frac{ \varphi_t^2}{|\Xet|^2} + \frac{\varphi_\theta^2}{|\Xeth|^2}\right)\Ne\,.
\end{split}
\end{equation}
\end{Lemma}

\Proof
Let $X$, $N$ and $\varphi$ be as in the statement of Proposition \ref{Prop:firstvariationmeancurv}. Set
\[
\begin{split}
J_{\varphi}(s):&\!\!=(X+s\varphi{\normal})_{t}\wedge (X+s\varphi{\normal})_{\theta}\\[3pt]
&\!\!=
\Xet \wedge \Xeth+ s \left[\varphi_\theta (\Xet \wedge \Ne) +  \varphi (\Xet \wedge \Neth) +\varphi_t (\Ne\wedge \Xeth) +\varphi (\Net \wedge \Xeth)\right]\\[3pt]
&\quad+  s^2 \left[\varphi_t \varphi (\Ne\wedge \Neth)+ \varphi \varphi_\theta (\Net\wedge \Ne)+ \varphi^2 (\Net\wedge \Neth)\right]\,.
\end{split}
\]
Using \eqref{eq:defnormal} and \eqref{eq:propXetwedgeXeth}, \eqref{eq:propNetwedgeNeth} we have that
\begin{equation}\label{eq1foglio}
\begin{split}
J_{\varphi}(s)&=\displaystyle |\Xet| |\Xeth| \Ne + s \left[-\varphi_\theta |\Xet| \widehat \Xeth - \varphi_t |\Xeth| \widehat \Xet -  \varphi \left(\frac{|\Xet|}{|\Xeth|} (\Xethth \cdot \Ne) + \frac{|\Xeth|}{|\Xet|} (\Xett \cdot \Ne)\right) \Ne \right]\\[6pt]
&\quad+  s^2 \left[\varphi_t \varphi \left(\frac{\Xethth\cdot \Ne} {|\Xeth|} \widehat \Xet -\frac{\Xe_{t\theta}\cdot \Ne}{|\Xet|} \widehat \Xeth \right) + \varphi \varphi_\theta  \left(-\frac{\Xe_{t\theta}\cdot \Ne}{|\Xeth|} \widehat \Xet +\frac{\Xett\cdot \Ne}{|\Xet|} \widehat \Xeth \right)\right.\\[6pt]
&\qquad\displaystyle \left.\quad \quad + \varphi^2  \frac{(\Xe_{tt} \cdot \Ne) (\Xe_{\theta\theta} \cdot \Ne)-(\Xe_{t \theta} \cdot \Ne)^2}{|\Xet|| \Xeth|} \Ne \right].
\end{split}
\end{equation}
Setting for brevity
\begin{equation}\label{eq:defAB}
\displaystyle f_{1}:= \displaystyle \frac{|\Xet|}{|\Xeth|} (\Xethth\cdot \Ne) + \frac{|\Xeth|}{|\Xet|} (\Xett \cdot \Ne)\quad\text{and}\quad
\displaystyle f_{2}:=\displaystyle \frac{(\Xett \cdot \Ne) (\Xethth \cdot \Ne)-(\Xetth \cdot \Ne)^2}{|\Xet||\Xeth|}\,,
\end{equation} 
from \eqref{eq1foglio} we deduce that
\begin{equation}\label{eq2foglio}
\begin{split}
&\displaystyle \frac{d}{ds}J_{\varphi}(s)=\displaystyle -\varphi_\theta |\Xet| \widehat \Xeth  - \varphi_t |\Xeth| \widehat \Xet -  \varphi f_{1} \Ne\\[6pt]
&\quad\displaystyle   + 2s \left[\left(\varphi_t \varphi \frac{\Xethth\cdot \Ne}{|\Xeth|} - \varphi \varphi_\theta \frac{\Xetth\cdot \Ne}{|\Xeth|}\right)\widehat\Xet+ \left(\varphi \varphi_\theta  \frac{\Xett\cdot \Ne}{|\Xet|} - \varphi_t \varphi  \frac{\Xetth\cdot \Ne}{|\Xet|} \right)\widehat\Xeth+ \varphi^2 f_{2} \Ne \right].
\end{split} 
\end{equation}
Now, we observe that
\[
\frac{d}{ds}|J_{\varphi}|=\frac{J_{\varphi}\cdot\frac{dJ_{\varphi}}{ds}}{|J_{\varphi}|}\quad\text{and}\quad
\frac{d^{2}}{ds^{2}}|J_{\varphi}|=\frac{\left|\frac{dJ_{\varphi}}{ds}\right|^{2}+{J_{\varphi} \cdot \frac{d^2J_{\varphi}}{ds^2} }-\left(\frac{d}{ds}|J_{\varphi}|\right)^{2}}{|J_{\varphi}|}\,.
\]
Using \eqref{eq2foglio}  we see that
\begin{equation}\label{eq:gprime0}
\left[\frac{d}{ds}|J_{\varphi}(s)|\right]_{s=0}= \Ne \cdot(- \varphi f_{1} \Ne)=-\varphi f_{1}\,.
\end{equation}
Moreover, observing that
\[
\left|\frac{d}{ds}J_{\varphi}(s)\right|^2_{s=0}= \varphi_\theta^2 |\Xet|^2+ \varphi_t^2 |\Xeth|^2+ \varphi^2 f_{1}^2\,,
\]
and since
\[
\Ne\cdot\left[\frac{d^2}{ds^2}J_{\varphi}(s)\right]_{s=0}=2 \varphi^2 f_{2}\,,
\]
as it follows by differentiating \eqref{eq2foglio}, then, after some elementary computations, we obtain
\begin{equation}\label{eq4foglio}
\left[\frac{d^{2}}{ds^{2}}\left|J_{\varphi}(s)\right|\right]_{s=0}=\displaystyle  \varphi_\theta^2 \frac{|\Xet|}{|\Xeth|}+ \varphi_t^2\frac{|\Xeth|}{|\Xet|}+  2 \varphi^2f_{2}\,.
\end{equation}
Finally, since $N_{\varphi}(s)=|J_{\varphi}(s)|^{-1}J_{\varphi}(s)$,
\[
\frac{d}{ds}N_{\varphi}(s)= \frac{\frac{d}{ds}J_{\varphi}(s)}{\left| J_{\varphi}(s)\right|} -\frac{J_{\varphi}(s)}{\left| J_{\varphi}(s)\right|^2} \frac{d}{ds}|J_{\varphi}(s)|\,,
\]
then, exploiting \eqref{eq2foglio}--\eqref{eq4foglio} we get \eqref{eq5foglio}. Moreover, since
\[
\frac{d^2}{ds^2}N_{\varphi}(s)=\frac{\frac{d^2}{ds^2}J_{\varphi}(s)}{\left| J_{\varphi}(s) \right|} -2 \frac{\frac{d}{ds}J_{\varphi}(s)}{\left| J_{\varphi}(s) \right|^2} \frac{d}{ds}|J_{\varphi}(s)|-\frac{J_{\varphi}(s)}{\left| J_{\varphi}(s) \right|^2} \frac{d^{2}}{ds^{2}}|J_{\varphi}(s)| +2 \frac{J_{\varphi}(s)}{\left| J_{\varphi}(s) \right|^3} \left(\frac{d}{ds}|J_{\varphi}(s)|\right)^2\,,
\]
then, recalling the definition of $f_{1}$ (see \eqref{eq:defAB}), using \eqref{eq:gprime0}, \eqref{eq4foglio}, and after some elementary computations, we obtain \eqref{eq6foglio}.
\qed

\textit{Proof of Propositions \ref{Prop:firstvariationmeancurv} and \ref{Prop:secvariationmeancurv}.}
Let $X$, $N$ and $\varphi$ be as in the statement of Proposition \ref{Prop:firstvariationmeancurv}. By definition of mean curvature one has
\[
\mathfrak{M}(\Xe+s\varphi \Ne)=\frac{{\E}(s){\N}(s)-2{\F}(s){\M}(s)+{\G}(s){\L}(s)}{2({\E}(s){\G}(s)-{\F}^{2}(s))}=:\frac{{P}(s)}{{Q}(s)}
\]
where the coefficients of the first and second fundamental form correspond to the parameterization $\Xe+s\varphi \Ne$. Then
\begin{equation}\label{eq:espfirstvar}
\left[\frac{d}{ds}\mathfrak{M}(\Xe+s\varphi \Ne)\right]_{s=0}=\frac{P'(0)Q(0)-P(0)Q'(0)}{(Q(0))^{2}}=\frac{P'(0)}{Q(0)}-\mathfrak{M}(X) \frac{Q'(0)}{Q(0)}\,,
\end{equation}
and 
\begin{equation}\label{eq:varsecmeancurv}
\left[\displaystyle\frac{d^2}{ds^2}\mathfrak{M}(\Xe+s\varphi \Ne)\right]_{s=0}
=\displaystyle\frac{P''(0)}{Q(0)} -2\frac{P'(0)Q'(0)}{(Q(0))^{2}}+2 \ \mathfrak{M}(X) \left(\frac{Q'(0)}{Q(0)}\right)^2 -\mathfrak{M}(X) \frac{Q''(0)}{Q(0)}\,.
\end{equation}
To conclude we need to compute $P(0)$ $P'(0)$, $P''(0)$, $Q(0)$, $Q'(0)$, $Q''(0)$. To this end we begin observing that from \eqref{eq:defnormal} and \eqref{eq:prop1NtNtheta} one has
\[
{\E}(s)=|(\Xe+s\varphi{\normal})_{t}|^{2}= |\Xet|^2 + 2 s \varphi (\Xet\cdot \Net)+ s^2 \left[ \varphi_t^2 + \varphi^2 |\Net|^2\right]\,,
\]
and thus
\begin{equation}\label{eq:coeffE}
{\E}(0)=|\Xet|^2\,,\quad {\E}^\prime(0)=2 \varphi \Xet\cdot \Net\,,\quad{\E}^{\prime\prime}(0)= 2\varphi_t^2 + 2\varphi^2 |\Net|^2\,.
\end{equation}
By the same argument, one has
\[
{\G}(s)=|(\Xe+s\varphi\Ne)_{\theta}|^{2}=|\Xeth|^2 + 2 s \varphi (\Xeth\cdot \Neth)+ s^2 \left[ \varphi_\theta^2 + \varphi^2 |\Neth|^2\right]\,,
\]
so that
\begin{equation}\label{eq:coeffG}
{\G}(0)=  |\Xeth|^2,\ \ {\G}^\prime(0)= 2 \varphi \Xeth\cdot \Neth\,,\quad{\G}^{\prime\prime}(0)=  2 \varphi_\theta^2 + 2\varphi^2 |\Neth|^2\,.
\end{equation}
Moreover, since
\[
{\F}(s) =\Xet \cdot \Xeth+  s \varphi [(\Xet \cdot \Neth)+(\Xeth\cdot\Net)]+ s^2 \left[\varphi_{t}  \varphi_{\theta} +\varphi^2 (\Net\cdot \Neth) \right]\,,
\]
and $X$ is orthogonal, by \eqref{identitN} it is immediate to check that 
\begin{equation}\label{eq:coeffF}
{\F}(0)=0\,,\quad{\F}^\prime(0)=-2\varphi (X_{t\theta} \cdot N) \,,\quad{\F}^{\prime\prime}(0)=  2\varphi_{t}  \varphi_{\theta} + 2\varphi^2 (\Net\cdot \Neth)\,.
\end{equation}
Let us treat the remaining coefficients. Firstly, by definition one has
\begin{equation}\label{eq7fogliobis}
\N(s)=\left[\Xethth + s \left( \varphi_{\theta\theta} \Ne+ 2 \varphi_\theta \Neth+\varphi  \Nethth\right)\right] \cdot N_{\varphi}(s)
\end{equation}
where ${\normal}_{\varphi}(s)$ is the normal vector corresponding to the parametric surface $X+s\varphi N$, defined in \eqref{N-varphi-s}. In particular, recalling that $N_{\varphi}(0)=N$, we have that
\begin{equation}\label{eq:coeffN}
\N(0)= \Xethth\cdot \Ne\,,
\end{equation}
and differentiating \eqref{eq7fogliobis}, taking into account Lemma \ref{lem:techNphi}, we infer that
\begin{equation}\label{eq:coeffNprime}
\N^\prime(0)=\displaystyle  \varphi_{\theta\theta} + \varphi (\Nethth\cdot \Ne) - \varphi_\theta \frac{\Xethth \cdot \Xeth}{|\Xeth|^2}  - \varphi_t \frac{\Xethth \cdot \Xet}{|\Xet|^2}
\end{equation}
and 
\begin{equation}\label{eq:coeffNsec}
\begin{split}
\N^{\prime\prime}(0)
&= \displaystyle -4\varphi_\theta^2  \frac{\Neth\cdot \Xeth}{|\Xeth|^2}  -4\varphi_\theta \varphi_t  \frac{\Neth\cdot  \Xet}{|\Xet|^2}  -2\varphi_\theta \varphi  \frac{\Nethth\cdot \Xeth}{|\Xeth|^2} \\[6pt]
&\quad \displaystyle -2\varphi_t \varphi  \frac{\Nethth\cdot \Xet}{|\Xet|^2} -2\varphi_t \varphi  \frac{(\Xethth\cdot \Xet)(\Xett\cdot \Ne)}{|\Xet|^4}\\[6pt]
 &\quad \displaystyle -2\varphi_\theta \varphi  \frac{(\Xethth\cdot \Xet)(\Xetth\cdot \Ne)}{|\Xet|^2|\Xeth|^2}-2\varphi_\theta \varphi  \frac{(\Xethth\cdot\Xeth)(\Xethth\cdot \Ne)}{|\Xeth|^4}\\[6pt]
&\quad \displaystyle -2\varphi_t \varphi  \frac{(\Xethth\cdot \Xeth)(\Xetth\cdot \Ne)}{|\Xet|^2|\Xeth|^2}- \varphi_t^2\frac{\Xethth\cdot \Ne}{|\Xet|^2}- \varphi_\theta^2\frac{\Xethth\cdot \Ne}{|\Xeth|^2}\,.
\end{split}
\end{equation}
Similarly, since $\L(s)=\left[\Xett+s\left(\varphi_{tt} \Ne+ 2\varphi_t \Net + \varphi \Nett \right)\right]\cdot N_{\varphi}(s)$, we immediately deduce that
\begin{equation}\label{eq:coeffL}
\L(0)= \Xett\cdot \Ne\,,
\end{equation}
\begin{equation}\label{eq:coeffLprime}
\displaystyle \L^\prime(0) =\displaystyle \varphi_{tt} + \varphi \Nett\cdot \Ne -\varphi_\theta \frac{\Xett\cdot \Xeth}{|\Xeth|^2}  - \varphi_t \frac{\Xett\cdot\Xet}{|\Xet|^2}\,,
\end{equation}
\begin{equation}\label{eq:coeffLsec}
\begin{split}
\displaystyle \L^{\prime\prime}(0)
&=\displaystyle -4\varphi_t^2  \frac{\Net\cdot \Xet}{|\Xet|^2}  -4\varphi_t \varphi_\theta   \frac{\Net\cdot  \Xeth}{|\Xeth|^2}  -2\varphi_\theta \varphi  \frac{\Nett\cdot \Xeth}{|\Xeth|^2} \\[6pt]
&\quad \displaystyle -2\varphi_t \varphi  \frac{\Nett\cdot \Xet}{|\Xet|^2} -2\varphi_t \varphi  \frac{(\Xett\cdot \Xet)(\Xett\cdot \Ne)}{|\Xet|^4}\\[6pt]
 &\quad\displaystyle -2\varphi_\theta \varphi  \frac{(\Xett\cdot \Xet)(\Xetth\cdot \Ne)}{|\Xet|^2|\Xeth|^2}-2\varphi_\theta \varphi  \frac{(\Xett\cdot\Xeth)(\Xethth\cdot \Ne)}{|\Xeth|^4}\\[6pt]
&\quad \displaystyle -2\varphi_t \varphi  \frac{(\Xett\cdot \Xeth)(\Xetth\cdot \Ne)}{|\Xet|^2|\Xeth|^2}- \varphi_t^2\frac{\Xett\cdot \Ne}{|\Xet|^2}- \varphi_\theta^2\frac{\Xett\cdot \Ne}{|\Xeth|^2}\,.
\end{split}
\end{equation}
Finally, we study the term $\M(s)$. By definition we have
\begin{equation}\label{eq:defMphi}
\M(s)= \left[\Xetth + s\left(\varphi_{t\theta} \Ne+\varphi_t \Neth + \varphi_\theta \Net + \varphi \Netth \right)\right] \cdot N_{\varphi}(s),
\end{equation}
and thus 
\begin{equation}\label{eq:coeffM}
\M(0)=\Xetth\cdot \Ne\,.
\end{equation}
Differentiating \eqref{eq:defMphi}, exploiting \eqref{eq:prop1NtNtheta} and using Lemma \ref{lem:techNphi}, then, after some elementary computations we get that 
%
\begin{equation}\label{eq:coeffMprime}
\displaystyle\M^\prime(0)=\displaystyle  \varphi_{t\theta} + \varphi \Netth\cdot \Ne - \varphi_\theta \frac{\Xetth \cdot \Xeth }{|\Xeth|^2} - \varphi_t \frac{\Xetth \cdot \Xet}{|\Xet|^2},
\end{equation}
and
\begin{equation}\label{eq:coeffMsec}
\begin{split}
\displaystyle\M^{\prime\prime}(0)&=\displaystyle  - 2\varphi_t \varphi_\theta\left(\frac{\Neth\cdot\Xeth}{|\Xeth|^2}+\frac{\Net\cdot\Xet}{|\Xet|^2}\right)   - 2\varphi_t^2 \frac{\Neth\cdot  \Xet}{|\Xet|^2}  - 2\varphi_\theta^2 \frac{\Net\cdot \Xeth}{|\Xeth|^2}\\[6pt]
&\quad\displaystyle  -  2\varphi \varphi_\theta \frac{\Netth\cdot \Xeth}{|\Xeth|^2}  -  2\varphi \varphi_t \frac{\Netth\cdot \Xet}{|\Xet|^2} -2 \varphi_t \varphi  \frac{(\Xett\cdot \Ne)(\Xetth\cdot \Xet)}{|\Xet|^4}\\[6pt]
&\quad\displaystyle   -2 \varphi_\theta \varphi  \frac{(\Xetth\cdot \Ne)(\Xetth\cdot \Xet)}{|\Xet|^2|\Xeth|^2}  -2 \varphi_\theta \varphi  \frac{(\Xethth\cdot \Ne) (\Xetth\cdot\Xeth)}{|\Xeth|^4}\\[6pt]
&\quad\displaystyle -2 \varphi_t \varphi  \frac{(\Xetth\cdot \Ne)(\Xetth\cdot \Xeth)}{|\Xet|^2|\Xeth|^2} - \varphi_t^2\frac{(\Xetth\cdot \Ne)}{|\Xet|^2}- \varphi_\theta^2\frac{(\Xetth\cdot \Ne)}{|\Xeth|^2}\,.
\end{split}
\end{equation}
At the end, putting together \eqref{eq:coeffE}--\eqref{eq:coeffF}, \eqref{eq:coeffN}--\eqref{eq:coeffNprime}, \eqref{eq:coeffL}--\eqref{eq:coeffLprime}, \eqref{eq:coeffM}--\eqref{eq:coeffMprime}, exploiting \eqref{identitN}--\eqref{idDerSecX2} and after some elementary compuations, we get
\begin{equation}\begin{split} 
\label{eqPQprime}
\displaystyle {P}(0)&=\displaystyle |\Xet|^2(\Xethth\cdot \Ne)+|\Xeth|^2 (\Xett\cdot \Ne)\,,\\[8pt]
\displaystyle{Q}(0)&=\displaystyle 2|\Xet|^2|\Xeth|^2\,,\\[8pt]
\displaystyle {P}'(0)&= \displaystyle  \varphi \left[-4(\Xett \cdot \Ne) (\Xethth \cdot \Ne) + 4(\Xetth\cdot \Ne)^2  -|\Neth|^2 |\Xet|^2 - |\Net|^2 |\Xeth|^2 \right] + \varphi_{\theta\theta} |\Xet|^2  \\
&\quad+ \varphi_{tt}  |\Xeth|^2\displaystyle   - \varphi_t \left[- \frac{1}{2}\frac{d |\Xeth|^2}{dt}+ \frac{1}{2}\frac{ |\Xeth|^2}{|\Xet|^2}\frac{d |\Xet|^2}{dt}\right] - \varphi_\theta \left[- \frac{1}{2}\frac{d |\Xet|^2}{d\theta} +  \frac{1}{2}\frac{ |\Xet|^2}{|\Xeth|^2}\frac{d |\Xeth|^2}{d\theta} \right]\,,\\
\displaystyle  {Q}'(0)&= \displaystyle \varphi \left[-4 (\Xett\cdot \Ne)|\Xeth|^2 -4 (\Xethth\cdot \Ne)|\Xet|^2\right]\,.
\end{split}
\end{equation}
Finally, from \eqref{eq:espfirstvar} and \eqref{eqPQprime}, taking into account that 
\[
\mathfrak{M}(X)=\frac{P(0)}{Q(0)}=\frac{1}{2}\left(\frac{\Xett\cdot \Ne}{|\Xet|^2}+\frac{\Xethth\cdot \Ne}{|\Xeth|^2}\right),
\]
we easily obtain \eqref{eq:formulaProp:firstvariationmeancurv}, and the proof of Proposition \ref{Prop:firstvariationmeancurv} is complete.
\medskip

Arguing in a similar way, from \eqref{eq:coeffE}--\eqref{eq:coeffF}, \eqref{eq:coeffN}--\eqref{eq:coeffLsec}, \eqref{eq:coeffM}--\eqref{eq:coeffMsec}, 
\eqref{identitN}--\eqref{idDerSecX2} and after a tedious computation one gets
\begin{equation}\label{eqPsecond}
\begin{split}
\displaystyle  P''(0)&= \displaystyle  \varphi_t^2 \left[ \Xethth \cdot \Ne +3 \frac{(\Xett \cdot \Ne) |\Xeth|^2}{|\Xet|^2} \right]+\varphi_\theta^2 \left[\Xett \cdot \Ne)+3 \frac{(\Xethth \cdot \Ne) |\Xet|^2}{|\Xeth|^2}\right]
\\[3pt] 
&\quad \displaystyle  +\varphi^2 \left[ 6(\Xethth \cdot \Ne) |\Net|^2 + 6(\Xett \cdot \Ne) |\Neth|^2-12(\Net\cdot \Neth)(\Xetth\cdot \Ne)\right]  \ +4\varphi_t\varphi_\theta (\Xetth\cdot \Ne)
\\[3pt]
&\quad \displaystyle+\varphi\varphi_t \left[-\frac{(\Xett\cdot \Ne)}{|\Xet|^2}\frac{d|\Xeth|^2}{dt} +2\frac{(\Xethth\cdot \Ne)}{|\Xet|^2}\frac{d|\Xet|^2}{dt} - 2 \Nethth\cdot \Xet - 2 \frac{(\Nett\cdot \Xet)|\Xeth|^2}{|\Xet|^2} \right.
\\[3pt]
&\quad \displaystyle \ \ \ \ \ \ \ \ \left. -  \frac{(\Xett\cdot \Ne)|\Xeth|^2}{|\Xet|^4}\frac{d|\Xet|^2}{dt}-  \frac{(\Xetth\cdot \Ne)}{|\Xeth|^2} \frac{d|\Xeth|^2}{d\theta}- 3 \frac{(\Xetth\cdot \Ne)}{|\Xet|^2} \frac{d|\Xet|^2}{d\theta} \right]
\\[3pt]
&\quad \displaystyle  +\varphi\varphi_\theta \left[-\frac{(\Xethth\cdot \Ne)}{|\Xeth|^2}\frac{d|\Xet|^2}{d\theta} +2\frac{(\Xett\cdot \Ne)}{|\Xeth|^2}\frac{d|\Xeth|^2}{d\theta} - 2 \Nett\cdot \Xeth  - 2 \frac{(\Nethth\cdot \Xeth)|\Xet|^2}{|\Xeth|^2} \right.
\\[3pt] 
&\quad \displaystyle\ \ \ \ \ \ \ \ \ \left.
 -  \frac{(\Xethth\cdot \Ne)|\Xet|^2}{|\Xeth|^4}\frac{d|\Xeth|^2}{d\theta} -  \frac{(\Xetth\cdot \Ne)}{|\Xet|^2} \frac{d|\Xet|^2}{dt}- 3 \frac{(\Xetth\cdot \Ne)}{|\Xeth|^2} \frac{d|\Xeth|^2}{dt} \right]
\\[3pt]
&\quad \displaystyle  -4\varphi\varphi_{\theta\theta}(\Xett\cdot \Ne) - 4\varphi\varphi_{tt} (\Xethth\cdot \Ne) + 8\varphi\varphi_{t\theta} (\Xetth\cdot \Ne)\,,
\end{split}
\end{equation}
\begin{equation}\label{eqQsecond}
\begin{split}
\displaystyle Q''(0)&=\displaystyle 4\varphi_t^2|\Xeth|^2+4\varphi_\theta^2|\Xet|^2
\\
 &\quad \displaystyle +\varphi^2 \left[ 4|\Net|^2 |\Xeth|^2+4|\Neth|^2 |\Xet|^2 + 16(\Xett\cdot \Ne)(\Xethth\cdot \Ne)-16(\Xetth\cdot \Ne)^2\right]\,.
\end{split}
\end{equation}
From \eqref{eq:varsecmeancurv}, \eqref{eqPQprime}, \eqref{eqPsecond}--\eqref{eqQsecond} and after elementary computations we finally obtain \eqref{eq:formulaProp:secvariationmeancurv}.
\qed

\begin{Remark}
\label{R:weight-C1}
The computation of the first derivative of $N_{\varphi}(s)$ stated in Lemma \ref{lem:techNphi} turns out to be useful to justify the weight $x_{a}^{-1}$ in the norm \eqref{eq:weighted-C2-norm} up to the first derivatives. More precisely, fixing $a_{0}\in(0,\frac{1}{2})$, assume that there exist  $\eps_{0}>0$ and, for every $a\in[-a_{0},a_{0}]\setminus\{0\}$, a continuous function $\omega_{a}=\omega_{a}(t)\colon\R\to\R$ such that for any $\varphi$ belonging to some convex neighborhood of $0$  in $C^{1}(\R^{2})$, endowed with its standard norm, one has  
\begin{equation}
\label{eq:corollary}
\left|
\frac{(X_{\eps,a}+\varphi N_{\eps,a})_{t}\wedge(X_{\eps,a}+\varphi N_{\eps,a})_{\theta}}{|(X_{\eps,a}+\varphi N_{\eps,a})_{t}\wedge(X_{\eps,a}+\varphi N_{\eps,a})_{\theta}|}
-N_{\eps,a}\right|\le \omega_{a}|\nabla\varphi|\quad\text{on~~}\R^{2}\,.
\end{equation}
Then there exists a constant $C_{0}>0$ independent of $a\in[-a_{0},a_{0}]\setminus\{0\}$ and $\eps\in(0,\eps_{0})$ such that $\omega_{a}\ge C_{0}x_{a}^{-1}$. Indeed, writing \eqref{eq:corollary} with $s\varphi$ instead of $\varphi$, dividing by $|s|$ and passing to the limit $s\to 0$, by \eqref{eq5foglio} we obtain
\[
\sqrt{
\left(\frac{\varphi_{\theta}}{|(X_{\eps,a})_{\theta}|}\right)^{2}+
\left(\frac{\varphi_{t}}{|(X_{\eps,a})_{t}|}\right)^{2}}\le \omega_{a}|\nabla\varphi|\,.
\]
Then, taking account of \eqref{eq:X_tX_thetaVW}, the conclusion follows.
\end{Remark}

\subsection{An expansion for the second variation of the mean curvature of Delaunay tori}
In this subsection we provide an expansion for the second variation of the mean curvature of Delaunay tori along normal variations.
\begin{Proposition}\label{prop:secvarexpansion}
Let $a\in[-\frac{1}{2},\infty)\setminus\{0\}$, $\eps>0$ small enough, let $X_{\eps,a}$ be the parametrization of a Delaunary torus given by \eqref{eq:eps-toroidal-unduloid} with unit normal $N_{\eps,a}$, and let $\varphi \in C^2(\R^2)$. Then
\[
\begin{split}
\displaystyle \frac{d^2}{ds^2}\left[\mathfrak{M}(X_{\eps,a}+s\varphi N_{\eps,a})\right]_{s=0} &=\displaystyle\varphi_t^2 \left(\frac{\gamma_a}{x_a^4}\right) + \varphi_\theta^2 \left(-\frac{\gamma_a}{x_a^4}\right) + \varphi^2\left[\left(\frac{x_a^2+\gamma_a}{x_a^2}\right)^3+\left(\frac{x_a^2-\gamma_a}{x_a^2}\right)^3\right] \\[6pt]
&\quad\displaystyle+2\varphi\varphi_{tt} \left(\frac{x_a^2+\gamma_a}{x_a^2}\right) + 2\varphi\varphi_{\theta\theta} \left(\frac{x_a^2-\gamma_a}{x_a^2}\right)+ \eps {\mathcal{Q}}_{\eps,a}(\varphi),
\end{split}
\]
where
\[
{\mathcal{Q}}_{\eps,a}(\varphi)=\bta_{\eps,a}^{(1)}\varphi_{t}^2+\bta_{\eps,a}^{(2)}\varphi_{\theta}^2+\bta_{\eps,a}^{(3)}\varphi^2+ \bta_{\eps,a}^{(4)}\varphi_{t}\varphi_\theta +\bta_{\eps,a}^{(5)}\varphi\varphi_t +\bta_{\eps,a}^{(6)}\varphi\varphi_\theta +\bta_{\eps,a}^{(7)}\varphi\varphi_{tt}+\bta_{\eps,a}^{(8)}\varphi\varphi_{\theta\theta}+\bta_{\eps,a}^{(9)}\varphi\varphi_{t\theta}
\]
is a quadratic form with coefficients $\bta_{\eps,a}^{(i)}=\bta_{\eps,a}^{(i)}(t,\theta)$ ($i=1,\ldots,9$) which are doubly periodic with respect to the rectangle $K_a=[-\tau_a, \tau_a]\times[-\pi,\pi]$. Moreover, fixing $a_{0}\in(0,\frac{1}{2})$, there exist $\eps_0>0$ and $C>0$ such that for any $a\in[-a_0,a_0]\setminus\{0\}$ and $\eps\in(0,\eps_0)$ one has
\[
\|x_a^3\bta_{\eps,a}^{(i)}\|_{C^0(K_a)}  \leq C, \ \ \forall i=1,\ldots,9\,.
\]
\end{Proposition}
\Proof
The result is a straightforward consequence of Proposition \ref{Prop:secvariationmeancurv}, Lemmata \ref{lem:expsecderXN}--\ref{lem:expNsecX}, taking into account \eqref{eq:xz}, \eqref{eq:X_tX_thetaVW}, \eqref{eq:expinvsqnormXt} and \eqref{eq:dersquarednormttheta}.
\qed

\noindent
\textbf{Acknowledgements.}
The first and the second author are members of the Gruppo Nazionale per l'Analisi Matematica, la Probabilit\`a e le loro Applicazioni (GNAMPA) of the Istituto Nazionale di Alta Matematica (INdAM). The third author has been supported by EPSRC research GrantEP/T008458/1. This paper has been written when the second author was at D\'epartement de Math\'ematique, Universit\'e Libre de Bruxelles and was also supported by FNRS-F.R.S.

\end{document}